\documentclass[titlepage,11pt]{article}
\usepackage{url}
\usepackage{latexsym}
\usepackage{amsfonts}
\usepackage{amsmath}
\newcommand\blackslug{\hbox{\hskip 1pt \vrule width 4pt height 8pt depth 1.5pt
        \hskip 1pt}}
\newcommand\bbox{\hfill \quad \blackslug \bigbreak}
\newcommand{\indentitem}{\setlength\itemindent{25pt}}
\def\d{\hbox{-}}
\def\c{\hbox{-}\cdots\hbox{-}}
\def\l{,\ldots,}

\title{Induced subgraphs of graphs with large chromatic number. \\
IV. Consecutive holes}
\author{Alex Scott\thanks{Supported by a Leverhulme Trust Research Fellowship}\\
Oxford University, Oxford, UK
\\
\\
Paul Seymour\thanks{Supported by ONR grant N00014-14-1-0084 and NSF
grant DMS-1265563.}\\
Princeton University, Princeton, NJ 08544, USA}

\date{January 17, 2015; revised \today}

\newtheorem{thm}{}[section]

\newcommand{\Proof}{\noindent{\bf Proof.}\ \ }

\begin{document}
\maketitle
\begin{abstract}
A {\em hole} in a graph is an induced subgraph which is a cycle of length at 
least four. 
We prove that for all $\nu>0$, every triangle-free graph with sufficiently large 
chromatic number contains holes of $\nu$ consecutive lengths.

\end{abstract}

\section{Introduction}
All graphs in this paper are finite and without loops or parallel edges. A {\em hole} in a graph is an induced subgraph which is a cycle of length
at least four, and a hole is {\em odd} if its length is odd. A {\em triangle}
in $G$ is a three-vertex complete subgraph, and a graph is {\em triangle-free}
if it has no triangle. In this paper we are concerned with the chromatic
number of triangle-free graphs that have no holes of certain specified lengths.

What can we say about the hole lengths in triangle-free graphs with large
chromatic number? There are three well-known conjectures of Gy\'arf\'as~\cite{gyarfas},
the third implying the first two, as follows:

\begin{thm}\label{gyarfasconj}
{\bf Conjecture: }For all $k, \ell$, there exists $n$ such that if $G$ has no clique of cardinality $k$ and has chromatic number at least $n$, then
\begin{itemize}
\item $G$ has an odd hole;
\item $G$ has a hole of length at least $\ell$; and
\item $G$ has an odd hole of length at least $\ell$.
\end{itemize}
\end{thm}
The first conjecture was proved in~\cite{oddholes}, and the second in~\cite{longholes}. 
There are a few other results about the lengths of holes in a graph $G$ with (sufficiently) large chromatic number:
\begin{itemize}
\item $G$ contains a large clique or an even hole~\cite{addario};
\item $G$ contains a large clique or a hole of length 5 or a long hole~\cite{threesteps};
\item $G$ contains a triangle or an odd hole of length at least seven~\cite{threesteps}; and
\item $G$ contains a triangle or a hole of length a multiple of three~\cite{stephan}.
\end{itemize}

Since this paper was submitted, there has been some further progress.  In joint work with Maria Chudnovsky and Sophie Spirkl~\cite{longoddholes},
we proved the third conjecture of Gy\'arf\'as.  Finally, in recent work \cite{holeresidues}, we proved the following result, which gives a common generalization of all the results mentioned above.

\begin{thm}\label{resresult}
For all $k, s,t$, there exists $n$ such that if $G$ has no clique of cardinality $k$ and has chromatic number at least $n$, then
$G$ has a hole of length $s$ modulo $t$.
\end{thm}

In this paper we consider the case $k=3$.  In this case, we show that a far stronger result holds.
The main result of this paper is:
\begin{thm}\label{mainthm}
For all integers $\nu>0$ there exists $n$ such that if $G$ is triangle-free with chromatic number at least $n$, then
for some $t$, $G$ has a hole of length $t+i$ for $1\le i\le \nu$.
\end{thm}
This contains as special cases the $k=3$ cases of all the results mentioned above.
We conjecture that the corresponding result is true for graphs with bounded clique number rather than just triangle-free graphs,
but so far we have made no progress in
proving this.


Let us mention in passing a much more general question, which seems to be interesting even though we cannot answer it.
Let us say a set $F$ of integers is {\em $k$-constricting} if there exists $n$ such that every 
graph with chromatic number at least $n$ contains either a clique with $k$ vertices or a hole with length in $F$. 
Say that $F$ is {\em constricting} if it is $k$-constricting for every $k$.
Which sets are constricting?
Certainly every constricting set is infinite, because there are graphs with arbitrarily
large chromatic number and arbitrarily large girth. 
On the other hand, a consequence of our main result is the following.

\begin{thm}\label{mainthm2}
Let $F$ be an infinite set of positive integers with bounded gaps.  Then $F$ is 3-constricting.
\end{thm}

As noted above, we conjecture that the following more general result should hold.

\begin{thm}{\bf Conjecture:}
Let $F$ be an infinite set of positive integers with bounded gaps.  Then $F$ is constricting.
\end{thm}

The only source we know for examples of infinite sets that are not constricting 
 is the following.
Let $G_1$ be the null graph; for each $i>1$; let $G_i$ be a triangle-free graph with girth
at least $2^{|V(G_{i-1})|}$ and chromatic number at least $i$; and let $F$ be the set of all cycle lengths 
that do not occur in any $G_i$. Then $F$ is not constricting, and yet $F$ has upper density 1. This shows
that not every infinite set is constricting, not even sets with upper density one.
Lower density seems to be closer to the truth. We have not been able to rule out the possibility a set is constricting if and only if it
has strictly positive lower density, although this does not seem likely.  It would be interesting to answer the following question.

 \begin{thm}{\bf Problem:}
Is there an infinite set $F$ of positive integers such that $F$ is constricting and has density 0?
\end{thm}

We conjecture that the answer is positive.  In fact, perhaps an even stronger statement holds.

 \begin{thm}{\bf Problem:}
Is there an infinite set $F=\{f_1,f_2,\dots\}$ of positive integers such that $F$ is constricting and $f_{i+1}-f_i\to\infty$ as $i\to\infty$?
\end{thm}

It would be very interesting to answer these questions even in the case of 3-constricting sets.

\section{Chromatic number and radius}

The proof of \ref{mainthm} 
breaks into three cases, depending on the chromatic number of the subgraphs within a fixed distance of a vertex (even
if we just want to prove the long odd holes conjecture). 
Let us describe this more exactly. If $X\subseteq V(G)$, the subgraph of $G$ induced on $X$ is denoted by $G[X]$,
and we often write $\chi(X)$ for $\chi(G[X])$. The {\em distance} (denoted by $d_G(u,v)$ or $d(u,v)$) between two vertices $u,v$
of $G$ is the length of a shortest path between $u,v$, or $\infty$ if there is no such path.
If $v\in V(G)$ and $\rho\ge 0$ is an integer, $N_G^{\rho}(v)$ or $N^{\rho}(v)$ denotes the set of all vertices with distance
exactly 
$\rho$  from $v$, and $N_G^{\rho}[v]$ or $N^{\rho}[v]$ denotes the set of all vertices with distance at most $\rho$ from $v$.
We denote the maximum over all $v\in V(G)$ of $\chi(N_G^{\rho}[v])$ by $\chi^{\rho}(G)$ (setting $\chi^{\rho}(G)=0$ for the null graph).

Since we are only concerned with triangle-free graphs, it follows that $\chi^1(G)\le 2$, but there may be vertices $v$ such that
$\chi(N^2_G[v])$ is large, and such vertices cause difficulties. If we can find an induced subgraph $H$ with large chromatic number such that
$\chi^2(H)$ is bounded, then we might as well replace $G$ by $H$. If we cannot find such a subgraph, then
we will prove that for all $\ell \ge 5$, $G$ has a hole of length $\ell $ (if its chromatic number is large enough in terms of $\ell $).

Next we assume $\chi^2(G)$ is bounded. If there is an induced subgraph $H$
with large chromatic number and with $\chi^3(H)$ bounded, we might as well pass to that; and if not, we prove
that $G$ contains holes of any fixed length (except very short ones) if $\chi(G)$ is large enough. And the same for $\chi^{\rho}(G)$
for all bounded $\rho$.

Finally, we assume $\chi^{\rho}(G)$ is bounded, for 
some appropriately large $\rho$. (We need $\rho$ to be exponentially large in terms of $\nu$.)
In that case we prove that $G$ contains holes of $\nu$ consecutive lengths (but the smallest of them might be 
arbitrarily large).

Let us say this more precisely.
Let $\nu\ge 0$; a {\em hole $\nu$-interval} in a graph $G$ is a sequence $C_1\l C_{\nu}$ of holes in $G$, such that
$|E(C_{i+1})|=|E(C_i)|+1$ for $1\le i<\nu$ (thus, $\nu$ holes with consecutive lengths).
Let $\mathbb{N}$ denote the set of nonnegative integers, and let $\phi:\mathbb{N}\rightarrow \mathbb{N}$ be a non-decreasing function.
For $\rho\ge 1$, let us say a graph $G$ is {\em $(\rho,\phi)$-controlled} if 
$\chi(H)\le \phi(\chi^{\rho}(H))$ for every induced subgraph $H$ of $G$. Roughly, this says that in every induced subgraph $H$ of $G$ with 
large chromatic number, there is a vertex $v$ such that $H[N^{\rho}_H[v]]$ has large chromatic number.

We will show the following three statements:
\begin{thm}\label{rad2}
Let $\phi:\mathbb{N}\rightarrow \mathbb{N}$ be a non-decreasing function; then for all $\ell\ge 5$ there exists $n$ such that
every $(2,\phi)$-controlled triangle-free graph with chromatic number more than $n$ has an $\ell$-hole.
\end{thm}
\begin{thm}\label{rad3}
Let $\rho>2$ and $\ell\ge 4\rho(\rho+2)$
be integers. For every non-decreasing function
$\phi:\mathbb{N}\rightarrow \mathbb{N}$
there is a non-decreasing function $\phi'$
with the following property.
Let $G$ be a $(\rho,\phi)$-controlled triangle-free graph. Then either $G$ is $(2,\phi')$-controlled or $G$ has an $\ell$-hole.
\end{thm}
\begin{thm}\label{bigrad}
Let $\nu\ge 2$; then there exist $\rho>0$ and a non-decreasing function $\phi$
with the following property. If $G$ is a triangle-free graph then either $G$ is $(\rho,\phi)$-controlled or $G$ admits a hole $\nu$-interval.
\end{thm}

\ref{rad2} might be true when $\ell = 4$ as well, but we have not been able to decide this.
\ref{rad2} is easy for $\ell\le 6$, and in another paper~\cite{threesteps} (with Maria Chudnovsky) we proved it for $\ell=7$,  expecting that
to be the easiest of the open cases. By a happy coincidence, $\ell=7$ turns out to be the one case that is not handled by
the proof method of the present paper. 
Let us see that these three together imply~\ref{mainthm}.

\bigskip

\noindent{\bf Proof of \ref{mainthm}, assuming \ref{rad2}, \ref{rad3}, \ref{bigrad}.}
Let $\nu\ge 2$, and let $\rho$ and $\phi$ be as in \ref{bigrad}.
Let ${\ell}_0 = 4\rho(\rho+2)$, and for $i = 1\l\nu-1$ let ${\ell}_i = {\ell}_0+i$. By \ref{rad3}, for each 
$i\in\{0\l \nu-1\}$ there is a function $\phi'$ as in \ref{rad3}
(with $\ell$ replaced by ${\ell}_i$); define $\phi_i = \phi'$. Thus $\phi_0\l \phi_{\nu-1}$ are all non-decreasing functions; define 
$$\psi(\kappa) = \max(\phi_0(\kappa)\l \phi_{\nu-1}(\kappa))$$
for $\kappa \ge 0$. Thus $\psi$ is non-decreasing. Now by \ref{rad2} (with $\phi$ replaced by $\psi$) for $\ell = 5\l \nu+4$ 
there exists $n$ as in \ref{rad2};
let $n_{\ell} = n$. Let $n=\max(n_5\l n_{\nu+4})$. 

We claim that every triangle-free graph with chromatic number more than $n$ admits a hole $\nu$-interval.
For let $G$ be such a graph, and suppose it admits no hole $\nu$-interval. From the choice of $\rho$ and $\phi$, it follows that
$G$ is $(\rho,\phi)$-controlled. For some $i\in \{0\l \nu-1\}$, $G$ has no ${\ell}_i$-hole; so from the choice of $\phi_i$,
$G$ is $(2,\phi_i)$-controlled and hence $(2,\psi)$-controlled. For some $\ell\in \{5\l \nu+4\}$, $G$ has no $\ell$-hole; 
and so from the choice of $n_{\ell}$,
$\chi(G)\le n_{\ell}\le n$. This proves \ref{mainthm}.~\bbox

The three statements \ref{rad2}, \ref{rad3}, \ref{bigrad} will be proved in separate parts of the paper. 
By far the most difficult is \ref{bigrad}. If all we want is the long odd holes conjecture, then we still need most of 
the two easier results \ref{rad3} and \ref{rad2}, but we could skip most of the proof of \ref{bigrad}; indeed, we need 
nothing after \ref{staircase}.

\section{Radius 2}

In this section we prove \ref{rad2}. We begin with the following:

\begin{thm}\label{2rad}
Let $\phi:\mathbb{N}\rightarrow \mathbb{N}$ be non-decreasing, and let $G$ be triangle-free and $(2,\phi)$-controlled.
\begin{itemize}
\item If $\chi(G)>\phi(2)$ then $G$ has a $5$-hole.
\item If $\chi(G)> \phi(3)$ then $G$ has a $6$-hole.
\item If $\chi(G)>\phi(\phi(\phi(2\phi(2)+2)+1)+1)$ then $G$ has a $7$-hole. 
\end{itemize}
\end{thm}
\Proof
The first statement was proved in~\cite{threesteps}, but we repeat the proof because it is easy. Suppose that $\chi(G)>\phi(2)$, and 
let $v$ be a vertex such that $\chi(G)\le \phi(\chi(N^2[v]))$. It follows that $\chi(N^2[v])>2$, and so there are two adjacent vertices
$x,y\in N^2(v)$. Since $G$ is triangle-free, $x,y,v$, together with two vertices of $N^1(v)$ adjacent to $x,y$ respectively, 
form a $5$-hole.

For the second statement, let $\chi(G)>\phi(3)$, and choose a vertex $v$ such that $\chi(G)\le \phi(\chi(N^2[v]))$. It follows that
$\chi(N^2[v])>3$, and so $\chi(N^2(v))>2$; and hence 
there is an odd hole $P$ in $G[N^2(v)]$. Let $P$ have vertices $p_1\d p_2\c p_n\d p_1$ in order, where $n\ge 5$.
Choose $S\subseteq N^1(v)$ minimal such that every vertex in $V(P)$ has a neighbour in $S$. 
Let $s_i\in S$ be adjacent to $p_i$ for $1\le i\le n$. (Possibly $s_1\l s_5$ are not all distinct.)
For each $s\in S$, some vertex in $P$ is adjacent to $s$ and to no other vertex in $S$, from the minimality
of $S$. Consequently we may assume that $p_3$ is adjacent to $s_3\in S$ and has no other neighbour in $S$. If $p_1$ is nonadjacent to 
$s_3$ then $v\d s_1\d p_1\d p_2\d p_3\d s_3\d v$ is a $6$-hole as required, so we may assume that $p_1$ is adjacent to $s_3$,
and similarly $p_5$ is adjacent to $s_3$. Hence $p_1,p_5$ are nonadjacent since $G$ is triangle-free, and so
$n\ge 7$. If $s_2, s_4$ are nonadjacent to $p_4,p_2$ respectively then 
$v\d s_2\d p_2\d p_3\d p_4\d s_4\d v$ is a $6$-hole, so we may assume that one of $s_2,s_4$ is adjacent both of $p_2,p_4$, say $s_2$.
But then $s_3\d p_1\d p_2\d s_2\d p_4\d p_5\d s_3$ is a $6$-hole.

The third statement is proved in~\cite{threesteps}. This proves \ref{2rad}.~\bbox

\bigskip

Let $X\subseteq V(G)$. A {\em $t$-trellis on $X$ in $G$} is a subgraph $H$ of $G$ with the following properties.
\begin{itemize}
\item $X\subseteq V(H)$, 
and $V(H)\setminus X$ consists of the disjoint union of 
four sets $\{a_1 \l a_t\}$, $\{b_1\l b_t\}$, $\{a_{x,j}: x\in X, 1\le j\le t\}$ and
$\{b_{x,j}:x\in X, 1\le j\le t\}$. 
\item
The edges of $H$ are as follows:
\begin{itemize}
\item $a_jb_j$ for $1\le j\le t$;
\item $xa_{x,j}$ and $xb_{x,j}$ for $x\in X$ and $1\le j\le t$; and
\item $a_{x,j}a_j$ and $b_{x,j}b_j$ for $x\in X$ and $1\le j\le t$.
\end{itemize}
(Thus, to construct $H$ we start with $K_{s,2t}$, with bipartition $X$ and $Y$ say, where $|X|=s$; subdivide all its edges; and then
add a matching pairing up the vertices in $Y$.)
\item 
For all distinct $u,v\in V(H)$, if $u,v$ are adjacent in $G$ and nonadjacent in $H$
then there exist $x,x'\in X$ and $j\in \{1\l t\}$ such that $\{u,v\} = \{a_{x,j},b_{x',j}\}$. (In particular, $X$ is stable.)
\end{itemize}

We also need a modification of this.
An {\em extended $t$-trellis on $X$ in $G$} is a subgraph $H$ of $G$ with the following properties.
\begin{itemize}
\item $X\subseteq V(H)$,
and $V(H)\setminus X$ consists of the disjoint union of
four sets $\{a_0, a_1,\l a_t\}$, $\{b_0,b_1\l b_t\}$, $\{a_{x,j}: x\in X, 0\le j\le t\}$ and
$\{b_{x,j}:x\in X, 0\le j\le t\}$, together with one more vertex $c_0$.
\item
The edges of $H$ are as follows:
\begin{itemize}
\item $a_0c_0$ and $c_0b_0$;
\item $a_jb_j$ for $1\le j\le t$;
\item $xa_{x,j}$ and $xb_{x,j}$ for $x\in X$ and $0\le j\le t$; and
\item $a_{x,j}a_j$ and $b_{x,j}b_j$ for $x\in X$ and $0\le j\le t$.
\end{itemize}
\item 
For all distinct $u,v\in V(H)$, if $u,v$ are adjacent in $G$ and nonadjacent in $H$
then there exist $x,x'\in X$ and $j\in \{0\l t\}$ such that $\{u,v\} = \{a_{x,j},b_{x',j}\}$. 
\end{itemize}

We need both these definitions; we will show that certain graphs contain extended trellises, and to do so
we first show they contain trellises, and then find the extension. 

\begin{thm}\label{allholes}
For every integer $\ell\ge 8$, there exists $t\ge 0$ with the following property. 
Let $G$ be a graph, let $X\subseteq V(G)$ with $|X|= t $, and let $H$ be an extended $t$-trellis on $X$.
Then $G$ has an $\ell$-hole.
\end{thm}
\Proof By Ramsey's theorem, there exists $t\ge 0$ such that if
$\mathcal{A}$ is the set of all triples $(i,i',j)$ with $1\le i<i'\le t$ and $1\le j\le t$,
and we partition $\mathcal{A}$ into two subsets $\mathcal{A}_1$, $\mathcal{A}_2$, then there exist
$R,S\subseteq \{1\l n\}$ with $|R|,|S|\ge \ell$, such that the triples
$(i,i',j)$ with $i<i'\in R$ and $j\in S$ either all belong to $\mathcal{A}_1$ or all belong to $\mathcal{A}_2$.
We claim that $n$ satisfies the theorem. 

For let $G,X, H$ be as in the theorem. Let $X=\{x_1\l x_t\}$, and let us write $a_{i,j}, b_{i,j}$ for $a_{x_i,j}$ and $b_{x_i,j}$ 
respectively. Let $\mathcal{A}_1$ be the set of all triples
$(i,i',j)$ with $1\le i<i'\le t$ and $1\le j\le t$ such that $a_{i,j}, b_{i',j}$ are nonadjacent,
and let $\mathcal{A}_2$ be the set of all such triples such that  $a_{i,j}, b_{i',j}$ are adjacent.
From the choice of $t$, we may assume that for some $k\in \{1,2\}$, $(i,i',j)\in \mathcal{A}_k$ for
all $i,i',j$ with $1\le i<i'\le \ell$ and $1\le j\le \ell$.

For $1\le i<\ell$ let $P_i$ be the path 
$x_i\d a_{i,i+1}\d a_{i+1}\d a_{i+1,i+1}\d x_{i+1}$.
If $k=1$ let $Q_i$ be the path
$x_i\d a_{i,i+1}\d a_{i+1}\d b_{i+1}\d b_{i+1,i+1}\d x_{i+1}$,
and if $k=2$ let $Q_i$ be the path 
$x_i\d a_{i,i+1}\d b_{i+1,i+1}\d x_{i+1}$.
Thus $P_i$ has length four, and $Q_i$ has length five if $k=1$, and three if $k=2$.

Suppose that $\ell$ is a multiple of four, say $\ell=4p$. Then 
the union of $P_1\l P_{p-1}$ and the path $x_1\d a_{1,1}\d a_1\d a_{p,1}\d x_{p}$
is a hole of length $\ell$ as required. Thus we may assume that $\ell$ is not a multiple of four.

If $k=2$, choose integers $p,q\ge 0$ such that 
$\ell=4p+3q$ and $q>0$;
then the union of $Q_i\; (1\le i<q)$, $P_i\;(q\le i < p+q)$, and 
$x_1\d a_{1,1}\d b_{p+q,1}\d x_{p+q}$ is the desired hole.

Thus we may assume that $k=1$. If $\ell\ne 11$, then, since $4$ does not divide $\ell$, 
$\ell$ can be expressed as $4p+5q$ where $p,q$ are nonnegative integers and $q>0$;
and the union of 
$Q_i\; (1\le i<q)$, $P_i\;(q\le i < p+q)$, and 
$x_1\d a_{1,1}\d a_1\d b_1\d  b_{p+q,1}\d x_{p+q}$ is the desired hole. 

Finally we may assume that $\ell=11$. If $a_{1,0},b_{2,0}$ are nonadjacent then the union of $Q_1$ and
$x_1\d a_{1,0}\d a_0\d c_0\d b_0\d b_{2,0}\d x_2$
is the desired hole; while if $a_{1,0},b_{2,0}$ are adjacent then the union of $P_2$,
$x_1\d a_{1,0}\d b_{2,0}\d x_2$, and $x_1\d a_{1,3}\d a_3\d a_{3,3}\d x_3$
is the desired hole. This proves \ref{allholes}.~\bbox

We remark that we only used the ``extended'' part of the trellis in \ref{allholes} for the case $\ell=11$. To prove the result just
for $\ell\ge 8$ and $\ell\ne 11$, the same proof would work for a (non-extended) trellis.

We also need another definition. Let $x\in V(G)$, let $N$ be some set of neighbours of $x$, and let $C\subseteq V(G)$
be disjoint from $N\cup \{x\}$, such that every vertex in $C$ is nonadjacent to $x$ and has a neighbour in $N$. In this situation
we call $(x,N)$ a {\em cover} of $C$ in $G$. For $C,X\subseteq V(G)$, a {\em multicover of $C$} in $G$ 
is a family $(N_x:x\in X)$
such that
\begin{itemize}
\item for each $x\in X, (x,N_x)$ is a cover of $C$;
\item for all distinct $x,x'\in X$, $x'$ has no neighbour in $\{x\}\cup N_x$ (and in particular all the sets
$\{x\}\cup N_x$ are pairwise disjoint).
\end{itemize}
If in addition we have 
\begin{itemize}
\item for all distinct $x,x'\in X$, no vertex in $N_{x'}$ has a neighbour in $N_x$,
\end{itemize}
we call $(N_x:x\in X)$ an {\em independent multicover}.

\begin{thm}\label{gettrellis}
For all $t,\kappa\ge 0$, there exist $\tau, m\ge 0$ with the following property.
Let $G$ be a triangle-free graph such that every induced subgraph of $G$ with chromatic number
more than $\kappa$ has a $5$-hole. 
Let $C\subseteq V(G)$ with chromatic number more than $\tau$; and let $(N_x:x\in X)$ be a multicover of $C$ in $G$ with $|X|\ge m$.
Then there exist $Y\subseteq X$ with $|Y|=t$ and an extended $t$-trellis on $Y$ in $G$.
\end{thm}
\Proof For $0\le s\le t$ let $m_s'= 5t\cdot 5^{t-s}$, and let $m' = m'_0$.
For $0\le s\le t$ let $m_s = 5t(20m')^{t-s}$, and let $m = m_0$.
Let $\tau_t' = \kappa+1$, and for $s=t-1\l 0$ let 
$$\tau_s'= 5(m_s'+1)+ 5^{m_s'} \tau_{s+1}'.$$
Let $\tau'=\tau_0$.
Let $\tau_t= \kappa+1$, and for $s=t-1\l 0$ let 
$$\tau_s= 5(m_s+1)+m_s^{m'+1}5^{m_s}\tau'+ 2^{m_s}5^{m_s} \tau_{s+1}.$$
Let $\tau=\tau_0$.
We claim that $\tau,m$ satisfy the theorem. 
Let 
$G$ be a triangle-free graph such that every induced subgraph of $G$ with chromatic number
more than $\kappa$ has a $5$-hole. We shall prove  the following, which implies the theorem:
\\
\\
(1) {\em Let $C\subseteq V(G)$ and let $(N_x:x\in X)$ be a multicover of $C$, such that either
\begin{itemize}
\item $\chi(C)>\tau$ and $|X|= m$, or
\item $\chi(C)>\tau'$ and $|X|= m'$ and $(N_x:x\in X)$ is independent.
\end{itemize}
Then there exist $Y\subseteq X$ with $|Y|=t$ and an extended $t$-trellis on $Y$ in $G$.}
\\
\\
If $X'\subseteq X$, and $N'_x\subseteq N_x$ for each $x\in X'$, and $C'\subseteq C$, and every vertex in $C'$ has a neighbour in $N_x'$
for each $x\in X'$, then $(N_x':x\in X')$ is a multicover of $C'$, 
and we say it is {\em contained in} $(N_x:x\in X)$.
Consequently, to prove (1), we may assume that:
\\
\\
(2) {\em Either
\begin{itemize}
\indentitem\item[{\bf (Case 1)}]  $\chi(C)>\tau$ and $|X|\ge m$ and there do not exist $C'\subseteq C$ with $\chi(C')>\tau'$ 
and $X'\subseteq X$ with $|X'|\ge m'$ and an
independent multicover $(N_x':x\in X')$ of $C'$ contained in $(N_x:x\in X)$, or
\indentitem\item[{\bf (Case 2)}]  $\chi(C)>\tau'$ and $|X|\ge m'$ and $(N_x:x\in X)$ is independent.
\end{itemize}
}

Now we construct a $t$-trellis on a subset of $X$ as follows (later we will enlarge it to an extended trellis). 
We begin with the $0$-trellis on $X$, $H_0$ say, and let $C_0 = C$.
Inductively, suppose that $s< t$, and we have constructed an $s$-trellis $H_s$ on a subset $X_s\subseteq X$,
with vertex set the disjoint union of $X_s$, $\{a_1,\l a_s\}$, $\{b_1\l b_s\}$, $\{a_{x,j}: x\in X_s, 1\le j\le s\}$ and
$\{b_{x,j}:x\in X_s, 1\le j\le s\}$ in the usual notation, and a subset $C_s\subseteq C$, satisfying:
\begin{itemize}
\item $a_j,b_j\in C$ for $1\le j\le s$;
\item $a_{x,j}, b_{x,j}\in N_x$ for each $x\in X_s$ and $1\le j\le s$; 
\item in case 1, $|X_s|= m_s$, and in case 2, $|X_s|= m_s'$;
\item no vertex in $V(H_s)$ has a neighbour in $C_s$;
\item for each $v\in C_s$ and each $x\in X_s$, there is a neighbour of $v$ in $N_x$ that has no neighbour in $V(H_s)$
except $x$; and
\item in case 1, $\chi(C_s)>\tau_s$, and in case 2, $\chi(C_s)> \tau_s'$.
\end{itemize}
For each $x\in X_s$, let $N_x'$ be the set of vertices in $N_x$ with no neighbour in $V(H_s)$ except $x$. Then 
$(N_x':x\in X_s)$ is a multicover of $C_s$, and is independent in case 2.

Since $\chi(C_s)>\tau_s'\ge \kappa$, there is a $5$-hole $P$ in $G[C_s]$, with vertices $p_1\d p_2\c p_5\d p_1$ say, in order.
For each $x\in X_s$, and $1\le i\le 5$, let $D_i(x)$ be the set of vertices in $N_x'$ adjacent to $p_i$, and select
$d_i(x)\in D_i(x)$. Thus the union of $V(P)$ and $\{d_i(x):1\le i\le 5,x\in X_s\}$ has cardinality at most $5(|X_s|+1)$,
and since $G$ is triangle-free, there exists $C_s^1\subseteq C_s$ with 
$\chi(C_s^1)\ge \chi(C_s)-5(|X_s|+1)$, such that
no vertex in $C_s^1$ is adjacent to any of the vertices $d_i(x)$ or to any vertex in $P$ (and in particular, 
$C_s^1\cap V(P)=\emptyset$).

For each $x\in X_s$, no vertex is in more than two of $D_1(x)\l D_5(x)$, because $G$ is triangle-free.
For each $v\in C_s^1$ and $x\in X_s$, since $v$ has a neighbour in $N_x'$, it follows that there exist
adjacent vertices $p_k,p_{k+1}$ of $P$ such that some neighbour of $v$ belongs to
$N_x'\setminus (D_k(x)\cup D_{k+1}(x))$ (reading subscripts modulo $5$); choose some such $k$ and define $c_x(v) = k$.
There are $5^{|X_s|}$ possibilities for the $X_s$-tuple $(c_x(v):x\in X_s)$, and so there exists $C_{s}^2\subseteq C_s^1$
with $\chi(C_s^2)\ge \chi(C_s^1)/5^{|X_s|}$, such that $c_x(v) = c_x(v')$ for all $x\in X_s$ and all $v,v'\in C_s^2$. Moreover,
since there are only five possibilities for $c_x(v)$, there exists $k\in \{1\l 5\}$ and $Y_s\subseteq X_s$
with $|Y_{s}|=|X_s|/5$ such that
$c_x(v)=k$ for all $x\in Y_{s}$ and $v\in C_s^2$. 
Thus $\chi(C_s^2)\ge (\chi(C_s)-5(|X_s|+1))/5^{|X_s|}$, and so in case 1
$$\chi(C_s^2) > (\tau_s-5(m_s+1))/5^{m_s}=m_s^{m'+1}\tau'+ 2^{m_s}\tau_{s+1},$$
and in case 2
$$\chi(C_s^2)> (\tau_s'-5(m_s'+1))/5^{m_s'}=\tau_{s+1}'.$$
Let $a_{s+1}=p_k$ and $b_{s+1}=p_{k+1}$, and for each $x\in Y_{s}$
let $a_{x,s+1} = d_k(x)$ and $b_{x,s+1}=d_{k+1}(x)$. To complete the inductive
definition it remains to define $X_{s+1}$ and $C_{s+1}$.

In case 2 we define $X_{s+1} = Y_s$ and $C_{s+1} = C_s^2$; so we assume we are in case 1.
The issue that we need to handle in this case is that for $v\in C_{s+1}$ and $x\in Y_s$, while we know that
$v$ has a neighbour $u\in N_x'$ that has no neighbours in $V(H_s)$ except $x$, it may be that every such neighbour $u$
is adjacent to one of $a_{x',s+1},b_{x',s+1}$ for some $x'\in Y_s$. We shall show that if this happens for ``many'' 
choices of $v$ then we can move into case 2.

Let $Z$ be the union of the sets $\{a_{x,s+1}, b_{x,s+1}\}$ over all $x\in Y_{s}$; then $|Z| = 2m_s/5\le m_s$. Let $z\in Z$, and let
$Y\subseteq Y_s$ with $|Y|=m'$.
Let $D_{z,Y}$ be the set of vertices $v\in C_s^2$ such that for each $x\in Y$ there exists a vertex in $N_x'$
adjacent to both $v,z$. For each $x\in Y$, let $N_x''$ denote the set of vertices in $N_x'$ adjacent to $z$; then
$(N_x'':x\in Y)$ is a multicover of $D_{z,Y}$; and it is independent, since $G$ is triangle-free. Since we are in
case 1, it follows that $\chi(D_{z,Y})\le \tau'$. Now let $D_z$ denote the set of vertices $v\in C_s^2$
such that for at least $m'$ values of $x\in Y_s$ there exists a vertex in $N_x'$
adjacent to both $v,z$; that is, $D_z$ is the union of the sets $D_{z,Y}$ over all choices of $Y$. Since  there are
only at most $m_s^{m'}$ choices of $Y$, it follows that
$\chi(D_z)\le m_s^{m'}\tau'$. Thus the union of the sets $D_z$ over all $z\in Z$ has chromatic number at most $m_s^{m'+1}\tau'$,
and so there exists $C_s^3\subseteq C_s^2$ with 
$$\chi(C_s^3)\ge \chi(C_2^2)-m_s^{m'+1}\tau'> 2^{m_s}\tau_{s+1},$$
such that for every $v\in C_s^3$,
and every $z\in Z$, there are fewer than $m'$ values of $x\in Y_{s}$ such that some vertex in $N_x'$ is
adjacent to both $v,z$. 

Fix $v\in C_s^3$ for the moment, and make a digraph $J_v$ with vertex set $Y_{s}$ in which  for distinct $x,y\in Y_{s}$,
$y$ is adjacent from $x$ in $J_v$  if some vertex in $N_y'$ is adjacent to $v$ and to one of $a_{x,s+1}, b_{x,s+1}$.
We have just seen that for all $v$, every vertex of the digraph $J_v$ has indegree in $J$ at most $2m'-2$. It follows
that in $J_v$, some vertex has indegree plus outdegree at most $4m'-4$, and the same holds for every nonnull 
subdigraph of $J_v$; and so
the undirected graph underlying $J_v$ can be $4m'$-coloured. Hence there is a subset $U_v$ say of $Y_{s}$
of cardinality $|Y_{s}|/(4m')=m_{s+1}$ such that no edge of $J_v$ has both ends in $U_v$. There are only $2^{|Y_{s}|}$
possibilities for $U_v$, and so there exists $C_s^4\subseteq C_s^3$ with 
$$\chi(C_s^4)\ge \chi(C_s^3)/2^{|Y_{s}|}>\tau_{s+1}$$
such that the sets $U_v$ are equal for all $v\in C_s^4$. Let $X_{s+1}$ be this common value of $U_v$, and let $C_{s+1}=C_s^4$.
This completes the definition of $C_{s+1}$ in case 1.

In both cases, the pairs $a_j,b_j (1\le j\le s+1)$ and the vertices
$a_{x,j}, b_{x,j} (x\in X_{s+1}, 1\le j\le s+1)$ define an $(s+1)$-trellis $H_{s+1}$ on $X_{s+1}$, and no vertex in $H_{s+1}$
has a neighbour in $C_{s+1}$, and for all $v\in C_{s+1}$ and $x\in X_{s+1}$, some neighbour of $v$ in $N_x$ has no neighbour
in $V(H_{s+1})$ except $x$. This
completes the inductive definition of $H_s$ and $C_s$ for $0\le s\le t$.

Thus there is a $t$-trellis on the set $X_t$, where $|X_t| = 5t$; next we need to convert it to an extended
$t$-trellis on a subset of $X_t$ of cardinality $t$. 
With the same notation as before (with $s=t$),
since $\chi(C_t)>\tau_t'> \kappa$, there is a $5$-hole $P$ in $G[C_t]$, with vertices $p_1\d p_2\c p_5\d p_1$ say, in order.
Let $x\in X_t$; a {\em handle} for $x$ means a $3$-vertex path $a\d c\d b$ of $P$ such that some vertex in $N_x'$ is adjacent to $a$, and not
to $b,c$, and some vertex in $N_x'$ is adjacent to $b$ and not to $a,c$. 
We claim that there is a handle for $x$.
Choose $S\subseteq N_x'$ minimal such that
every vertex in $V(P)$ has a neighbour in $S$. For $1\le i\le 5$, choose $s_i\in S$ adjacent to $p_i$.
Suppose first that some $s_1\in S$ has only one neighbour in $V(P)$, say $p_1$.
Then no other vertex in $S$ is adjacent to $p_1$, from the minimality of $S$, and since $s_3$ is nonadjacent to $p_2$
it follows that $p_1\d p_2\d p_3$ is a handle for $x$. We may assume therefore that each $s_i$
has at least two (and hence exactly two) neighbours in $V(P)$. Let $s_1$ be adjacent to $p_1,p_4$ say. From the minimality
of $S$, one of $p_1,p_4$ has no more neighbours in $S$, say $p_1$. But then again $p_1\d p_2\d p_3$ is a handle for $x$.
This proves the claim that for each $x\in X_t$ there is a handle for $x$. Since there are only five possibilities
for handles, there exists $X_0\subseteq X_t$ with $|X_0|=|X_t|/5=t$ such that every vertex in $X_0$ has the same handle,
say $a_0\d c_0\d b_0$. For each $x\in X_0$ let $a_{x,0}\in N_x'$ be adjacent to $a_0$ and not to $b_0,c_0$, and let $b_{x,0}$
be adjacent to $b_0$ and not to $a_0, c_0$. Then 
the pairs $a_j,b_j (1\le j\le t)$, the path $a_0\d c_0\d b_0$, and the vertices
$a_{x,j}, b_{x,j} (x\in X_{s+1}, 0\le j\le s+1)$ define an extended $t$-trellis on $X_0$. This proves \ref{gettrellis}.~\bbox

\bigskip

From \ref{allholes} and \ref{gettrellis} we deduce:
\begin{thm}\label{usecover}
For all $\kappa\ge 0$ and $\ell\ge 8$, there exist $\tau, m\ge 0$ with the following property.
Let $G$ be a triangle-free graph such that every induced subgraph of $G$ with chromatic number
more than $\kappa$ has a $5$-hole.
Let $C\subseteq V(G)$ with chromatic number more than $\tau$; and let $(N_x:x\in X)$ be a multicover of $C$ with $|X|\ge m$.
Then $G$ has an $\ell$-hole.
\end{thm}

Let $G$ be a graph and let $t\ge 0$ be an integer. A {\em $t$-cable} in $G$ consists of:
\begin{itemize}
\item $t$ distinct vertices $x_1\l x_t$, pairwise nonadjacent;
\item for $1\le i\le t$, a subset $N_i$ of the set of neighbours of $x_i$, such that the sets $N_1\l N_t$ are pairwise disjoint;
\item for $1\le i\le t$, disjoint subsets $Z_{i,i+1}\l Z_{i,t}, Y_{i}$ of $N_i$; and
\item a subset $C\subseteq V(G)$ disjoint from $\{x_1\l x_t\}\cup N_1\cup\cdots\cup N_t$
\end{itemize}
satisfying the following conditions:
\begin{itemize}
\item for $1\le i\le t$, every vertex in $C$ has a neighbour in $Y_{i}$, and has no neighbours in $Z_{i,j}$ for $i+1\le j\le t$, and is nonadjacent to $x_i$;
\item for $i<j\le t$, $x_i$ has no neighbours in $N_j$;
\item for $i<j<k\le t$, there are no edges between $Z_{i,j}$ and $N_k$;
\item for all $i<j\le t$, either
\begin{itemize}
\item $Z_{i,j}=\emptyset$ and $x_j$ has no neighbours in $Y_i$, or
\item every vertex in $N_j$ has a neighbour in $Z_{i,j}$ and has no neighbours in $Y_i$.
\end{itemize}
\end{itemize}

We call $C$ the {\em base} of the $t$-cable, and say $\chi(C)$ is the {\em chromatic number} of the $t$-cable.
Given a $t$-cable in this notation, let $I\subseteq \{1\l t\}$;
then (after appropriate renumbering) the vertices $x_i\;(i\in I)$, the sets $N_{i}\; (i\in I)$, the sets $Z_{i,j}\; (i,j\in I)$, the sets $Y_i\; (i\in I)$
and $C$ define an $|I|$-cable; we call this a {\em subcable}.

Thus there are two types of pair $(i,j)$ with $i<j\le t$, and we aim next to apply Ramsey's theorem on these pairs to get a large subcable
where all the pairs have the same type. 
Two special kinds of $t$-cables are therefore of interest: {\em type 1} $t$-cables, where for all $i<j\le t$, $Z_{i,j}=\emptyset$ and $x_j$ has no neighbours in $Y_i$, and
{\em type 2} $t$-cables, where for all $i<j\le t$, every vertex in $N_j$ has a neighbour in $Z_{i,j}$ and has no neighbours in $Y_i$.
A type 1 $t$-cable with base $C$ is just a 
multicover of $C$ in disguise, so from \ref{usecover} we have:

\begin{thm}\label{usetype1}
For all $\kappa\ge 0$ and $\ell\ge 8$, there exist $\tau, m\ge 0$ with the following property.
Let $G$ be a triangle-free graph such that every induced subgraph of $G$ with chromatic number
more than $\kappa$ has a $5$-hole.
If $G$ admits a type 1 $m$-cable with chromatic number more than $\tau$,
then $G$ has an $\ell$-hole.
\end{thm}

We need a similar theorem for type 2 cables.
\begin{thm}\label{usetype2}
Let $G$ be a triangle-free graph.
For all $\ell\ge 5$, 
if $G$ admits a type 2 $(\ell-3)$-cable with nonnull base,
then $G$ has an $\ell$-hole.
\end{thm}
\Proof Let $t=\ell-3$ (and so $t\ge 2$) and assume $G$ contains a type 2 $t$-cable with nonnull base. In the usual notation, let $v\in C$. 
Since every vertex in $C$ has a neighbour in $Y_t$, there exists $y_t\in Y_t$ adjacent to $v$.
Since every vertex in $N_t$ has a neighbour in $Z_{t-1,t}$, there exists $z_{t-1}\in Z_{t-1,t}$ adjacent to $y_t$. Similarly
for $i = t-2,t-3\l 1$ there exists $z_i\in Z_{i,i+1}$ such that $z_{i+1}$ is adjacent to $z_i$. Thus $z_1\d z_2\c z_{t-1}\d y_t$ is a path.
It is induced; for if $i,j\le t$ and $j\ge i+2$ then $z_i$ has no neighbour in $N_j$, since $z_i\in Z_{i,i+1}$. Since $x_1$
is adjacent to $z_1$ and to none of $z_2\l z_{t-1}, y_t$ (because $t\ge 2$ and $x_1$ has no neighbours in $N_j$ for $j>1$), 
and $v$ is adjacent to $y_t$ and nonadjacent to $x_1,z_1\l z_{t-1}$, it follows that
$$x_1\d z_1\d z_2\c z_{t-1}\d y_t\d v$$ 
is an induced path. Now $v$ has a neighbour $y_1\in Y_1$; and we claim that $y_1$ is nonadjacent to $z_1\l z_{t-1}, y_t$.
Certainly $y_1,z_1$ are nonadjacent, since they are both adjacent to $x_1$ and $G$ is triangle-free. For $2\le j\le t-1$,
$y_1$ is nonadjacent to $z_j$ since every vertex in $N_j$ has no neighbours in $Y_1$. For the same reason, $y_1$ is nonadjacent to $y_t$, since $t>1$.
Consequently 
$$x_1\d z_1\d z_2\c z_{t-1}\d y_t\d v\d y_1\d x_1$$ 
is a hole of length $t+3=\ell$. This proves \ref{usetype2}.~\bbox

We deduce:

\begin{thm}\label{usecable}
For all $\kappa\ge 0$ and $\ell\ge 8$, there exist $t,\tau\ge 0$ with the following property.
Let $G$ be a triangle-free graph such that every induced subgraph of $G$ with chromatic number
more than $\kappa$ has a $5$-hole. If $G$ admits a $t$-cable with chromatic number more than $\tau$ then $G$ has an $\ell$-hole.
\end{thm}
\Proof
Let $m$, $\tau$ be as in \ref{usetype1}. Let $n=\ell-3$.
Let $t$ equal the Ramsey number $R(m,n)$; that is, the smallest integer $t$ such for for every partition of the edges of $K_t$ into two sets, there is either
a $K_m$ subgraph with all edges in the first set, or a $K_n$ with all edges in the second. We claim that $t,\tau$ satisfy the theorem.

For let $G$ admit a $t$-cable with base $C$ and chromatic number more than $\tau$.
By Ramsey's theorem either 
\begin{itemize}
\item there exists $I\subseteq \{1\l t\}$ with $|I|=m$ such that for all $i,j\in I$ with $i<j$, every vertex in $N_j$ has a neighbour in $Z_{i,j}$ 
and has no neighbours in $Y_i$, or
\item there exists $I\subseteq \{1\l t\}$ with $|I|=n$ such that for all $i,j\in I$ with $i<j$, $Z_{i,j}=\emptyset$ and $x_j$ has no neighbours in $Y_i$.
\end{itemize}
Thus either there is an $m$-subcable of type 1, or an $n$-subcable of type 2, with base $C$ in each case. In the first case the result follows from
\ref{usetype1}, and in the second from \ref{usetype2}. This proves \ref{usecable}.~\bbox

\begin{thm}\label{getcable}
Let $\phi:\mathbb{N}\rightarrow\mathbb{N}$ be non-decreasing, and let $t,\tau\ge 0$.
Then there exists $\tau'$ with the following property.
Let $G$ be a triangle-free graph such that
$G$ is $(2,\phi)$-controlled and $\chi(G)>\tau'$. Then $G$ admits a $t$-cable with chromatic number more than $\tau$.
\end{thm}
\Proof
Let $\tau_t=\tau$, and for $s=t-1\l 0$ let $\tau_s=\phi(2^s\tau_{s+1}+1)$; and let $\tau'=\tau_0$. We claim that $\tau'$ satisfies the theorem.
For let $G$ be a triangle-free graph such that
$G$ is $(2,\phi)$-controlled and $\chi(G)>\tau'$. Consequently $G$ admits a $0$-cable with chromatic number more than $\tau_0$. We claim that 
for $s=1\l t$, $G$ admits an $s$-cable with chromatic number more than $\tau_s$. For suppose the result holds for some $s<t$; we prove it also holds for $s+1$.
In the usual notation, since $\chi(C)>\tau_s=\phi(2^s\tau_{s+1}+1)$, there exists $x_{s+1}\in C$ such that $\chi(N^2_{G[C]}[x_{s+1}])>2^s\tau_{s+1}+1$, and hence 
$\chi(N^2_{G[C]}(x_{s+1}))> 2^s\tau_{s+1}$.
Let $D=N^2_{G[C]}(x_{s+1})$. For each $v\in D$, and $1\le i\le s$, if some neighbour of $v$ in $Y_i$ is nonadjacent to $x_{s+1}$ define $c_i(v) = 1$, and otherwise
define $c_i(v)=2$. There are only $2^s$ possibilities for the $s$-tuple $(c_1(v)\l c_s(v))$, and so there exists $C'\subseteq D$ with $\chi(C')\ge \chi(D)/2^s>\tau_{s+1}$
and an $s$-tuple $(c_1\l c_s)$ such that $c_i(v)=c_i$ for all $v\in C'$ and $1\le i\le s$. 

Let $N_{s+1} = Y_{s+1}'$ be the set of neighbours of $x_{s+1}$ in $C$.
For $1\le i\le s$ define $Z_{i,s+1}, Y_i'\subseteq Y_i$ as follows: 
\begin{itemize}
\item if $c_i = 1$, let $Y_i'$ be the set of vertices in $Y_i$ nonadjacent to $x_{s+1}$,
and let $Z_{i,s+1}=\emptyset$ 
\item if $c_i=2$, let $Y_i'$ be the set of vertices in $Y_i$ adjacent to $x_{s+1}$, and let $Z_{i,s+1}$ be the set of vertices 
in $Y_i$ nonadjacent to $x_{s+1}$.
\end{itemize}
Note that in the second case, no vertex in $Z_{i,s+1}$ has a neighbour in $C'$, and no vertex in $Y_i'$ has a neighbour in $Y_{s+1}'$.
It follows that $x_1\l x_{s+1}$, the sets $N_1\l N_{s+1}$, the sets $Z_{i,j}$ for $1\le i<j\le s+1$, the sets $Y_i'$ for $1\le i\le s+1$, and $C'$, define 
an $(s+1)$-cable with chromatic number more than $\tau_{s+1}$.

This proves that $G$ admits a $t$-cable with chromatic number more than $\tau_t=\tau$, and so proves \ref{getcable}.~\bbox

\bigskip

Let us put these pieces together to prove \ref{rad2}, which we restate:
\begin{thm}\label{rad22}
Let $\phi:\mathbb{N}\rightarrow \mathbb{N}$ be a non-decreasing function; then for all $\ell\ge 5$ there exists $n$ such that
every $(2,\phi)$-controlled triangle-free graph with chromatic number more than $n$ has an $\ell$-hole.
\end{thm}
\Proof
If $l\le 7$ the result follows from \ref{2rad}, so we may assume that $l\ge 8$. 
Let $t,\tau$ be as in \ref{usecable}, taking $\kappa = \phi(2)$; and let $\tau'$ be as in \ref{getcable}. Let $n=\tau'$.
We claim that $n$ satisfies the theorem. For let $G$ be a $(2,\phi)$-controlled triangle-free graph with chromatic number more than $n$. 
By \ref{2rad}, every induced subgraph of $G$
with chromatic number more than $\kappa$ has a $5$-hole. By
\ref{getcable}, $G$ admits a $t$-cable with chromatic number more than $\tau$; and by \ref{usecable}, $G$ has an $\ell$-hole. This proves \ref{rad22}.~\bbox

\bigskip

The second conjecture of \ref{gyarfasconj} is proved in~\cite{longholes}, but if we just wanted to prove it for triangle-free graphs,
rather than the full strength of \ref{mainthm}, the remainder of the paper is not needed; let
us explain why. The following is proved in~\cite{threesteps} (the proof just takes a few lines):

\begin{thm}\label{longholelemma}
Let $\ell\ge 3$ and $\kappa\ge 1$ be integers, and let $G$ be a graph with no hole of length more than $\ell$, such that
$\chi(N(v)),\chi(N^2(v))\le \kappa$ for every vertex $v$.
Then $\chi(G)\le (2\ell-2)\kappa$.
\end{thm}

\bigskip

For each $\kappa\ge 0$, let $\phi(\kappa) = (2\ell-2)\kappa$. It follows from \ref{longholelemma} that if $G$ has no
hole of length more than $\ell$, and $H$ is an induced subgraph
of $G$ with $\chi(H)>\phi(\kappa)$, then $\chi(N^2_H[v])>\kappa$  for some vertex $v$ of $H$; that is,
$G$ is $(2,\phi)$-controlled. Then from \ref{rad22} it follows that $\chi(G)$ is bounded, which proves the second 
assertion of \ref{gyarfasconj} for triangle-free graphs. Indeed, we don't even need
all of \ref{rad22}; instead of an $\ell$-hole, we are content with a hole of length at least $\ell$, and with this modification \ref{rad22} is easier to prove.
For instance, we could get by with trellises instead of extended trellises, since holes of length 11 are of no significance, and indeed
we  could just use 1-subdivisions
of a large $K_{n,n}$ instead of trellises,  since we are not picky about the exact length of the hole.

Trellises give us a long odd hole, but this does not prove the third conjecture of \ref{gyarfasconj}, since we needed to use
\ref{longholelemma}. If
our goal is the long odd holes conjecture, there will be parts of the proof we can skip, but not yet.

\section{Bounded radius}

In this section we prove \ref{rad3}, which we restate, somewhat reformulated:

\begin{thm}\label{boundedrad}
Let $\phi:\mathbb{N}\rightarrow\mathbb{N}$ be non-decreasing, and let $\rho>2$ and $\ell\ge 4\rho(\rho+2)$ 
be integers. There is a non-decreasing function $\phi':\mathbb{N}\rightarrow \mathbb{N}$, 
with the following property.
Let $G$ be a triangle-free graph with no $\ell$-hole such that 
$G$ is $(\rho,\phi)$-controlled.
Then $G$ is $(2,\phi')$-controlled.
\end{thm}

\ref{boundedrad} follows immediately from the following.

\begin{thm}\label{increaserad}
Let $\phi:\mathbb{N}\rightarrow\mathbb{N}$ be non-decreasing, and let $\rho>2$  and $\ell\ge 4\rho(\rho +2)$           
be integers. There is a non-decreasing function $\phi':\mathbb{N}\rightarrow \mathbb{N}$, 
with the following property.
Let 
$G$ be a triangle-free graph with no $\ell$-hole such that 
$G$ is $(\rho,\phi)$-controlled.
Then $G$ is $(\rho-1,\phi')$-controlled.
\end{thm}
\Proof
Let $\ell = 2\alpha\rho + \beta$, where $\alpha\ge 0$ is an integer and $0\le \beta<2\rho$. Since $\ell\ge 4\rho(\rho+2)$, it follows that
$\alpha\ge 2\rho+4$. For $\kappa\in \mathbb{N}$, 
define $\mu_{\alpha+2}(\kappa) = \phi(0)+1$, and for $h=\alpha+2\l 2$ define 
$$\mu_{h-1}(\kappa)= (\rho+1)\kappa +\phi(\phi(\mu_h(\kappa)+\kappa) +(2\rho+2)\kappa),$$
and $\mu_0(\kappa)=\phi(\mu_1(\kappa)+\kappa)$. Define $\phi'(\kappa) = \mu_0(\kappa)$. 
We see that $\phi'$ is non-decreasing.

Let
$G$ be a triangle-free graph with no $\ell$-hole such that
$G$ is $(\rho,\phi)$-controlled.
We will show that $G$ is $(\rho-1,\phi')$-controlled. Let $\kappa\in \mathbb{N}$, such that
$\chi^{\rho-1}(G)\le \kappa$; we must show that $\chi(G)\le \mu_0(\kappa)$. (If so, then the same argument
applied to every induced subgraph $H$ of $G$ and every $\kappa$ shows that $G$ is $(\rho-1,\phi')$-controlled.) Suppose not.

Let $v\in V(G)$. Let $T$ be a path $v=t_0\d t_1\c t_{\rho}$, such that $d_G(v,t_{\rho})=\rho$. For the moment fix such a path $T$.
Let us say a path $P$ is a {\em $(v,T)$-extension} if it has the following properties, where $P$ has vertices 
$p_0\d p_1\c p_n$ in order:
\begin{itemize}
\item $P$ is induced, and $p_0=t_\rho$, and $n\ge \rho$;
\item $d_G(v,p_i)=\rho$ for $0\le i\le n$;
\item $d_G(t_i,p_j)\ge \rho$ for $0\le i\le \rho$ and $\rho\le j\le n$; and
\item $d_G(p_i,p_n)\ge \rho$ for $0\le i\le n-\rho$.
\end{itemize}
(1) {\em If $P$ as above is a $(v,T)$-extension, then $P\cup T$ is an induced path of length $\rho+n$.}
\\
\\
Because $T$ is induced since  $d_G(v,t_{\rho})=\rho$, and $P$ is induced by hypothesis. Moreover $V(P)\cap V(T)=\{t_\rho\}$
since $d_G(v,t_i)<\rho$ for $0\le i<\rho$, and $d_G(v,p_i)=\rho$ for $0\le i\le n$. Suppose that some $t_i$ is adjacent to some $p_j$,
where $i<\rho$ and $j>0$.
Since $d_G(v,p_j)=\rho$ and $d_G(v,t_i)=i<\rho$, it follows that $i=\rho-1$. Now $j\ne 1$ since $G$ is triangle-free,
so $j\ge 2$. Since $d_G(t_{\rho-1},p_k)\ge \rho$ for $\rho\le k\le n$, it follows that $j<\rho$. Then the path 
$t_{\rho-1}\d p_j\d p_{j+1}\c p_{\rho}$ has length $\rho-j+1<\rho$, a contradiction since $d_G(t_{\rho-1}, p_{\rho})\ge \rho$. This proves (1).

\bigskip
Let $P,P'$ both be $(v,T)$-extensions. We say they are {\em parallel} if
the last three vertices of $P$ are the 
same as the last three of $P'$, and in particular the last vertices of $P,P'$ are equal.
\\
\\
(2) {\em Let $P_1\l P_k$ be $(v,T)$-extensions, pairwise parallel. Then there exists $s\in \{2\rho,2\rho-2,2\rho-4\}$
such that $G$ has holes of lengths $|E(P_1)|+s\l |E(P_k)|+s$.}
\\
\\
Let $z$ be the common last vertex of $P_1\l P_k$, and choose a path 
$Z$ between $v,z$ of length $\rho$. Since $T\cup Z$ is connected, there is an
induced path $Q$ between $t_\rho,z$ with $V(Q)\subseteq V(T\cup Z)$. Let us first examine the length of $Q$.
Let $Z$ have vertices $z_0\d z_1\c z_{\rho}$, where $z_0=v$ and $z_{\rho}=z$.
If no vertex in $\{z_1\l z_{\rho}\}$ has a neighbour in $\{t_1\l t_{\rho}\}$, then the two sets are disjoint, and
$Q=T\cup Z$
and hence has length $2\rho$. We assume then that some $z_j\in \{z_1\l z_{\rho}\}$
is adjacent to some $t_i\in \{t_1\l t_{\rho}\}$. Since $d_G(t_i,z)\ge \rho$ from the definition of
a $(v,T)$-extension, the path $t_i\d z_j\d z_{j+1}\c z_{\rho}$ has length at least $\rho$, and so $j=1$. Since 
$z_j$ is adjacent to $t_0=v$, and $G$ is triangle-free, it follows that $i\ge 2$. Since $d_G(v,t_{\rho})=\rho$,
it follows that $i=2$. So there is only one such edge, and in particular the two sets $\{z_1\l z_{\rho}\}$, 
$\{t_1\l t_{\rho-1},t_{\rho}\}$ are disjoint, and $Q$ has length $2\rho-2$. We have proved then that $Q$ has length $2\rho$ or $2\rho-2$.

Now let $P$ be one of $P_1\l P_k$, and let $P$ have vertices $p_0\d p_1\c p_n$ in order. Thus $p_0=t_{\rho}$ and $p_n=z_{\rho}=z$.
Both $P,Q$ are induced, and their interiors are disjoint, since every vertex $x$ of the interior of $Q$ belongs to one
of $V(Z)\setminus \{z\}, V(T)\setminus \{t_{\rho}\}$ and hence satisfies $d_G(v,x)<\rho$, while $d_G(v,x)=\rho$ for every vertex $x$ of the 
interior of $P$. Suppose then that some vertex $x$ in the interior of $Q$ has a neighbour $p_j\in\{p_1\l p_{n-1}\}$.
From (1) it follows that $x\notin V(T)$, and so $x\in \{z_1\l z_{\rho-1}\}$.
Since $d_G(v,p_j)=\rho$, it follows that $d_G(v,x)=\rho-1$, and so $x=z_{\rho-1}$. Consequently $d_G(p_j,p_n)\le 2$,
and so $j>n-\rho$ from the final condition in the definition of a $(v,T)$-extension. Since $d_G(p_{n-\rho},p_n)\ge \rho$
from the same condition, it follows that the path $p_{n-\rho}\d p_{n-\rho+1}\c p_j\d z_{\rho-1}\d p_n$
has length at least $\rho$, and so $j\ge n-2$. Now $j\ne n-1$ since $G$ is triangle-free, and $j\ne n$ by its definition,
so $j=n-2$. 

Consequently there is at most one edge joining the interiors of $P,Q$, and any such edge is between $z_{\rho-1}$ and $p_{n-2}$.
Let $s=|E(Q)|$ if there is no such edge, and $|E(Q)|-2$ if there is such an edge. In either case $G$ has a hole
of length $|E(P)|+s$. Moreover, since the final three vertices of $P_1\l P_k$ are the same, it follows that
$G$ has a hole of length $|E(P_i)|+s$ for $1\le i\le k$.
This proves (2).

\bigskip

Since $\chi(G)>\mu_0(\kappa)$, there exists $z_0$ such that $\chi(N^{\rho}_G[z_0]) >\mu_1(\kappa)+\kappa$, and hence 
$\chi(N^{\rho}_G(z_0)) >\mu_1(\kappa)$. 
Let $H_0 = G$ and let $T_0$ be the one-vertex 
subgraph with vertex $z_0$. For $1\le h\le \alpha+2$, we define 
$y_h,y'_h, S_h, z_h, T_h, M_h, H_h$ as follows. Assume we have defined $H_{h-1}$, $T_{h-1}$ and $z_{h-1}$ such that 
$\chi(N^{\rho}_{H_{h-1}}(z_{h-1}))>\mu_{h-1}(\kappa)$ and $T_{h-1}$ is an induced path of $G$ with at most $\rho+1$ vertices and with one end $z_{h-1}$.
Let $M_h$ be the subgraph induced on the set of all vertices $v$ of $H_{h-1}$ that satisfy
\begin{itemize}
\item $d_{H_{h-1}}(z_{h-1},v)= \rho$; and
\item $d_{G}(x,v)\ge \rho$ for every vertex $x$ of $T_{h-1}$.
\end{itemize}
Since $\chi(N^{\rho-1}[x])\le \kappa$  for each vertex $x$ of $T_{h-1}$, and $\chi(N^{\rho}_{H_{h-1}}(z_{h-1}))>\mu_{h-1}(\kappa)$, 
it follows
that 
$$\chi(M_h)>\mu_{h-1}(\kappa)-(\rho+1)\kappa=\phi(\phi(\mu_h(\kappa)+\kappa) +(2\rho+2)\kappa).$$ 
Since $G$ is $(\rho,\phi)$-controlled, there is a 
vertex $y_h\in M_h$ such that 
$$\chi(N^{\rho}_{M_h}[y_h])>\phi(\mu_h(\kappa)+\kappa) +(2\rho+2)\kappa,$$ 
and hence with
$$\chi(N^{\rho}_{M_h}(y_h))>\phi(\mu_h(\kappa)+\kappa) +(2\rho+1)\kappa.$$
Let $S_h$ be a path of $H_{h-1}$ of length $\rho$ between $z_{h-1}$ and $y_h$. Let $y'_h$ be adjacent to $y_h$ in $M_h$.
Let $S'_h$  be a path of $H_{h-1}$ of length $\rho$ between $z_{h-1}$ and $y'_h$. 
Let $H_h$ be the subgraph induced on the set of all vertices $v$ of $M_h$ with the following properties:
\begin{itemize}
\item $d_{M_h}(y_h,v)=\rho$; and
\item $d_G(x,v)\ge \rho$ for every $x\in V(S_h)\cup V(S_h')$.
\end{itemize}
Since $\chi(N^{\rho}_{M_h}(y_h)))>\phi(\mu_h(\kappa)+\kappa) +(2\rho+1)\kappa$, and 
$\chi(N^{\rho-1}[x])\le \kappa$  for each vertex $x$ of $V(S_h\cup S_h')$, and there are
at most $2\rho+1$ such vertices $x$, it follows that $\chi(H_h)>\phi(\mu_h(\kappa)+\kappa)$. Consequently there exists $z_h\in H_h$
such that $\chi(N^{\rho}_{H_h}[z_h])>\mu_h(\kappa)+\kappa$, and hence with $\chi(N^{\rho}_{H_h}(z_h))>\mu_h(\kappa)$. 
Let $T_h$ be a path of $M_h$ of length $\rho$ between $y_h, z_h$. This completes
the inductive definition of $y_h,y'_h, S_h, z_h, T_h, M_h, H_h$ for $1\le h\le \alpha+2$. 
\\
\\
(3) {\em For $1\le h\le \alpha+2$, $S_h\cup T_h$ is an induced path $L_h$ between $z_{h-1}, z_h$ of length $2\rho$. Also
there is an induced path $L_h'$ between $z_{h-1}, z_h$ with $V(L_h')\subseteq V(S_h'\cup T_h)$ of length $2\rho-1$ or $2\rho+1$.}
\\
\\
The first claim follows from (1). For the second, the graph formed by the union of $S_h'$, $T_h$ and 
the edge $y_hy'_h$ is a path, but it might not be induced.
If it is induced, it has length $2\rho+1$ as required; and since $S_h'$ and $T_h$ are both induced paths, we may assume
that some vertex $a$ of $S'_h$ is adjacent to some vertex $b$ of $T_h$, where $(a,b)\ne (y'_h, y_h)$.
Since every vertex of $S_h'$ has distance at most $\rho-2$ from $z_{h-1}$ except the last two, and every vertex
of $T_h$ has distance at least $\rho$ from $z_{h-1}$, it follows that $a$ is either $y'_h$ or its neighbour in $S'_h$.
Now $d_G(y_h', z_h) = \rho$, so $y'_h$ has no neighbour in $T_h$ except for $y_h$ (because $y'_h$ is not adjacent to
the second vertex of $T_h$ since $G$ is triangle-free). Thus $a$ is the penultimate vertex of $S'_h$.
Consequently $b\ne y_h$ since $G$ is triangle-free, and since $d_G(a,z_h)\ge \rho$, $a$ has no neighbour in $T_h$
different from the second vertex of $T_h$. We deduce that $b$ is indeed the second vertex of $T_h$; and so
there is an induced path between $z_{h-1}, z_h$ of length $2\rho-1$ with vertex set
a subset of $V(S_h'\cup T_h)$. This proves (3).

\bigskip

Let there be $q$ values of $h\in \{4\l \alpha+2\}$ such that $L_h'$ has length $2\rho-1$.
For $4\le h\le \alpha+2$, choose $L_h''\in \{L_h, L_h'\}$; then $L_4''\cup L_5''\cup\cdots\cup L_{\alpha+2}''$ is an induced
path between $z_3$ and $z_{\alpha+2}$, and it is a $(y_3,T_3)$-extension, for every choice of $L_4'',L_5''\l L_{\alpha+2}''$. 
Moreover, all these $(y_3,T_3)$-extensions
are parallel (since the last $\rho$ vertices of $L_{\alpha+2}, L'_{\alpha+2}$ are the same). These paths
have lengths every integer between $2\rho(\alpha-1) -q$ and $(2\rho+1)(\alpha-1)-q$, that is, every integer between
$\ell-\beta-q-2\rho$ and $\ell +\alpha -\beta -q -2\rho-1 $.
From (2), $G$ has holes of every length between
$\ell-\beta-q$ and $\ell+\alpha-\beta -q-5$. Since $G$ has no $l$-hole,
it follows that $\ell+\alpha-\beta -q-5<\ell$, that is, $\alpha\le \beta+q+4$. But by concatenating each of the paths
$L_4''\cup L_5''\cup\cdots\cup L_{\alpha+2}''$ with $L_3$, we obtain a $(y_2,T_2)$-extension of length exactly
$2\rho$ more; and so there are $(y_2,T_2)$-extensions of all lengths between $\ell-\beta-q$ and $\ell+\alpha-\beta -q -1$.
Hence by (2) there are holes in $G$ of all lengths between $\ell-\beta -q+2\rho$ and $\ell+\alpha-\beta -q+2\rho-5$. 
Since $\beta+q\ge \alpha-4\ge 2\rho$, it follows that
$\ell-\beta -q+2\rho\le \ell$. Consequently $\ell+\alpha-\beta -q+2\rho-5<\ell$, since there is no $\ell$-hole, 
that is,  $\alpha+2\rho\le\beta+q+4$. Similarly, by concatenating all these $(y_2,T_2)$-extensions with $L_2$, we obtain
$(y_1,T_1)$-extensions of all lengths between $\ell-\beta-q +2\rho$ and $\ell+\alpha-\beta -q +2\rho-1$. By (2), there are holes of
all lengths between $\ell-\beta-q +4\rho$ and $\ell+\alpha -\beta -q +4\rho -5$. But $\ell-\beta-q +4\rho\le \ell$, since
$\beta+q \ge \alpha+2\rho-4\ge 4\rho$, and yet
$$\ell+\alpha-\beta -q +4\rho -5=\ell+2\rho-3+ (\alpha-1-q)+(2\rho-1-\beta) \ge \ell$$
since $q\le \alpha-1$ and $\beta\le 2\rho-1$. Consequently there is an $\ell$-hole, a contradiction.
This proves
\ref{increaserad} and hence \ref{boundedrad}.~\bbox

\section{Showers}

Now we come to the third and most complicated part of the proof: proving \ref{bigrad}. This will occupy the remainder of the paper.

What can we prove about hole lengths
if $\chi^{\rho}(G)$ is bounded for some large fixed $\rho$?
In \ref{boundedrad} we were able to guarantee the presence of a hole of any desired length (almost), but in these new 
circumstances that becomes impossible; for any fixed $\rho\ge 0$ and $\ell\ge 2$, there are graphs with arbitrarily large $\chi$, and 
girth more than $\max(\ell,\rho/2)$;
which implies that $\chi^{\rho}(G)$ is at most $2$, and yet they have no $\ell$-hole.
We will show the following, a reformulation of \ref{bigrad}.

\begin{thm}\label{localbound}
Let $\nu\ge 2$ and $\kappa\ge 0$ be integers, and let $G$ be a triangle-free graph such that
$\chi^{\rho}(G)\le \kappa$, where $\rho = 3^{\nu+2}+4$.
If $G$ admits no hole $\nu$-interval then $\chi(G)$ is bounded.
\end{thm}

The proof will need a number of steps and preliminary lemmas. We begin with some definitions.
A {\em levelling} in $G$ is a sequence of pairwise disjoint subsets $(L_0, L_1\l L_k)$ of $V(G)$
such that
\begin{itemize}
\item $|L_0|=1$;
\item for $1\le i\le k$ every vertex in $L_i$ has a neighbour in $L_{i-1}$;
\item for $0\le i<j\le k$, if $j>i+1$ then no vertex in $L_j$ has a neighbour in $L_i$.
\end{itemize}
We call $L_k$ the {\em base} of the levelling.
The {\em chromatic number} of a levelling is the chromatic number of its base.
We observe first:
\begin{thm}\label{biglevel}
For any integer $\tau\ge 0$, if $\chi(G)>2\tau$ then $G$ admits a levelling with chromatic number more than $\tau$.
\end{thm}
\Proof
Choose a component $C$ of $G$ with chromatic number equal to that of $G$, and let $z$ be a vertex in that component.
For each $i\ge 0$, let $L_i$ be the set of vertices $v$ of $C$ such that $d_C(z,v)=i$, and choose $j$ such that
$L_0\cup\cdots\cup L_j = V(C)$. If $\chi(L_k)\le \tau$ for all $k$ with $0\le k\le j$, then $\chi(C)\le 2\tau$ 
(take two disjoint sets of colours
both of size $\tau$, and use them for the even and odd levels alternately), which is impossible; so there exists
$k$ such that $\chi(L_k)> \tau$. Then $(L_0\l L_k)$ is the desired levelling. This proves \ref{biglevel}.~\bbox

If $(L_0\l L_k)$ is a levelling in $G$, we call the unique vertex in $L_0$ the {\em head} of the levelling, and we call $L_0\cup\cdots\cup L_k$
the {\em vertex set} of the levelling. A path $P$ of $G[V]$ (where $V$ is the vertex set of the levelling) with ends $u,v$ is 
{\em monotone} (with respect to the given levelling) if
there exist $h,j$ with $0\le h,j\le k$, such that $u\in L_h, v\in L_j$, and $P$ has length $|j-h|$; and therefore $P$ has exactly one 
vertex in $L_i$ for each $i$ between $h,j$, and has no other vertices. 

There is a notational problem with levellings: that while it seems most natural to number levels starting with the head
as level zero, most of the action will be at or close to the base $L_k$, and we constantly have to refer to the parameter $k$.
To obviate this, let us say a vertex $v$ of the vertex set has  {\em height} $k-i$ if $v\in L_i$ where $0\le i\le k$.
Thus vertices in $L_k$ have height zero.

We have shown that, if we start with a triangle-free graph of large $\chi$, we can choose a levelling in it with base of large
$\chi$; and by replacing the base by one of its components with maximum chromatic number, we could choose the levelling 
such that the base is connected. This, however, is awkward to maintain, and not really necessary. 
All we really need is that the base has large $\chi$,
and is included in a connected set which has no further neighbours in higher parts of the levelling. So we will modify the
definition of a levelling to allow this. In addition, our main strategy to find a hole $\nu$-sequence is to fix some vertex 
in the base, which is joined to the head by a ``recirculator'' (a private path whose internal vertices have no neighbours elsewhere
in the levelling), and find holes of many different lengths all containing this recirculator; that is, we want to find
many paths of different lengths between the head of the shower and some fixed vertex of the base. Those two considerations
motivate the following definition.

A {\em shower} in $G$ is a sequence $(L_0, L_1\l L_k,s)$ where $L_0, L_1\l L_k$ are pairwise disjoint subsets
of $V(G)$ and $s\in L_k$, 
such that
\begin{itemize}
\item $|L_0|=1$;
\item for $1\le i< k$ every vertex in $L_i$ has a neighbour in $L_{i-1}$;
\item for $0\le i<j\le k$, if $j>i+1$ then no vertex in $L_j$ has a neighbour in $L_i$; and
\item $G[L_k]$ is connected.
\end{itemize}
(We suggest that the reader picture a shower with $L_0$ on top and $L_k$ at the bottom, in order to make
sense of the terminology to come.)
The differences between a shower and a levelling are that, first, not every vertex in $L_k$ needs to have a neighbour in $L_{k-1}$
(and indeed, there may be no edges between $L_{k-1}$ and $L_k$, although such showers will not be of interest);
second, that $G[L_k]$ is connected; and third, the distinguished vertex $s$. We call $L_0\l L_k$ the {\em levels} of the shower, and
$s$ the {\em drain} of the shower. We define ``head'',
``base'', ``vertex set'', ``monotone'', ``height'' for showers just as for levellings. The set of vertices in $L_k$ with a neighbour in $L_{k-1}$
is called the {\em floor} of the shower. (It is the floor, and subsets of the floor, whose chromatic number will
concern us.) If $\mathcal{S}=(L_0, L_1\l L_k,s)$ is a shower,
and $uv$ is an edge with $u \in L_i$ and $v\in L_{i+1}$ for some $i$ with $0\le i<k$, we say that $u$ is an {\em $\mathcal{S}$-parent} or just {\em parent}
of $v$, and $v$ an {\em $\mathcal{S}$-child} or just {\em child} of $u$.

If $\mathcal{S}=(L_0\l L_k,s)$ is a shower, with head $z_0$ and vertex set $V$,
a {\em recirculator} for $\mathcal{S}$
is an induced path $R$
with ends $s,z_0$ such that no internal vertex of $R$ belongs to $V$
and no internal vertex of $R$ has any neighbours in $V\setminus \{s,z_0\}$. 
The {\em distance} $d_G(P_1,P_2)$ between two nonnull subgraphs $P_1,P_2$ of $G$ 
is the minimum of $d_G(v_1,v_2)$ over all $v_1\in V(P_1)$ 
and $v_2\in V(P_2)$. 

\begin{thm}\label{doubleshower}
Let $\tau, \kappa\ge 0$ be integers.
Let $G$ be a graph such that $\chi^{8}(G)\le \kappa$.
Let $(L_0\l L_k)$ be a levelling in $G$, where $\chi(L_k)>22\tau+2\kappa$.
Then there is a shower $(V_0\l V_n,s)$ in $G$, with floor of chromatic number more than $\tau$, and with a recirculator, such that 
\begin{itemize}
\item $V_n\subseteq L_k$, and $V_{n-1}\subseteq L_{k-1}$; and
\item $V_0\l V_{n-2}\subseteq L_0\cup\cdots\cup L_{k-2}$.
\end{itemize}
\end{thm}
\Proof
By replacing $L_k$ by the vertex set of a component of $G[L_k]$ with maximum chromatic number, 
we may assume that $G[L_k]$ is connected.
A {\em stake} is a monotone path with an end in $L_k$.
Since $\chi(L_k)>\kappa $, there exist two vertices of $L_k$ with distance more than $8$.
It follows that there are two stakes both of length three with distance at least three. Consequently we can
choose two stakes $P,Q$ with the following properties:
\begin{itemize}
\item $P$, $Q$ have the same length $k-h\ge 3$; 
\item $d_G(P,Q)\ge 3$;
\item subject to these two conditions, $h$ is minimum.
\end{itemize}
Let $P$ have vertices $p_k\d p_{k-1}\c p_{h}$ and $Q$ have vertices $q_k\d q_{k-1}\c q_{h}$, where $p_i,q_i\in L_i$
for $h\le i\le k$. Let $p_{h-1}, q_{h-1}$ be parents of $p_{h}, q_{h}$ respectively.
From the minimality of $h$,
either 
\begin{itemize}
\item $p_{h-1}, q_{h-1}$ are adjacent, or
\item some vertex is adjacent to $p_{h-1}$ and to at least one of $q_{h-1},q_{h},q_{h+1}$, or
\item some vertex is adjacent to $q_{h-1}$ and to at least one of $p_{h-1},p_{h},p_{h+1}$.
\end{itemize}
In each case there is a connected induced subgraph $M$ with $V(M)\subseteq L_0\cup\cdots\cup L_h\cup \{p_{h+1},q_{h+1}\}$, 
with at most seven vertices,
and with $p_{h+1}, p_h,p_{h-1},q_{h+1}, q_h,q_{h-1}\in V(M)$; and if there is a vertex in $V(M)\setminus V(P\cup Q)$, 
then it belongs to $L_{h-2}\cup L_{h-1}\cup L_h$, 
and has a neighbour in $\{p_{h+1}, p_h,p_{h-1}\}$ and one in 
$\{q_{h+1}, q_h,q_{h-1}\}$. Consequently, $p_{h+2},\ldots,p_k$ have no neighbours 
in $V(M)\setminus \{p_{h+1}\}$, and $q_{h+2},\ldots,q_k$ have no neighbours 
in $V(M)\setminus \{q_{h+1}\}$.

Let $X$ be the set of vertices $x\in L_{k-1}$ such that there is a path $R$ from $x$ to $p_{h+1}$ satisfying:
\begin{itemize}
\item $R$ has length at most $k-h+8$;
\item every internal vertex of $R$ belongs to $L_0\cup \cdots\cup L_{k-2}$; and
\item no vertex of $R\setminus p_{h+1}$ equals or is adjacent to any vertex in $\{p_{h+2}\l p_k\}$.
\end{itemize}
Define $Y\subseteq L_{k-1}$ similarly with $P,Q$ exchanged.
\\
\\
(1) {\em Every vertex $v \in L_k$ with $d_G(v,p_k), d_G(v,q_k)\ge 7$ has a neighbour in $X\cup Y$.}
\\
\\
Let $v\in L_k$ with $d_G(v,p_k), d_G(v,q_k)\ge 7$, and let $r_0\d r_1\c r_k = v$ be a path between $r_0\in L_0$ and $v=r_k$. 
We claim that $r_{k-1}\in X\cup Y$.
From the minimality of $h$,
one of $r_{h-1}\l r_k$ has distance at most two from one of $p_{h-1}\l p_k$. Choose $j$ maximum such that 
$r_j$ has distance at most two from some vertex $u$ say of $P\cup Q\cup M$. Thus $j\ge h-1$. If $j=k$, then
$u\notin V(M)\setminus V(P\cup Q)$ because $k-h\ge 3$, and so $u$ is one of $p_k,p_{k-1},p_{k-2},q_k,q_{k-1},q_{k-2}$; which is impossible
since $d_G(v,p_k),d_G(v,q_k)\ge 7$. Thus $j<k$. From the maximality of $j$, it follows that $d_G(r_j,u)=2$,
and none of $r_j\l r_k$ equals or is adjacent to any vertex in $P\cup Q\cup M$. From the symmetry we may assume that
$u\in V(Q)\cup (V(M)\setminus V(P\cup Q))$. Let $w$ be a vertex adjacent to both $u,r_j$. If $u\in L_k\cup L_{k-1}$
then $k-j\le 3$, and so $d_G(v,q_k)\le 6$, a contradiction; and if $u\notin L_k\cup L_{k-1}$ and $w\in  L_k\cup L_{k-1}$
then $u=q_{k-2}$ and $k-j\le 2$, and again $d_G(v,q_k)\le 6$, a contradiction. So $u,w\notin L_k\cup L_{k-1}$.
If $w$ has a neighbour in $\{p_{h+2}\l p_k\}$, then $w\in L_{h+1}\cup\cdots\cup L_k$, and so $u\in V(Q)$, contradicting that $d_G(P,Q)\ge 3$. Thus $w$ has
no neighbour in $\{p_{h+2}\l p_k\}$.

Now there is a path of $M\cup Q$ between $u$ and $p_{h+1}$. 
If $u\notin V(Q)$ then this path has length at most three, and its union with the path $r_{k-1}\d r_{k-2}\c r_j\d w\d u$ is of length 
at most $k-1-j+5\le k-h+5$, since $j\ge h-1$, and so $r_{k-1}\in X$ as required.
If $u\in V(Q)$, then $u$ is one of $q_{j-2}, q_{j-1}, q_j, q_{j+1}, q_{j+2}$,
and so some path of $M\cup Q$ between $u$ and $p_{h+1}$ has length at most $(j+2)-(h+1)+6$, and its union with the
path  $r_{k-1}\d r_{k-2}\c r_j\d w\d u$ has length at most 
$$(j+2)-(h+1)+6 + (k-1-j)+2=k-h+8,$$
and again $r_{k-1}\in X$.
This proves (1).

\bigskip

Now, since $\chi^8(G)\le \kappa$, the set of vertices $v\in L_k$ such that $d_G(v,p_k)\le 6$ or $d_G(v,q_k)\le 6$
has chromatic number at most $2\kappa$; and since $\chi(L_k)> 22\tau+2\kappa$, there exists a subset
$Z_0\subseteq L_k$ with $\chi(Z_0)> 22\tau$ such that $d_G(v,p_k), d_G(v,q_k)\ge 7$ for each $v\in Z_0$.
Every vertex in $Z_0$ has a neighbour in $X\cup Y$, by (1); so we may assume that there exists $Z_1\subseteq Z_0$
with $\chi(Z_1)>11\tau$, such that every vertex in $Z_1$ is adjacent to a vertex in $X$. For each vertex $x\in X$,
there is a path $R$ as in the definition of $X$; let $R_x$ be a shortest such path. Then $R_x$ has length at
most $k-h+8$, and at least $(k-1)-(h+1)$; so there are eleven possibilities for its length, the numbers between
$k-h-2$ and $k-h+8$. For each $c$ with $k-h-2\le c\le k-h+8$, let $X_c$ be the set of vertices $x\in X$ such that $R_x$
has length $c$. Then there exist $c$ and 
$Z_2\subseteq Z_1$ with $\chi(Z_2)\ge \chi(Z_1)/11> \tau$, such that
every vertex in $Z_2$ has a neighbour in $X_c$. Moreover we may choose $Z_2$ such that $G[Z_2]$ is connected.
Let $V$ be the union of the vertex sets of all the paths $R_x\;(x\in X_c)$.
Note that $V\subseteq L_0\cup \cdots\cup L_{k-1}$. For $0\le i\le c$, let $V_i$ be the set of vertices $u\in V$ such that
the shortest path of $G[V]$ between $u,p_{h+1}$ has length $i$. Then $(V_0\l V_c)$ is a levelling. Moreover,
$V_c=X_c$, and so no vertex in
$L_k$ has a neighbour in $V_0\l V_{c-1}$. Define $V_{c+1} = Z_2$; then also $(V_0\l V_{c+1})$ is a levelling.

Now no neighbour of $p_{k-1}$ belongs to $Z_0$, and hence there are no edges between $\{p_{h+2}\l p_{k-1}\}$
and $V_1\cup \cdots\cup V_{c+1}$. Since $G[L_k]$ is connected and $p_{k-1}$ has a neighbour in $L_k$, there is a path
$G[L_k]$ between a vertex adjacent to $p_{k-1}$ and a vertex with a neighbour in $Z_2=V_{c+1}$.
Choose a minimal such path, $D$, and let $s$ be its end adjacent to $p_{k-1}$. 
Then $(V_0\l V_c, V_{c+1}\cup V(D), s)$ is a shower, since $G[Z_2]$ is connected and hence so is $G[V_{c+1}\cup V(D)]$;
and its floor includes $Z_2$ and hence has chromatic number more than $\tau$;
and $p_{h+1}\d p_{h+2}\c p_{k-1}\d s$ is a recirculator for it.
This proves \ref{doubleshower}.~\bbox

\bigskip

Let $\mathcal{S}$ be a shower with head $z_0$, drain $s$ and vertex set $V$.
An induced path of $G[V]$ between $z_0,s$ is called a {\em jet} 
of $\mathcal{S}$.
The set of all lengths of jets of $\mathcal{S}$ is called the {\em jetset} of $\mathcal{S}$.
If $\mathcal{A}$
is a subset of the jetset of $\mathcal{S}$, then for each $a\in \mathcal{A}$ there is a jet $J_a$ with length $a$, and we say the
set of jets $\{J_a:\;a\in \mathcal{A}\}$ {\em realizes} $\mathcal{A}$.
For $\nu\ge 2$, we say a shower $\mathcal{S}$ is {\em $\nu$-complete} if there are $\nu$ consecutive integers in its jetset, and
{\em $\nu$-incomplete} otherwise. (Later we shall give a meaning to ``$1$-complete'', but at this stage it is not needed.)
We deduce:

\begin{thm}\label{firstshower}
Let $\tau, \kappa\ge 0$ and $\nu\ge 2$ be integers.
Let $G$ be a graph such that 
\begin{itemize}
\item $\chi^{8}(G)\le \kappa$;
\item $\chi(G)>44\tau+4\kappa$; and
\item $G$ admits no hole $\nu$-interval.
\end{itemize}
Then there is a $\nu$-incomplete shower in $G$ with floor of chromatic number more than $\tau$.
\end{thm}
\Proof
By \ref{biglevel} there is a levelling $(L_0\l L_k)$ with chromatic number more than $22\tau+2\kappa$.
By \ref{doubleshower},  there is a shower $\mathcal{S}$, with a recirculator, and
with floor of chromatic number more than $\tau$. Since the union of the recirculator with any jet is a hole, and
$G$ admits no hole $\nu$-interval, it follows that $\mathcal{S}$ is not $\nu$-complete. This proves \ref{firstshower}.~\bbox

Thus, in order to prove \ref{localbound}, it suffices to show that if $\nu,\kappa$, $G$ are as in the hypothesis of \ref{localbound} 
then the floor of every $\nu$-incomplete shower in $G$ has bounded chromatic number, and this is what we shall do.

\section{Stabilizing a shower}

A levelling  $(L_0\l L_k)$ or shower $(L_0\l L_k,s)$ is {\em stable} if $L_0\l L_{k-1}$ are stable; and
for $\lambda\ge 0$ an integer, it is {\em $\lambda$-stable} if $k\ge \lambda$ and $L_i$ is stable
for $k-\lambda\le i\le k-1$.
We would like to prove that there exists a stable shower
(still with floor of large $\chi$, but not as large as before), by converting the shower given by \ref{firstshower}. 
This will take several steps.
First we show how to convert a $\nu$-incomplete shower into a $\nu$-incomplete 
$\lambda$-stable shower (for any fixed $\lambda$).

If $\mathcal{S}$ is a levelling $(L_0\l L_k)$ or a shower $(L_0\l L_k,s)$, and 
there is a monotone path $P$ with ends $u,v$, and $u\in L_i$ and $v\in L_j$ where $j\ge i$, we say that
$v$ is a {\em $\mathcal{S}$-descendant} (or just {\em descendant}) of $u$ and $u$ is an {\em $\mathcal{S}$-ancestor} 
(or just {\em ancestor}) of $v$. If $X\subseteq L_0\cup\cdots\cup L_k$,
we denote by 
$\theta(X)$ or $\theta_{\mathcal{S}}(X)$ the chromatic number of the set of vertices
in $L_k$ with an ancestor in $X$.

\begin{thm}\label{basestable}
Let $\tau, \lambda \ge 0$ and $\nu\ge 2$ be integers, and let $\mu = (\lambda+1)(\nu-1)+1$.
Let $G$ be a triangle-free graph, and let $\mathcal{S}$ be
a $\nu$-incomplete shower in $G$, with floor of chromatic number more than $\nu \tau^{1+\mu}$, and 
 with levels $L_0\l L_k$, where $k\ge \mu$.
Then there is a $\lambda$-stable $\nu$-incomplete shower 
with floor of chromatic number more than $\tau$, and with levels $L_0'\l L_h'$, such that $0\le k-h\le \mu-\lambda-1$ and 
$L_i'\subseteq L_i$ for $0\le i< h$.
\end{thm}
\Proof
We may assume that for $0\le i<k$, every vertex in $L_i$ has a neighbour in $L_{i+1}$; for a vertex in $L_i$ without this
property could be deleted.
Let $z_0\in L_0$.  For $1\le j\le \nu$, let $h_j= k-1-(\lambda+1)(\nu-j)$; and for $1\le j< \nu$, let $I_j=\{i:h_{j}<i<h_{j+1}\}$.
(Thus the sets $I_j$ have cardinality $\lambda$, and there is an integer $h_j$
between $I_{j-1}$ and $I_{j}$ that belongs to neither, that we use as insulation.)
For $1\le j\le \nu$, let $T_j$ be the set of vertices $v\in L_{h_j}$ such that there are $j$ induced paths
between $v$ and $z_0$, each with interior in $L_1\cup \cdots\cup L_{h_j-1}$, of lengths $h_j, h_j+1\l h_j+j-1$.
\\
\\
(1) {\em $T_{\nu} = \emptyset$.}
\\
\\
Because suppose that $v\in T_{\nu}$. Then there are $\nu$ induced paths
between $v$ and $z_0$, each with interior in $L_1\cup \cdots\cup L_{k-2}$, of lengths $k-1,k\l k+\nu-2$,
say $R_1\l R_{\nu}$.
Let $s$ be the drain of $\mathcal{S}$; and choose a minimal path $Q$ between $s,v$ with interior in $L_k$.
Then for $1\le i\le \nu$, the union of $Q$ and $R_i$ is a jet, contradicting that the shower is $\nu$-incomplete.
This proves (1).

\bigskip
Since $T_1 = L_{h_1}$ it follows that
$$\theta(T_1)>
\nu \tau^{1+\mu}\ge
\tau^{k+1-h_2}+\tau^{k+1-h_{3}}+\cdots+ \tau^{k+1-h_{\nu}},$$
and so there exists $j\in \{1\l \nu\}$ maximum such that
$$\theta(T_j)>\tau^{k+1-h_{j+1}}+\tau^{k+1-h_{j+2}}+\cdots+ \tau^{k+1-h_{\nu}};$$
and $j<\nu$ by (1).
From the maximality of $j$ it follows that
$\theta(T_j)-\theta(T_{j+1})> \tau^{k+1-h_{j+1}}$.
Let $S_{j+1}$ be the set of vertices in $L_{h_{j+1}}\setminus T_{j+1}$
that have ancestors in $T_j$. For $h_j<i<h_{j+1}$ let $M_i$ be the set of vertices in $L_i$ with an ancestor in $T_j$
and a descendant in $S_{j+1}$.
\\
\\
(2) {\em $M_i$ is stable for $h_j<i<h_{j+1}$.}
\\
\\
For suppose that $x,y\in M_i$ are adjacent. Since $G$ is triangle-free, $x,y$ have no common parents and no common children.
Let $x',y'\in T_j$ be ancestors of $x,y$ respectively (possibly equal).
Let $z\in S_{j+1}$ be a descendant of $x$. Now there are induced paths from $y'$ to $z_0$
with interior in $L_1\cup \cdots\cup L_{h_j-1}$, of lengths $h_j, h_j+1\l h_j+j-1$. For each of these paths, its union
with a path of length $i-h_j$ between $y$ and $y'$, a path of length $h_{j+1}-i$ between $z$ and $x$, and the edge $xy$,
makes an induced path between $z,z_0$, of lengths $h_{j+1}+1\l h_{j+1}+j$. But also there is an induced path between $z, z_0$ of length
$h_{j+1}$, since $z\in L_{h_{j+1}}$; and so $z\in T_{j+1}$, a contradiction. This proves (2).

\bigskip
Now every vertex in $L_k$ with an ancestor in $T_j$ has an ancestor in $S_{j+1}\cup T_{j+1}$. Since
$\theta(T_j)-\theta(T_{j+1})> \tau^{k+1-h_{j+1}}$,
it follows that
$\theta(S_{j+1})> \tau^{k+1-h_{j+1}}$.
By setting $h=h_{j+1}$ and $M_h=S_{j+1}$, we have shown that:
\\
\\
(3) {\em There exist $h$ with $0 \le k-h\le \mu-\lambda-1$, and subsets $M_i\subseteq L_i$ for $h-\lambda\le i\le h$,
 with the following properties:
\begin{itemize}
\item $\theta(M_h)> \tau^{k+1-h}$;
\item $M_i$ is stable for $h-\lambda\le i< h$; and
\item every vertex in $M_{i+1}$ has a neighbour in $M_{i}$ for $h-\lambda\le i< h$.
\end{itemize}
}
Choose such a value of $h$, maximal. Suppose first that $\chi(M_h)\le \tau$.
Since 
$$\theta(M_h)>\tau^{k-h+1}\ge \tau\ge \chi(M_h)$$ 
it follows that $h\ne k$.
Take a partition of $M_h$
into $\tau$ stable sets; then for one of these sets, say $M_{h}'$, $\theta(M_h')\ge \theta(M_h)/\tau> \tau^{k-h}$.
Let $M_{h+1}$ be the set of vertices in $L_{h+1}$ with a neighbour in $M_h$; then $\theta(M_{h+1})=\theta(M_h')>\tau^{k-h}$,
contrary to the maximality of $h$. This proves that $\chi(M_h)> \tau$.

Let $Z=L_h\cup \cdots\cup L_k$; then $G[Z]$ is connected since $G[L_k]$ is connected and
for $0\le i<k$, every vertex in $L_i$ has a neighbour in $L_{i+1}$. Consequently
$$(L_0\l L_{h-\lambda-1}, M_{h-\lambda}\l M_{h-1}, Z,s)$$
is a shower $\mathcal{S}'$ say.
Its floor includes $M_h$ and so has chromatic number more than $\tau$.
Moreover, every jet for $\mathcal{S}'$ is also a jet for $\mathcal{S}$; and so
$\mathcal{S}'$ is $\nu$-incomplete.
This proves \ref{basestable}.~\bbox

\section{U-bends}

For $\nu \ge 2$, a  shower $(L_0\l L_k,s)$ is a {\em $\nu$-sprinkler} if 
\begin{itemize}
\item $G[L_k]$ is a path with one end $s$ and with at least $\nu$ vertices; let its vertices be $v_1\c v_n$ in order, 
where $v_1 = s$ and $n\ge \nu$;
\item for $i = 1\l n-\nu$, no vertex in $L_{k-1}$ is adjacent to $v_i$; and
\item for $i = n-\nu+1\l n$, some vertex in $L_{k-1}$ is adjacent to $v_i$ and to no other vertex in $L_{k}$.
\end{itemize}
Every $\nu$-sprinkler is therefore $\nu$-complete. We call $\{v_i:n-\nu+1\le i\le n$\} its {\em floor}.

We need another object, a ``u-bend'', which is not exactly a shower; and also something which is partway to a u-bend, which 
we call a ``w-bend''. We start with the latter.
Let $(L_0\l L_k)$ be a levelling in $G$ with vertex set $V$,
and let $U$ be an induced path of $G$. Suppose that
\begin{itemize}
\item $G[L_k]$ is an induced path;
\item $V\cap V(U)=\emptyset$;
\item $U$ has ends $w,s$, and there is at least one vertex in $L_{k-1}$ adjacent to $w$ and to a vertex in $L_k$; and
\item there are no edges between $V(U)$ and $V\setminus L_{k-1}$, and no vertex 
in $L_{k-1}$ has a neighbour in $L_k$ and a neighbour 
in $V(U)\setminus \{w\}$.
\end{itemize}
In this case, we call $(L_0\l L_k, U)$ a {\em w-bend}, and call $s$ its {\em drain}; and any induced path of $G[V\cup V(U)]$ between
the vertex in $L_0$ and the drain is called a {\em jet} of the w-bend. We call $L_k$ its {\em floor}. (Since $(L_0\l L_k)$
is a levelling, every vertex in $L_k$ has a neighbour in $L_{k-1}$.)
Let $G[L_k]$ have ends $v_1,v_2$; then $d_G(v_1,v_2)$ is called the {\em size} of the w-bend.
If in addition:
\begin{itemize}
\item $w$ has a unique neighbour in $L_{k-1}$, say $v$;
\item $v$ has a unique neighbour in $L_k$, and this neighbour is an end of the path $G[L_k]$; and
\item every vertex in $L_{k-1}$ has a neighbour in $L_k$;
\end{itemize}
then we call $(L_0\l L_k, U)$ a {\em u-bend}.
We need a containment relation for these objects:
\begin{itemize}
\item
Let $\mathcal{S}=(L_0\l L_k,s)$ and $\mathcal{S}'=(L_0'\l L_k',s')$ be showers.
We say that $\mathcal{S}'$ is {\em contained in} $\mathcal{S}$ if
they have the same drain, and $L_i'\subseteq L_i$ for $0\le i\le k$.
\item
Let $\mathcal{S}=(L_0\l L_k,s)$ be a shower, and let $\mathcal{S}'=(L_0'\l L_k',U)$ be a w-bend.
We say that $\mathcal{S}'$ is {\em contained in} $\mathcal{S}$ if
they have the same drain, and $L_i'\subseteq L_i$ for $0\le i\le k$, and $V(U)\subseteq L_k$.
\item
Let $\mathcal{S}=(L_0\l L_k,W)$ be a w-bend, and let $\mathcal{S}'=(L_0'\l L_k',U)$ be a u-bend.
We say that $\mathcal{S}'$ is {\em contained in} $\mathcal{S}$ if
they have the same drain, and $L_i'\subseteq L_i$ for $0\le i\le k$, and $U=W$.
\end{itemize}
In all three cases, every jet of $\mathcal{S}'$ is a jet
of $\mathcal{S}$.

We need to show that certain showers contain u-bends, and it is easier to show that they contain w-bends. Let us see first
that that is enough, because a w-bend contains a u-bend (and containment is clearly transitive).

\begin{thm}\label{wtou}
Let $(L_0\l L_k,W)$ be a w-bend in a triangle-free graph $G$, with size at least $2p+4$. Then it contains a u-bend with size at least $p$.
\end{thm}
\Proof
Let $W$ have ends $w,s$ where $s$ is the drain. 
Let $G[L_k]$ have vertices $v_0\c v_n$ say, in order.
Since $d_G(v_0,v_n)\ge 2p+4$, we may assume
by exchanging $v_0, v_n$ if necessary that $d_G(w,v_0)\ge p+2$.
Let $Y$ be the set of vertices in $L_{k-1}$ adjacent to $w$ and to a vertex in $L_k$. By hypothesis, $Y\ne \emptyset$.
Choose $i\le n$ minimum such that $v_i$ has a neighbour in $Y$, say $v$. Since $d_G(w,v_0)\ge p+2$, and $d_G(w,v_i)=2$,
it follows that $d_G(v_i,v_0)\ge p$. Let $L_{k-1}'$ consist of all vertices in $L_{k-1}$ with a neighbour in $\{v_0\l v_{i-1}\}$,
together with $v$. Then $v$ is the unique neighbour of $w$ in $L_{k-1}'$; and so 
$$(L_0\l L_{k-2}, L_{k-1}', \{v_0\l v_i\}, W)$$
is a u-bend contained in $(L_0\l L_k,W)$, and its size is at least $p$. This proves \ref{wtou}.~\bbox

\begin{thm}\label{getubend}
Let $\nu\ge 2$ be an integer, and let $\mu\ge 1$.
Let $\mathcal{S}$ be a shower in a triangle-free graph $G$.
Let $P$ be an induced path of $G$ with $V(P)$ a subset of the floor of $\mathcal{S}$, 
with ends $w_1,w_2$ such that $d_G(w_1,w_2)\ge 2(\mu+\nu)$.
Then $\mathcal{S}$ contains either:
\begin{itemize}
\item a $\nu$-sprinkler with floor a subset of $V(P)$, or
\item a u-bend with size at least $\mu$
and with floor a subset of $V(P)$.
\end{itemize}
\end{thm}
\Proof 
Let $\mathcal{S}=(L_0\l L_k,s)$, and let
$L_{k-1}^1$ be the set of vertices in $L_{k-1}$ with a neighbour in $V(P)$.
If $s\in V(P)$, let $u=s$ and let $D$ be the one-vertex path with vertex $s$.
If $s\notin V(P)$, then since 
$G[L_k]$ is connected, there is an induced path $D$ of $G[L_k]$ between $s$
and a vertex with a neighbour in $V(P)$; choose a minimal such path $D$, with ends
$s,u$ say. From the minimality of $D$, no vertex in 
$D\setminus \{u\}$ has a neighbour in $V(P)$.

Suppose that some vertex of $D\setminus \{u\}$ has a neighbour in $L_{k-1}^1$; and choose such a vertex, $w$
say, such that the subpath $D'$ of $D$ between $w,s$ is minimal. Then 
$$(L_0\l L_{k-2},L_{k-1}^1, V(P), D')$$
is a w-bend contained in $\mathcal{S}$, of size at least $2(\mu+2)$ (since $\nu\ge 2$), and the result follows from \ref{wtou}.
We may therefore assume that there are no edges between $D\setminus \{u\}$ and $L_{k-1}^1$.

Let $Y$ be the set of vertices in $L_{k-1}^1$ that are adjacent to $u$.
Now no vertex of $D$ except possibly $u$ has a neighbour in $L_{k-1}^1$; and 
$u$ has at least one neighbour in $V(P)\cup Y$.
Let $P$ have vertices $v_0\c v_n$ in order. By hypothesis, $d_G(v_0,v_n)\ge 2(\mu+\nu)$, so
by exchanging $v_0,v_n$ if necessary, we may assume that $d_G(u,v_0)\ge \mu+\nu$.
Choose $i$ minimum such that $v_i$ has a neighbour in $Y\cup \{u\}$. 

Suppose first that $v_i$ has a neighbour in $Y$. Choose such
a neighbour $v$ say, and let $L_{k-1}^2$ be the set of vertices in $L_{k-1}$ with a neighbour in $\{v_0\l v_{i-1}\}$, together with 
$v$. Now $v_i$ is not adjacent to $u$ (since $G$ is triangle-free); and $d_G(v_0,v_i)\ge d_G(v_0,u)-2\ge \mu$; so
$$(L_0\l L_{k-2}, L_{k-1}^2, \{v_0\l v_i\}, D)$$
is a u-bend contained in $\mathcal{S}$ with size at least $\mu$, as required.

We may assume then that $v_i$ has no neighbour in $Y$, and therefore $v_i$ is adjacent to $u$. 
In summary, no vertex in $L_{k-1}$
has a neighbour in $V(D)$ and a neighbour in $\{v_0\l v_i\}$; and 
there are no edges between $V(D)$ and $\{v_0\l v_i\}$ except the edge $uv_i$. 
Since
$d_G(v_0,u)\ge \mu+\nu$, it follows that $i\ge \mu+\nu-1$, and so $i-\nu+1\ge \mu$.

Suppose next that there exists
a vertex in $L_{k-1}$ adjacent to at least two of $v_{i-\nu+1}\l v_i$.
Choose $j$ with $i-\nu+3\le j\le i$ maximum such that some vertex in $L_{k-1}$ is adjacent to $v_j$ and to 
one of $v_0\l v_{j-2}$; choose $h$ with $0\le h\le j-2$ minimum such that some vertex in $L_{k-1}$ is adjacent to 
$v_h,v_j$; and choose $v\in L_{k-1}$ adjacent to $v_h, v_j$. Let $L_{k-1}^3$ be the set of vertices
in $L_{k-1}$ with a neighbour in $\{v_0\l v_{h-1}\}$, together with $v$. Then since there is a path between $u,v_h$
(via $v$) of length $i-j+3\le \nu$, it follows that $d_G(u,v_h)\le \nu$, and so 
$$d_G(v_h,v_0)\ge d_G(u,v_0)-\nu\ge \mu.$$
Let $D_2$ be the path formed by the union of $D$ and the path $u\d v_i\c v_j$.
Then 
$$(L_0\l L_{k-2}, L_{k-1}^3, \{v_0\l v_h\}, D_2)$$
is a u-bend contained in $\mathcal{S}$, of size at least $\mu$, as required.

We may therefore assume that no vertex in $L_{k-1}$ is adjacent to more than one of $v_{i-\nu+1}\l v_i$.
Let $L_{k-1}^4$ be the set of vertices in $L_{k-1}$ with a neighbour in  $ \{v_{i-\nu+1}\l v_i\}$.
Every vertex in $ \{v_{i-\nu+1}\l v_i\}$ has a neighbour in $L_{k-1}^4$, and $u$ has no neighbour in $L_{k-1}^4$, so
$$(L_0\l L_{k-2}, L_{k-1}^4, V(D)\cup \{v_{i-\nu+1}\l v_i\},s)$$
is a $\nu$-sprinkler contained in $\mathcal{S}$. This proves \ref{getubend}.~\bbox

\section{Jets of a shower}

Let $L_0\l L_k$ be the levels of a shower or w-bend, and let $J$ be a jet. Then at least one vertex of $J$ belongs to $L_{k-1}$;
and we define the {\em tail} of $J$ to be the minimal subpath of $J$ between $L_{k-1}$ and the drain.
For $\lambda\ge 0$,
we say that $J$ is {\em $\lambda$-monotone} if 
$\lambda\le k$, and 
$J$ contains exactly one vertex of $L_i$ for $0\le i< k-\lambda$. In every jet $J$, at least $k-1$ edges do not belong to its tail
and have an end not in $L_k$.
We say the {\em waste} of $J$ is $\mu$ if there are $k-1+\mu$ edges of $J$ that do not belong to its tail and have an end 
not in $L_k$; and $J$ is {\em $\mu$-wasteful} if its waste is at most $\mu$. 
Thus the waste is nonnegative.

A set of integers $\mathcal{A}$ is {\em dense} if for all $a_1,a_2\in \mathcal{A}$ with $a_1<a_2$, there does not
exist $b$ with $a_1<b<a_2$ such that $b,b+1\notin \mathcal{A}$; that is, there are no two consecutive numbers
both missing from $\mathcal{A}$ between the first and last members of $\mathcal{A}$.
If $\mathcal{A},\mathcal{B}$ are sets of integers, we define $\mathcal{A}+\mathcal{B} = \{a+b:\;a\in \mathcal{A}, b\in \mathcal{B}\}$.
Thus if $\mathcal{A}$ is dense, then for any integer $t$, $\mathcal{A}+\{t,t+1\}$ is a 
set of consecutive integers of cardinality at least $|\mathcal{A}|+1$.

Any subset of the floor of a shower is called a {\em mat}; and for a w-bend, we define its floor to be its only mat.
The {\em size} of a mat $M$ is the maximum of
$d_G(w_1,w_2)$ over all pairs $w_1,w_2$ of vertices in the same component of $G[M]$.
If $M$ is a mat for a shower or w-bend $\mathcal{S}$, a jet $J$ is
an {\em $M$-jet} if there is no edge of $J$ with an end in $L_k\setminus M$ and an end in $L_{k-1}$.  We define the {\em $M$-jetset}
as the set of all lengths of $M$-jets. A w-bend $(L_0\l L_k,U)$ is {\em $\lambda$-stable} if $k\ge \lambda$ and $L_i$ is stable
for $k-\lambda\le i\le k-1$.
In this section we prove the following. (Note that the next result immediately implies the long odd holes conjecture, via 
\ref{getubend}, so if we only wanted the long odd holes conjecture we could stop here.)

\begin{thm}\label{staircase}
Let $\nu\ge 2$ be an integer, and let $G$ be a triangle-free graph. If $\mathcal{S}$ is a $\nu$-stable shower or w-bend in $G$, 
and $M$ is a mat for $\mathcal{S}$ of size at least $3^{\nu+2}$, then there is a set $\mathcal{A}$ of integers, realized by
a set of $(\nu+1)$-monotone, $3\nu^2$-wasteful $M$-jets, such that 
$|\mathcal{A}|\le \nu+1$, and $\mathcal{A}$ includes a dense subset of cardinality $\nu$, and there are two members of $\mathcal{A}$ that differ by $1$ or $3$.
\end{thm}
\Proof
We proceed by induction on $\nu$. 
Thus we assume that either $\nu=2$ or 
the result holds for $\nu-1$. We claim we may assume:
\\
\\
(1) {\em There is a u-bend $\mathcal{S}_1 = (L_0\l L_k,U)$ contained in $\mathcal{S}$, and with $L_k\subseteq M$, 
of size at least $3^{\nu+2}/2-\nu$.}
\\
\\
Assume first that $\mathcal{S}$ is a $\nu$-stable shower in $G$, and $M$ is a mat of size at least $3^{\nu+2}$.
Let $P$ be an induced path of $G[M]$ with ends $w_1,w_2$, where $d_G(w_1,w_2)\ge 3^{\nu+2}$.
If $\mathcal{S}$ contains
a $\nu$-sprinkler with floor a subset of $V(P)$, then the theorem holds, so we assume not. By \ref{getubend} with
$\mu =  3^{\nu+2}/2-\nu$, it follows that $\mathcal{S}$ contains a u-bend as in the claim.
Next we assume that $\mathcal{S}$ is a  w-bend, of size at least $3^{\nu+2}$; then the claim follows from \ref{wtou}.
This proves (1).

\bigskip
Let $\mathcal{S}_1=(L_0\l L_k,U)$ as in (1). Let $U$ have ends $u,s$ where $s$ is the drain. Let $q_0$ be the unique neighbour of
$u$ in $L_{k-1}$; and let $D$ be the path formed by adding the edge $uq_0$ to $U$. 
There is an induced path $q_0\d q_1\c q_n$ such that $\{q_1\l q_n\}= L_k$; and
every vertex in $L_{k-1}$ has a neighbour in $L_k$. Also, $d_G(q_1,q_n)\ge 3^{\nu+2}/2-\nu$, and
so $d_G(q_0,q_n)\ge 3^{\nu+2}/2-\nu-1$.
We may assume that for $0\le i\le k-1$
every vertex in $L_i$ has a neighbour in $L_{i+1}$ (because any other vertex could be removed).
Let $V=L_0\cup\cdots\cup L_k$.

We recall that for $v\in V$, its {\em height} $h(v) =k-i$ where $v\in L_i$; 
and we define the {\em reach} of $v$ to be the maximum $i\ge 1$ such that $q_i$ is a descendant of $v$. (Since
every vertex in $V$ has a descendant in $L_k$, this is well-defined.) Next we show that we may assume that:
\\
\\
(2) {\em For $1\le m\le n$ there do not exist induced paths $R_1,R_2$ of $G[V]$ between $q_0$ and $q_{m}$ with the following properties:
\begin{itemize}
\item $|E(R_1)|+1=|E(R_2)|\le 2\nu+2$; and
\item for all $j$ with $m<j\le n$, $q_j$ has no neighbour in  $V(R_1\cup R_2)\setminus \{q_{m}\}$.
\end{itemize}
}
For suppose that such $m,R_1,R_2$ exist. Since $R_1,R_2$ both have length at most $2\nu+2$ and have ends in $L_k$ and $L_{k-1}$, 
it follows
that every vertex of $R_1\cup R_2$ has height at most $\nu+1$. Indeed, if $y\in V(R_1\cup R_2)$ then there is a subpath of one of
$R_1,R_2$ between $y$ and $q_{m}$, which must have length at least $h(y)$, and since $R_1,R_2$ both have length at most $2\nu+2$, it follows
that $d_G(y,q_0)\le 2\nu+2-h(y)$. Consequently, if
$x\in V$ has a neighbour (say $y$) in $R_1\cup R_2$ then 
$$d_G(x,q_0)\le d_G(y,q_0)+1\le 2\nu-h(y)+3\le 2\nu -h(x)+4.$$
It follows that for every descendant in $L_k$ of such a vertex $x$, its distance from $q_0$ is at most $d_G(x,q_0)+h(x)\le 2\nu+4$.
Since 
 $$d_G(q_0,q_n)\ge 3^{\nu+2}/2-\nu-1> 2\nu+4,$$ 
there exists $m'<n$
such that $d_G(q_0,q_{m'})=2\nu + 4$, and $d_G(q_0,q_j)>2\nu+4$ for all $j$ with $m'<j\le n$. 
Since $q_{m+1}$ has a neighbour in $R_1$, it follows that $d_G(q_{m+1},q_0)\le 2\nu+4$, and so $m'\ge m+1$.
For $0\le i<k$ let $L_i'$ be the set of all vertices in $L_i$ with a descendant in $\{q_j:\;m'<j\le n\}$. It follows that
$$(L_0'\l L_{k-1}',  \{q_j:\;m\le j\le n\}, q_{m})$$ 
is a shower $\mathcal{S}'$ say. It is $\nu$-stable, since 
$L_i'\subseteq L_i$ for $0\le i < k$. 
(It is not contained in $\mathcal{S}$ since the drain is different.)
Let its vertex set be $V'$. If $v\in V'\setminus \{q_{m}\}$, and $v\in L_k$, then $v$ has no neighbour in 
$V(R_1\cup R_2)\setminus \{q_{m}\}$
from the properties of $R_1,R_2$; and if $v\notin L_k$, then 
$v$ has a descendant in $\{q_j:\;m'<j\le n\}$, which therefore
has distance in $G$ more than $2\nu+4$ from $q_0$, and again $v$ has no neighbour in $R_1\cup R_2$.
Thus there are no edges between $V'\setminus \{q_m\}$ and $V(R_1\cup R_2)$ except the edge $q_{m}q_{m+1}$.

Now
$$d_G(q_n,q_{m'+1})\ge d_G(q_n,q_0)-(2\nu+5)\ge 3^{\nu+2}/2-\nu-1 -(2\nu+5) \ge  3^{\nu+1}.$$
If $\nu>2$, then from the inductive hypothesis on $\nu$, applied to $\mathcal{S}'$ and the mat $M'=\{q_{m'+1}\l q_n\}$, we deduce that 
there is a dense subset $\mathcal{A}$ of the $M'$-jetset of $\mathcal{S}'$ of cardinality $\nu-1$, realized by 
a set of $M'$-jets of $\mathcal{S}'$ that are $\nu$-monotone and $3(\nu-1)^2$-wasteful.
If $\nu=2$, let $\mathcal{A}$ be a singleton set containing the length of a $0$-monotone, $0$-wasteful $M'$-jet of $\mathcal{S}'$. In either case,
let $J$ be an $M'$-jet in this set. Its tail has exactly one edge
not in the path $q_{m}\d q_{m+1}\c q_{n}$, and so at most $3(\nu-1)^2+1+(k-1)$ edges of $J$ have an end not in $L_k$.
Moreover,
both $J\cup R_1\cup D$ and $J\cup R_2\cup D$ are jets of $\mathcal{S}_1$, and they are both
$(\nu+1)$-monotone (since every vertex of $R_1\cup R_2$ has height at most $\nu+1$). Since $R_1,R_2$
have length at most $2\nu+2$, it follows that these two jets both have waste 
at most $3(\nu-1)^2 + 1 + 2\nu+2 \le 3\nu^2$. 
Let $|E(R_1)|+|E(D)|=t$; then
$|E(R_2)|+|E(D)| = t+1$, so
for each $a\in \mathcal{A}$, both $a+t, a+t+1$ belong to the jetset of $\mathcal{S}_1$, and
so $\mathcal{A}+\{t,t+1\}$ is a subset of 
the jetset of $\mathcal{S}_1$, and hence of the $M$-jetset of $\mathcal{S}$, and this is a set of at least $\nu$ consecutive integers.
And this set is realized by $M$-jets of $\mathcal{S}$ that are $(\nu+1)$-monotone and have waste at most $3\nu^2$.
Thus in this case the theorem holds. Consequently we may assume that no such $m,R_1,R_2$ exist. This proves (2).

\bigskip
  
For each vertex $v\in V$ with reach $r<n$, let $f(v)\in V$ be defined as follows. There is a monotone path between
$v$ and $q_r$; let $X$ be the set of all vertices $x$ such that $x$ is adjacent to a vertex in a monotone path between $v$ and $q_r$.
Consequently $q_{r+1}\in X$, and so there exists $x\in X$ with reach greater than $r$. Choose such a vertex $x$ with maximum
reach, and define $f(v) = x$. If $v$ has reach $n$ let $f(v)=v$.

Let $v_1=q_0$, and for $1\le i\le \nu-1$ let $v_{i+1}=f(v_{i})$. We need to establish several properties of this sequence.
Let $t\le \nu$ be maximum such that $v_t\ne v_{t-1}$. Thus either $t=\nu$ or $v_t$ has reach $n$. For $1\le i\le t$, 
$r_i$ be the reach of $v_i$; then $r_1 = 1$, and $r_i<r_{i+1}$ for $1\le i<t$.
For $1\le i\le t$ let $P_i$ be a monotone path between $v_i$ and $q_{r_i}$ such that if $i<t$ then
$v_{i+1}$ has a neighbour in $P_i$.
The paths $P_1\l P_{t}$ are pairwise vertex-disjoint, because
the reach of every vertex in $P_i$ is precisely $r_i$, and $r_1\l r_{t}$ are all different.
For $1\le i<t$ let $B_i$ be an induced path of $G[V(P_i)\cup \{v_{i+1}\}]$ between $v_i$ and $v_{i+1}$. 
Thus
for $1\le i\le t$,  $B_1\cup B_2\cup \cdots\cup B_{i-1}\cup P_{i}$ is a path, say $C_i$, between $v_1$ and $q_{r_i}$.
In particular, $B_i$ has length at least one, so there is a unique vertex $y_i$ of $B_i$ adjacent to $v_{i+1}$. 
For $1\le i\le t$, let $\epsilon_i = 1$ if $v_{i+1}, y_i\in L_k$, and $2$ otherwise.
\\
\\
(3) {\em  $t=\nu$; for $1\le i<\nu$, $B_i$ has length $h(v_i)-h(v_{i+1})+\epsilon_i$; and for $1\le i\le \nu$, $C_i$
is an induced path of length $1+\sum_{1\le j< i}\epsilon_j$.}
\\
\\
Let $1\le i<t$. Since $h(y_i)\le h(v_i)$, and $h(v_{i+1})\le h(y_i)+1$, it follows that $h(v_{i+1})\le h(v_i)+1$; and since $h(v_1) = 1$,
it follows inductively that $h(v_i)\le i$ for $1\le i \le t$. Consequently for $1\le i< t$, $y_i$
has height at most $\nu-1$; and since the levelling is $\nu$-stable, it follows that $y_i,v_{i+1}$
do not have the same height unless they both have height zero. 
Moreover, $v_{i+1}$ is not a child of $y_i$, since the reach of $v_{i+1}$ is greater than the reach of $y_i$;
so we have proved that either $v_{i+1}$ is a parent of $y_i$, or $v_{i+1}, y_i$ both have height zero. 
It follows that 
the length of $B_i$ equals $h(v_i)-h(v_{i+1})+\epsilon_i$, for all $i<t$.

For $1\le i\le t$, the path 
$B_1\cup B_2\cup \cdots\cup B_{i-1}$ therefore has length
$$1- h(v_{i})+\sum_{1\le j< i}\epsilon_j,$$
and since $P_{i}$ has length $h(v_{i})$,
it follows that $C_i$ has length $1+\sum_{1\le j< i}\epsilon_j$. Since this quantity is less than $2\nu$, and 
$d_G(u,q_n)\ge 3^{\nu+2}>  2\nu$,
it follows that $r_i<n$. In particular, $r_t<n$, and so $t=\nu$.

We claim that 
for $1\le i\le \nu$, the path $C_i$ is induced;
and prove this by induction on $i$. Certainly $C_1$ is induced, so we may assume inductively that $i<\nu$ and 
$C_{i}$ is induced, and we prove that $C_{i+1}$ is induced. Now $C_{i+1}$ is obtained from a subpath of $C_i$
by adding the edge $y_iv_{i+1}$ and the path $P_{i+1}$; so it suffices to check that  there are no edges between
$B_1\cup B_2\cup \cdots\cup B_{i}$ and $P_{i+1}$ except the edge $y_iv_{i+1}$. Suppose then that $y\in V(B_j)$
for some $j\le i$, and $x\in V(P_{i+1})$, and $xy$ is an edge. Since the reach of $x$ equals $r_{i+1}$,
it follows that $x$ has no neighbour in any of $P_1\l P_{i-1}$, and so $y\in V(P_{i})$. Since also $y\in V(B_j)$
for some $j\le i$, it follows that $y\in V(B_i\cap P_i)$. Since $B_i$ is induced and we may assume that
$(x,y)\ne (v_{i+1}, y_i)$, it follows that $x\ne v_{i+1}$, and so $h(v_{i+1})>0$ and $h(x)<h(v_{i+1})$.
Since $h(v_{i+1})>0$, also $v_{i+1}$ is a parent of $y_i$, and so $h(x)\le h(y_i)$. But $h(y)\ge h(y_i)$, and since
the levelling is $\nu$-stable and $xy$ is an edge, it follows that $y$ is a parent of $x$. But this is impossible
since the reach of $x$ is greater than the reach of $y$. This proves that each $C_i$ is induced, and so completes the proof of (3).

\bigskip
For $1\le j\le n$, let $A_j$ be a monotone path between $q_{j}$ and the shower head $z_0$. Thus $A_j$ has length $k$.
For $1\le i\le \nu$, the reach of every vertex in $A_{r_i+1}$ is at least $r_i+1$, and so is greater than the reach of every vertex
in $C_i$; and so there is a path $J_i$ formed by the union of $D$, $C_i$, the edge $q_{r_i}q_{r_i+1}$, and $A_{r_i+1}$.
\\
\\
(4) {\em For $1\le i\le \nu$ the path $J_i$ is induced.}
\\
\\
Suppose that some $J_t$ is not induced, where $1\le t\le \nu$. 
Consequently some vertex $x$ of $A_{r_t+1}$ is adjacent to some vertex $y$ of $C_t$, and $(x,y)\ne (q_{r_t+1},q_{r_t})$.
Choose such a pair $x,y$ with $x$ of minimum height.
Since $y$ has height at most $\nu$, it follows that $h(x)\ne h(y)$; and $x$ is not a child of $y$ since the reach of $x$
is greater than the reach of $y$. Thus $x$ is a parent of $y$. Let $y\in V(P_j)$ where $j\le t$. Since $x$ has a neighbour in $P_j$,
it follows that the reach of $x$ is at most $r_{j+1}$; and so $r_t< r_{j+1}$. Consequently $t<j+1$, and since $j\le t$ it follows that
$j=t$, and so $y\in V(P_t)$. Let $a$ be the vertex of $A_{r_t+1}$ of height $1$.
Now there are two cases. First suppose that $a$ is nonadjacent to $q_j$ for $r_t+2\le j\le n$. Let 
$R_1$ be the path formed by the union of $C_t$ and the edge $q_{r_t}q_{r_t+1}$, and let $R_2$ be the path formed by the union
of the subpath of $C_t$ between $q_0,y$, the edge $xy$, and the subpath of $A_{r_t+1}$ between $x, q_{r_t+1}$. Note that $R_1$
is induced by (3), and $R_2$ is induced since we chose $xy$ with $x$ of minimum height. Also $R_1$ has length at
most $2\nu$, and $R_2$ has length one more. This is therefore impossible by (2).
Consequently there exists $j>r_t+1$ adjacent to $a$; choose such a value of $j$, maximum.
Let $R_2$ be the path formed by the union of $C_t$ and the path $q_{r_t}\d q_{r_t+1}\d a\d q_j$, and let
$R_1$ be the path formed by the union of the subpath of $C_t$ between $q_0,y$, the edge $xy$, the subpath of $A_{r_t+1}$ between $x,a$,
and the edge $aq_j$. In this case $R_2$ has length at most $2\nu+2$, and $R_1$ has length one less. Since $j\ge r_t+3$ (because
$G$ is triangle-free) it follows that both paths are induced, and again this contradicts (2). Thus there is no such $t$. This proves (4).

\bigskip

Since each $J_i$ is induced, it is therefore a jet for the u-bend $\mathcal{S}_1$ (and hence an $M$-jet for $\mathcal{S}$), of length 
$k+1+\sum_{1\le j< i}\epsilon_j+|V(D)|$, and with tail the path $D$; 
and since $J_i\setminus V(D)$ has length at most $k+2\nu$, and all vertices of $B_i$
have height at most $\nu$, it follows that $J_i$ is $\nu$-monotone and $2\nu$-wasteful (and hence $3\nu^2$-wasteful). The shortest of these jets is 
$J_1$, and it has length $k+1+|V(D)|$. Let $A_0$ be a monotone path
between $q_0$ and $z_0$; then $A_0\cup D$ is also an $M$-jet, of length $k-2+|V(D)|$ (so, three
less than the length of $J_1$). Consequently these $M$-jets realize a subset of the $M$-jetset satisfying the theorem.
This proves \ref{staircase}.~\bbox
\bigskip

(We no longer need u-bends or sprinklers after this point.)
The previous result will have several applications later in the paper. First, let us use it to convert a $\lambda$-stable shower
into a stable shower.

\begin{thm}\label{topstable}
Let $\kappa,\tau\ge 0$ and $\nu\ge 2$ be integers, and let $\rho= 3^{\nu+2}$.
Let $G$ be a triangle-free graph such that $G$ has no hole $\nu$-interval, and 
$\chi^{\rho}(G)\le \kappa$.
If $G$ admits a $\nu$-incomplete $(\nu+2)$-stable shower with floor of chromatic number more than $\kappa+\tau$, then 
$G$ admits a $\nu$-incomplete stable shower with floor of chromatic number more than $\tau$.
\end{thm}
\Proof Let $\mathcal{S}$ be a  $\nu$-incomplete $(\nu+2)$-stable shower $(L_0\l L_k,s)$ in $G$. Thus $k\ge \nu+2$.
Let $j=k-\nu-2$; then $L_i$ is stable for $j\le i< k$.
Let $L_0=\{z_0\}$.
Let $X$ be the set of all vertices $v\in L_j$ such that there is an induced path $P_v$ of $G$ between $v,z_0$ with length
$j+1$, such that every vertex in $P_v$ different from $v$ belongs to one of $L_0\l L_{j-1}$. Let $Y=L_{j}\setminus X$.
Let $X'$ be the set of vertices in $L_k$ with an ancestor in $X$, and $Y'$ the set of vertices in $L_k$ with an ancestor in $Y$.
Thus $X'\cup Y'$ is the floor of $\mathcal{S}$.

Suppose that $\chi(G[X'])>\kappa$. 
For $j\le i\le k-1$, let $L_i'$ be the set of vertices in $L_i$ with an ancestor
in $X$. Then
$$(L_0\l L_{j-1},L_{j}'\l L_{k-1}', L_k,s)$$
is a $\nu$-stable shower
$\mathcal{S}_1$ say, and its floor includes $X'$. This is contained in $\mathcal{S}$, so
 $\mathcal{S}_1$ is $\nu$-incomplete. Since $\chi(G[X'])>\kappa$, there exist $w_1,w_2\in X'$, in the same component of $G[X']$, with 
$d_G(w_1,w_2)>\rho\ge 3^{\nu+2}$. By \ref{staircase}
there is a dense subset $\mathcal{A}$ of the jetset of $\mathcal{S}_1$ of cardinality $\nu$, 
and a set 
$\{J_a:\;a\in \mathcal{A}\}$ 
of $(\nu+1)$-monotone jets for $\mathcal{S}_1$ realizing $\mathcal{A}$.
Thus for each $a\in \mathcal{A}$, $J_a$ contains exactly one vertex of $L_i$ for $0\le i< j$, and exactly one vertex in $L_j'=X$,
say $x$. 
The subpath of $J_a$ between $x,z_0$ has length $j$, and so the subpath
$R_a$ say of $J_a$ between $x,s$ has length $|E(J_a)|-j$. By definition of $X$, the path $P_x$ exists and has length $j+1$;
and since both $R_a$, $P_x$ have exactly one vertex in $L_j$, their union $R_a\cup P_x$ is an induced path between $s,z_0$
of length exactly one more than the length of $J_a$. Now both $J_a$ and $R_a\cup P_x$ are jets of
$\mathcal{S}_1$ and hence of $\mathcal{S}$. Thus
$\mathcal{A}+\{0,1\}$
is a subset of the jetset of $\mathcal{S}$. But this set consists of at least $\nu+1$ consecutive integers, since $\mathcal{A}$ is dense
of cardinality $\nu$; and this is impossible since $\mathcal{S}$ is not $\nu$-complete.
This proves that $\chi(G[X'])\le \kappa$.

Consequently $\chi(G[Y'])> \tau$. 
For $0\le i\le j$, let $L_i'$ be the set of vertices in 
$L_i$ with a descendant in $Y$, and for $j+1\le i\le k$, let $L_i'$ be the set of vertices in $L_i$ with an ancestor in $Y$.
Then $(L_0'\l L_{k-1}', L_k,s)$ is a shower $\mathcal{S}'$ say, with floor of chromatic number
more than $\tau$ since its floor includes $Y'$. This is contained in $\mathcal{S}$,
so $\mathcal{S}'$ is $\nu$-incomplete. We claim that $\mathcal{S}'$ is stable. For certainly
$L_{j}'\l L_{k-1}'$ are stable, since $\mathcal{S}$ is $(\nu+2)$-stable. Suppose that $0\le i\le j-1$ and $y,y'\in L_{i}'$
are adjacent. Since $y$ has a descendant in $Y$, there is a path between
$y$ and $Y$ of length $j-i$; and since $y'\in L_{i}'$, there is a path between $y', z_0$ of length $i$. Since $G$ is triangle-free,
$yy'$ is the only edge between these two paths; and so their union, together with this edge, is an induced path between $y,z_0$
of length $j+1$, contradicting that $y\notin X$.  This proves that $\mathcal{S}'$ is stable; and so the theorem holds.
This proves \ref{topstable}.~\bbox

\bigskip

We deduce:

\begin{thm}\label{getstable}
Let $\tau \ge 0$ and $\nu\ge 2$ be integers, and let $\rho= 3^{\nu+2}$.
Let $G$ be a triangle-free graph, such that $G$ has no hole $\nu$-interval, and
$\chi^{\rho}(G)\le \kappa$.
If $\chi(G)>44 \nu (\kappa+\tau)^{(\nu+1)^2} + 4\kappa$, then
$G$ admits a $\nu$-incomplete stable shower with floor of chromatic number more than $\tau$.
\end{thm}
\Proof
By \ref{firstshower}, there is a $\nu$-incomplete shower $(L_0\l L_k,s)$ in $G$, with floor of chromatic number more than 
$\nu(\kappa+\tau)^{(\nu+1)^2}$.
Then $k> \rho$, since $\nu(\kappa+\tau)^{(\nu+1)^2}\ge \kappa$. Since $\rho\ge (\nu+3)(\nu-1)+1$,
\ref{basestable} (with $\lambda = \nu+2$) implies that there is a $(\nu+2)$-stable $\nu$-incomplete shower in $G$ with 
floor of chromatic number more than $\kappa+\tau$, 
so the result follows from \ref{topstable}. This proves \ref{getstable}.~\bbox

\bigskip

The reason for controlling the waste of the jets that are output by \ref{staircase} is that a jet with bounded waste can be covered by a bounded number of monotone paths. More precisely:

\begin{thm}\label{wastecover}
Let $\mathcal{S}=(L_0\l L_k, s)$ be a shower in a graph $G$, and let $J$ be a $\mu$-wasteful jet of $\mathcal{S}$. 
Then there is a set of
at most $\mu+1$ monotone paths of $\mathcal{S}$
such that every vertex of $J$ in $L_0\cup\cdots\cup L_{k-1}$ belongs to one of these paths.
\end{thm}
\Proof Choose $d\in V(J)$ such the tail $T$ of $J$ has ends $d,s$. Then no vertex of $T$ belongs to $L_0\cup\cdots\cup L_{k-1}$
except $d$. Let $P$ be the subpath of $J$ between $z_0,d$, where $z_0\in L_0$. 
At most $k-1+\mu$ edges of $P$ have an end not in $L_k$, since the waste of $J$ is at most $\mu$. Let us say the {\em height} of an edge
$uv$ of $P$ is the maximum of the heights of $u,v$. Thus at most $k-1+\mu$ edges of $P$ have nonzero height. As $P$ is traversed
starting from $d$, the number of edges in it that have height at least 2 and different from the heights of 
all previous edges is at least
$k-1$, since the difference of the heights of $z_0,d$ is $k-1$; and so there are at most $\mu$ edges of $P$
that have height 1 or the same nonzero height as some earlier edge. 
By removing all such edges, we decompose $P$ into at most $\mu+1$ paths each of which is either monotone or a path of $G[L_k]$; and
every vertex of $P$ in $L_0\cup\cdots\cup L_{k-1}$ belongs to 
one of these monotone paths.
This proves \ref{wastecover}.~\bbox

\section{Stable showers}

From now on, there is no need to consider general showers; 
we might as well just concern ourselves with stable showers, in view of \ref{getstable}.
To complete the proof of \ref{localbound}, we only need to show that if $\nu,\kappa, G$ satisfy the hypotheses
of \ref{localbound} then every $\nu$-incomplete stable shower in $G$ has floor with bounded $\chi$, and that is the goal 
of the remainder of the paper.

We are concerned with a triangle-free graph which admits no hole $\nu$-interval; and
we will not need to use induction on $\nu$ any more; so from now on we shall fix $\nu\ge 2$, to avoid having to carry it along.
We might as well also set $\rho = 3^{\nu+2}+4$, for the remainder of the paper, and fix $\kappa\ge 0$. Let us say a graph $G$
is a {\em candidate} if $G$ is triangle-free, and admits no hole $\nu$-interval, and $\chi^{\rho}(G)\le \kappa$.
Our eventual goal is to prove that every stable shower in every candidate has floor of bounded $\chi$.

Let  $\mathcal{S}$ be a stable shower, with vertex set $V$, and let $M$ be a mat. For $X\subseteq V$, we denote
the set of vertices in $M$ with an ancestor in $X$ by $M(X)$
(and we write $M(v)$ for $M(\{v\})$).

We already defined ``containment'' for showers, but now we need a different inclusion relation.
Let $\mathcal{S}=(L_0\l L_k, s)$ be a stable shower, and let 
$\mathcal{S}'=(L_0'\l L_{k'}')$ be a levelling,
both in a graph $G$. We say that $\mathcal{S}'$ is a 
{\em sublevelling}
of $\mathcal{S}$ if
$k'\le k$, and $L_{i}'\subseteq L_{i+k-k'}$ for $0\le i\le k'$.

If $\mathcal{S}=(L_0\l L_k, s)$ is a shower, we define $U(\mathcal{S})$ to be $L_0\cup \cdots\cup L_{k-1}$. 
(Note that this is different from $V(\mathcal{S})$, as we do not include $L_k$.)

\begin{thm}\label{splitshower}
Let $G$ be a candidate.
Let $\mathcal{S}=(L_0\l L_k,s)$ be a stable shower in $G$, and let $z_1,z_2\in U(\mathcal{S})$, either equal or nonadjacent.
For $i = 1,2$, let $\mathcal{S}_i$ be a 
sublevelling of $\mathcal{S}$ with vertex set $V_i$ and head $z_i$ respectively, disjoint except possibly $z_1=z_2$ (more precisely,
$V_1\cap V_2=\{z_1\}\cap \{z_2\}$), and let $M_i$ be 
the base of $\mathcal{S}_i$. Let $\chi(M_1)> \kappa$.
Then either
\begin{itemize}
\item there are $\nu$ induced paths
$Q_0\l Q_{\nu-1}$ of $G[V_1\cup V_2\cup L_k]$ between $z_1,z_2$, such that $|E(Q_i)|=|E(Q_0)|+i$ for $0\le i<\nu$ (and in particular
$z_1,z_2$ are distinct and nonadjacent); or
\item $\chi(M_2\setminus M_2(X))\le 2\kappa$, where $X$ denotes the set of vertices in $V_2\setminus \{z_1\}$ that have
a neighbour 
in $V_1\setminus \{z_2\}$; and 
if $\chi(M_2)> 2\kappa$ 
then there is a
monotone path $R$ of $G[V_1]$ between $z_1$ and $M_1$ such that 
$(\nu+1)(3\nu^2+1)\chi(M_2(X(R))) \ge \chi(M_2) - 2\kappa$,
where 
$X(R)$ denotes the set of vertices in $V_2\setminus \{z_1\}$
that have a neighbour in $V(R)\setminus \{z_2\}$.
\end{itemize}
\end{thm}
\Proof
Choose a component of $G[M_1]$ with maximum chromatic number;
and since this chromatic number is larger than $\kappa$, it follows that there are two vertices of this component with distance
more than $\rho$ (in $G$). Consequently there is a path $P_1$ with $V(P_1)\subseteq M_1$ joining two vertices
with distance at least $3^{\nu+2}$ (in $G$).
Choose a minimal such path $P_1$, and let $w_1$ be one of its ends. From the minimality of
$P_1$ it follows that $d_G(w_1,v)\le 3^{\nu+2}$
for every vertex $v$ of $P_1$. (That concludes the role of $M_1$ in this proof.)

We may assume that $\chi(M_2)> 2\kappa$, since otherwise the second bullet of the theorem holds.
Let $C_2$ be a connected induced subgraph of $G[M_2\setminus N^{\rho}[w_1]]$ with $\chi(C_2)>\kappa$.
In addition, choose $C_2$ with 
$V(C_2)\cap M_2(X)=\emptyset$ if possible, where $X$ denotes the set of vertices in 
$V_2\setminus \{z_1\}$ that have a neighbour       
in $V_1\setminus \{z_2\}$.
Every path of $G$ between $C_2$ and $P_1$ has length at least $3$, since $\rho\ge 3^{\nu+2}+3$.
(Now we are finished with $M_2$.)

Let $z_1\in L_{h_1}$, and for $h_1\le i\le k$ let $L_i^1$ be the set of all $\mathcal{S}_1$-descendants of $z_1$ in $L_i$ 
with an $\mathcal{S}_1$-descendant in $P_1$. Thus $L_k^1=V(P_1)$.
Let $V_1' = L_{h_1}^1\cup\cdots\cup L_{k}^1$.
Since $G[L_k]$ is connected, there is a path
of $G[L_k]$ between $V(P_1)$ and $C_2$; let $D$ be a minimal path of $G[L_k]$ such that
one end (say $d_1$) has a neighbour in $V(P_1)\cup L^1_{k-1}$ and the other (say $d_2$) has a neighbour in $C_2$.
\\
\\
(1) {\em There is a set $\mathcal{A}_1$ of integers, of cardinality at most $\nu+1$,
including a dense subset of cardinality $\nu$, and
containing two integers $x,y$ with $y-x\in \{1,3\}$, such that the following holds. For each $a\in \mathcal{A}_1$ there is an induced path $J_a$ of $G$
between $d_1,z_1$ of length $a$, such that
\begin{itemize}
\item $V(J_a)\subseteq V_1'\cup\{d_1\}$; and 
\item there is a set of $3\nu^2+1$ monotone paths of $G[V_1']$ between $V(P_1)$ and $z_1$,
such that every vertex of $V(J_a)\setminus (V(P_1)\cup \{d_1\})$ belongs to one of these paths.
\end{itemize}
}
\noindent Let $D_1$ be the one-vertex path with vertex $d_1$. If $d_1$ has no neighbour in $V(P_1)$, then
$$(L_{h_1}^1\l L_{k}^1, D_1)$$ 
is a w-bend $\mathcal{S}_1'$ of size at least $3^{\nu+2}$; and otherwise
$(L_{h_1}^1\l L_{k}^1\cup \{d_1\}, d_1)$ is a shower  $\mathcal{S}_1'$. In either case we can apply \ref{staircase}
to  $\mathcal{S}_1'$, and deduce that
there is a subset $\mathcal{A}_1$ of the jetset of $\mathcal{S}_1'$, of cardinality at most $\nu+1$, including a dense subset of
cardinality $\nu$, and containing two integers $x,y$ with $y-x\in \{1,3\}$; and realized by
a set of jets of $\mathcal{S}_1'$ that are $3\nu^2$-wasteful.
By \ref{wastecover}, this proves (1).

\bigskip
Since $|\mathcal{A}_1|\le \nu+1$, there is a set of at most $(\nu+1)(3\nu^2+1)$ monotone paths of $G[V_1']$ between $V(P_1)$ and $z_1$
such that, if $Y$ denotes the set of vertices
in these paths, then $V(J_a)\subseteq Y\cup V(P_1)\cup \{d_1\}$ for each $a\in \mathcal{A}_1$. Let $X'$ 
denote the set of vertices in $V_2\setminus \{z_1\}$
that have a neighbour in $Y\setminus \{z_2\}$.
Let $h_2$ be such that $z_2\in L_{h_2}$, and for $h_2\le i\le k$ let $L_i^2$ be the set of vertices $v\in L_i$
such that there is a monotone path of $G[V_2\setminus X']$ between $z_2, V(C_2)$ containing $v$.
It follows that no vertex in $Y\setminus \{z_2\}$ has a neighbour in $(L_{h_2}^2\cup \cdots\cup  L_{k}^2)\setminus \{z_1\}$,
and $V(C_2)\setminus M_2(X')\subseteq L^2_k$.
\\
\\
(2) {\em If $\chi(L_k^2)> \kappa$ then the first bullet of the theorem holds.}
\\
\\
For then there exists
an induced path $P_2$ of $G[L_k^2]$ with ends at distance at least $\rho$.
Since $d_2$ has a neighbour in $C_2$, it follows that $G[V(C_2)\cup V(D)]$ is connected.
Thus
$$(\{z_2\}, L_{h_2+1}^2\l L_{k-1}^2, V(C_2)\cup V(D), d_1)$$
is  a shower $\mathcal{S}_2'$, and $L_k^2$ is a mat; and by \ref{staircase},
there is a dense subset $\mathcal{A}_2$ of the $L^2_k$-jetset of $\mathcal{S}_2'$ of cardinality $\nu$. 
We claim that $\mathcal{A}_1+\mathcal{A}_2$ contains a set $\mathcal{B}$ of $\nu$ consecutive integers. 
To see this, suppose first that there are two
consecutive integers $a,a+1\in \mathcal{A}_2$. Let $\mathcal{A}'$ be a dense subset of $\mathcal{A}_1$ of cardinality $\nu$;
then $\mathcal{A}'+\{a,a+1\}$ consists of at least $\nu+1$ consecutive integers, all contained in $\mathcal{A}_1+\mathcal{A}_2$
as required. We may assume that that no two members of $\mathcal{A}_2$ are consecutive. Since $\mathcal{A}_2$ is dense
of cardinality $\nu$, there exists $s$ such that $s,s+2,s+4\l s+2(\nu-1)\in \mathcal{A}_2$. But there exist $x,y\in \mathcal{A}_1$
with $y-x\in \{1,3\}$; and then
$$\{s,s+2,s+4\l s+2(\nu-1)\}+\{x,y\}$$
contains $\nu$ consecutive integers (indeed, almost $2\nu$). This proves that $\mathcal{B}$ exists.

If $z_1=z_2$
then for every $J_a\:(a\in \mathcal{A}_1)$
and every $L^2_k$-jet of $\mathcal{S}_2'$, their union is a hole; and so $G$ has holes of every length in $\mathcal{B}$,
and so has a hole $\nu$-interval, which is impossible since $G$ is a candidate. Thus $z_1\ne z_2$, and so they are nonadjacent;
but then for every $J_a\:(a\in \mathcal{A}_1)$
and every $L^2_k$-jet of $\mathcal{S}_2'$, their union is an induced path between $z_1,z_2$, and so the first bullet of the theorem holds. 
This proves (2).

\bigskip
We may therefore assume that 
$\chi(L_k^2)\le \kappa$, and we will show that the second bullet of the theorem holds. 
Consequently $V(C_2)\not\subseteq L_k^2$, and therefore
$V(C_2)\cap M_2(X')\ne \emptyset$.
From the choice of $C_2$, it follows that $\chi(M_2\setminus M_2(X))\le 2\kappa$
(for otherwise we could have chosen $C_2$ with $V(C_2)\subseteq M_2\setminus (M_2(X)\cup N^{\rho}[w_1])$). This proves
the first statement of the second bullet.

We can choose $C_2$ such that $\chi(C_2)\ge \chi(M_2)-\kappa$, and since 
$V(C_2)\setminus M_2(X')\subseteq L^2_k$ and $\chi(L_k^2)\le \kappa$, it follows that 
$\chi(V(C_2)\cap M_2(X'))\ge \chi(M_2)-2\kappa$, and consequently $\chi(V(C_2)\cap M_2(X))\ge \chi(M_2)-2\kappa$.
Thus, one of the $(\nu+1)(3\nu^2+1)$ monotone paths satisfies the second statement of the second bullet.
This proves \ref{splitshower}.~\bbox

(We will not need w-bends after this point.)
There is a special case of \ref{splitshower} that we will use several times, and we extract it to make application easier.
\begin{thm}\label{splitshower2}
Let $G$ be a candidate.
Let $\mathcal{S}=(L_0\l L_k,s)$ be a stable shower in $G$, and let $z_1\in U(\mathcal{S})$.
Let $A,B$ be disjoint sets of children of $z_1$.
Let $M$ be a mat for $\mathcal{S}$, and suppose that $\chi(M(B)\setminus M(A))>\kappa$.
Then $\chi(M(A)\setminus M(B))\le 2\kappa$, and 
if $\chi(M(A))> 2\kappa$,
then there is a
monotone path $R$ between $B, M(B)\setminus M(A)$ such that 
$(\nu+1)(3\nu^2+1)\chi(M(X))\ge \chi(M(A)) - 2\kappa$,
where $X$ denotes the set of vertices with 
a parent in $V(R)$
and an $\mathcal{S}$-ancestor in $A$ and an $\mathcal{S}$-descendant in $M(A)$.
\end{thm}
\Proof
Let $\mathcal{S}_1$ be the maximal sublevelling 
of $\mathcal{S}$ with head $z_1$ such that for every vertex $v\ne z_1$ of its vertex set, $v$ has an $\mathcal{S}$-ancestor in 
$B$ and has an $\mathcal{S}$-descendant in $M(B)\setminus M(A)$ (and hence has no $\mathcal{S}$-ancestor in $A$);
and let $\mathcal{S}_2$ be the maximal 
sublevelling of $\mathcal{S}$ with head $z_1$ such that every vertex of its vertex set except $z_1$ has an $\mathcal{S}$-ancestor 
in $A$ and an $\mathcal{S}$-descendant in $M(A)$.
Let their vertex sets be $V_1,V_2$ respectively. Then $V_1\cap V_2 = \{z_1\}$, and no vertex in $V_1$ has a parent in 
$V_2\setminus \{z_1\}$. 
Also $M(B)\setminus M(A)$ is 
the base of $\mathcal{S}_1$,
and $M(A)$ is the base of $\mathcal{S}_2$. 
By \ref{splitshower}, the result follows. This proves \ref{splitshower2}.~\bbox

\begin{thm}\label{bigson}
Let $G$ be a candidate.
Let $(L_0\l L_k,s)$ be a stable shower in $G$, let $z_1\in L_0\cup\cdots\cup L_{k-1}$, 
let $Y$ be a subset of the set of children of $z_1$,
and let $M\subseteq L_k$ such that every vertex in $M$ has an ancestor in $Y$. 
If $$\chi(M)>  ((\nu+1)(3\nu^2+1)+7)\kappa,$$
then there exists $z_2\in Y$ such that $\chi (M(z_2))\ge \chi(M) - ((\nu+1)(3\nu^2+1)+7)\kappa$.
\end{thm}
\Proof 
Let $\tau=\chi(M)-((\nu+1)(3\nu^2+1)+7)\kappa$, and 
choose $A\subseteq Y$ minimal such that
$\chi(M(A))\ge 2\kappa + \tau$. Suppose that $\chi(M(A))\ge 3\kappa + \tau$, and choose $z_2\in A$;
then from the minimality of $A$, 
$\chi(M(A\setminus \{z_2\}))< 2\kappa + \tau$, and so
$$\chi(M(z_2)\setminus M(A\setminus \{z_2\}))> \kappa.$$ 
By \ref{splitshower2} applied to $A\setminus \{z_2\}$ and $\{z_2\}$,
it follows that $\chi(M(A)\setminus M(z_2))\le 2\kappa$; and since
$\chi(M(A))\ge 2\kappa + \tau$, it follows that $\chi(M(z_2))\ge \tau$, as required.

We may assume therefore that $\chi(M(A))<  3\kappa + \tau$.
Let  $\mu = (\nu+1)(3\nu^2+1)\kappa$; then
$$\chi(M) =  \tau + \mu + 7\kappa\ge \chi(M(A)) + \mu + 4\kappa,$$ 
so we may
choose $B\subseteq Y$ with $A\subseteq B$, minimal such that 
$\chi(M(B))> \chi(M(A))+ \mu+2\kappa.$
Again, by the same argument, we may assume that
$$\chi(M(B))\le  \chi(M(A))+ \mu+3\kappa< \tau+ \mu+6\kappa;$$
and since $\chi(M)=  \tau + \mu+7\kappa,$
it follows that 
$\chi(M\setminus M(B))> \kappa.$
Since $\chi(M(B))> 2\kappa$, \ref{splitshower2} applied to the mat $M\setminus M(A)$ and the sets $B, Y\setminus B$,
implies that
there is a
monotone path $R$ between $Y\setminus B, M\setminus M(B)$ such that
$$(\nu+1)(3\nu^2+1)\chi(M(X)\setminus M(A))\ge \chi(M(B)\setminus M(A)) - 2\kappa > \mu= (\nu+1)(3\nu^2+1)\kappa,$$
and so $\chi(M(X)\setminus M(A))>\kappa$,
where $X$ denotes the set of vertices with
a parent in $V(R)$
and an $\mathcal{S}$-ancestor in $B$ and an $\mathcal{S}$-descendant in $M(B)\setminus M(A)$.
Let $z_2$ be the vertex of this path in $Y\setminus B$; then 
$$\chi(M(z_2)\setminus M(A))>\kappa.$$ 
By \ref{splitshower2} again, applied
to the mat $M$ and the sets $A, \{z_2\}$, it follows that $\chi(M(A)\setminus M(z_2))\le 2\kappa$. Since
 $\chi(M(A))\ge 2\kappa+\tau$, it follows that
$\chi(M(z_2))\ge  \tau$. This proves \ref{bigson}.~\bbox

\section{What is going on?}

It might be helpful at this stage to make some general remarks about where the proof is going. Look at some vertex $z_i\in L_i$, 
such that the set of descendants of $z_i$ has large chromatic number. By \ref{bigson} there is a child $z_{i+1}$ of $z_i$
whose descendants have chromatic number almost as large (reduced by an additive constant); and \ref{nosplit} tells us that the
set of vertices in the base that are descendants of $z_i$ and not of $z_{i+1}$ has bounded chromatic number. This suggests that 
we start with $z_0\in L_0$, and generate a sequence $z_1\l z_t$ as above, until it stops. This sequence induces a useful
partition of $V(G)$; for each $i$ we look at the descendants of $z_i$ that are not descendants of $z_{i+1}$. If $v$ is a descendant of
$z_i$ and not $z_{i+1}$, say the ``reach" of $v$ is $i$. So the set of vertices in $M$ with any given reach has bounded chromatic
number, and we would like to exploit the partition given by the reach numbers. (This is the start of the proof of \ref{noknobbly2}.)

For each vertex $v\in M$ with reach $i$, there is an induced path from $v$ to $z_i$ such that all its vertices have reach $i$; we call 
such paths ``vertical'' (they are monotone, but not all monotone paths are vertical). Follow this vertical path from $v$ until
it first contains a neighbour $x$ of some $z_j$; then $x$ might have just a parent in $\{z_0\l z_t\}$, or just a child in this set,
or both. Thus $M$ is divided into three parts, and we will bound their chromatic numbers separately. For vertices $v$ such that
the corresponding $x$ has both a parent and child, it follows that $x$ is adjacent to $z_i,z_{i+2}$, and so $x$ is in some sense
a copy of $z_{i+1}$; and to handle such vertices $x$, we replace the sequence $z_0,z_1\l z_t$ by a sequence of sets of vertices,
each complete to the next.
This sequence of sets is called a ``wand'', but for this sketch let us confine ourselves to wands
where all the sets are singletons. (We only use the more general wands once, in the proof of \ref{noknobbly2}.) It remains to handle the $x$'s with only a parent in $\{z_0\l z_t\}$, and those with only
a child. 

The ones with only a child are suggestive. Suppose that the corresponding set of vertices in $M$ has large chromatic number.
Delete everything except $z_0\l z_t$ and the vertical paths that lead to vertices $x$ 
of this ``only-a-child'' type. Then we get a new shower, still with big $\chi$, and all distances from $z_0$ to vertices not in the wand 
are two more than before. (This is called
``raising the wand''). We would like to say that
if in this smaller shower we can guarantee some jetset $\mathcal{A}$ (up to shifting), then in the original shower we
can guarantee $\mathcal{A}+\{0,2\}$. Unfortunately this does not seem to be true; but if in the smaller shower there is a wand that 
can be raised to get a third shower still with big $\chi$, this third shower has the property we want.
Since in the third shower we can at least get two jets whose lengths differ by 1 or 3 by \ref{staircase}, we can now get 
two jets that differ by one in the first shower. 
If in this third shower we can again find a wand giving us the same situation, we could get three jets of consecutive 
lengths in the original shower, and this cannot go on arbitrarily, or we would get many jets of consecutive lengths in
the big shower and win. More precisely,
let $\sigma\le \nu$ be maximum such that every stable shower with large enough $\chi$ has $\sigma$ jets of consecutive lengths. 
Our goal is 
to prove that $\sigma= \nu$, so we assume not, and assume we have a shower with large $\chi$ in which there are no $\sigma+1$
jets of consecutive lengths. Then it follows that raising any wand gives a shower in which raising another wand gives a shower
with bounded $\chi$. So we might as well assume that that we have a shower with large $\chi$ in which raising any wand gives
bounded $\chi$. The details are in \ref{knoblemma3}.
This is how we manage the ``only-a-child'' type $x$'s.

To handle the ``only-a-parent'' $x$'s is more complicated. The idea is that we partition the set of possible reach values into
a few intervals, such that the vertices in $M$ with reach in each interval have large chromatic number. The vertices
in each interval can all be accessed from the corresponding $z_i$ by a vertical path, and the vertical paths for different
intervals are disjoint, and we know a great deal about the edges between them. (In particular, since we are in the ``only-a-parent''
case, nothing bad happens very close to the wand.) That allows us to apply \ref{splitshower} to
obtain a contradiction. For instance, suppose we divide into two intervals, splitting $M$ into two large $\chi$ subsets.
We apply \ref{splitshower}. The second outcome of \ref{splitshower}, involving a monotone path $R$, is impossible, because
the vertices of $R$ would have larger reach than the vertices in the other sublevelling (the set $V_2$ of \ref{splitshower}) and so
all vertices in $V_2$ with a neighbour in $R$ would have a child and not a parent in $R$, and then we could treat $R$ as a 
wand and raise it to get a contradiction. Thus the first outcome of \ref{splitshower} must always hold, and we are equipped
with a set of paths joining the two shower heads with many different but similar lengths. We can do this simultaneously
with different pairs of sets if we partition $M$ into several parts instead of just two; and we 
can chain two of these
objects together, to get many paths of consecutive lengths, in such a way that these paths can be completed to holes of many
consecutive lengths. This is the argument of section 11.

\section{Shower completeness}

To go further we use a global induction that we explain next. For $n\ge 2$, a set of integers is {\em $n$-solid} if some subset consists
of $n$ consecutive integers. It is {\em $1$-solid} if it contains
two integers that differ by 1 or 3. 
A key observation is that if a set $\mathcal{A}$ of integers is $n$-solid where $n>0$, then $\mathcal{A}+\{0,2\}$
is $(n+1)$-solid.
Let us say a shower is {\em $n$-complete} over a mat $M$ if its $M$-jetset is $n$-solid. (For $n\ge 2$ this agrees with our earlier definition.)
Now  \ref{staircase} implies that in every candidate, all stable showers with 
a mat $M$ of large enough chromatic number are $1$-complete over $M$; and as we have seen, to finish the proof of our main theorem \ref{localbound}
we only need to show that 
all stable showers with a mat $M$ of large enough chromatic number are $\nu$-complete over $M$. The induction just mentioned is that we 
assume that for some $\sigma>0$, all stable showers with a mat $M$ of large enough chromatic number are $\sigma$-complete over $M$; and
we will prove the same with $\sigma$ replaced by $\sigma+1$.

For $\sigma>0$, let us say an integer $\zeta\ge 0$ is a {\em sidekick} for $\sigma$ if 
for every candidate $G$, and every stable shower $\mathcal{S}$ in $G$, $\mathcal{S}$ is $\sigma$-complete over $M$ 
for every mat $M$ for $\mathcal{S}$
with chromatic number more than $\zeta$.

Next we need another inclusion relation for showers, as follows. Let $\mathcal{S}=(L_0\l L_k, s)$ be a stable shower, 
with vertex set $V$,
and let $\mathcal{S}'=(L_0'\l L_{k'}', s')$ be a 
shower, both in a graph $G$. Let $P$ be an induced path of $G[V]$ between $L_0, L_0'$. Suppose that
\begin{itemize}
\item $s = s'$;
\item $L_0'\l L_{k'-1}'\subseteq L_0\cup \cdots\cup L_{k-1}$;
\item $L_{k'}'\subseteq L_k$; and
\item no vertex of $P$ belongs to $L_1'\cup\cdots\cup  L_{k'}'$, and no vertex of $P$ has a neighbour in this set
except the vertex in $L_0'$.
\end{itemize}
In this situation we say that $\mathcal{S'}$ is {\em included} in $\mathcal{S}$, and $P$ is a {\em pipe}. 
Note that there may be vertices $u,v$ such that $v$ is a child of $u$ in $\mathcal{S}$, and $u$ is a child of $v$
in $\mathcal{S}'$.
Nevertheless, it follows that $\mathcal{S'}$
is a stable shower, because the subgraph induced on $L_0\cup\cdots\cup L_{k-1}$ is bipartite. 

Let $\mathcal{S}'$ be included in $\mathcal{S}$, with a pipe $P$. For every jet
$J$ of $\mathcal{S}'$, $J\cup P$ is a jet of $\mathcal{S}$; and consequently, if the jetsets of the two showers are
$\mathcal{A}, \mathcal{A}'$ respectively then $\mathcal{A}'+\{|E(P)|\}\subseteq \mathcal{A}$. Thus if $\mathcal{S}'$ is $n$-complete
for some $n$, then so is $\mathcal{S}$. If $M, M'$ are mats for $\mathcal{S}, \mathcal{S}'$ respectively, and $M'\subseteq M$,
then for every $M'$-jet $J$ of $\mathcal{S}'$, $J\cup P$ is an $M$-jet of $\mathcal{S}$; and so the same relation
holds between the $M$- and $M'$-jetsets of the two showers. Note that the floor of $\mathcal{S}'$ is a subset of the floor of 
$\mathcal{S}$, but for an individual vertex $v$, there may be $\mathcal{S}'$-descendants of $v$ that are not $\mathcal{S}$-descendants.
(This is not the case for sublevellings.)

Let $\mathcal{S}'$ be included in $\mathcal{S}$. We say a {\em switch} for $\mathcal{S}'$ in $\mathcal{S}$
is a pair $(P_1,P_2)$ of pipes such that $|E(P_2)|=|E(P_1)|+2$.

\begin{thm}\label{noswitch}
Let $\zeta$ be a sidekick for $\sigma$.
Let $\mathcal{S}$ be a stable shower in a candidate $G$, and let $\mathcal{S}$ include a shower $\mathcal{S}'$.
Let $M, M'$ be mats for $\mathcal{S}$, $\mathcal{S}'$ respectively, with $M'\subseteq M$.
If $\mathcal{S}$ is not $(\sigma+1)$-complete over $M$, and $\chi(M')>\zeta$, then 
there is no switch for $\mathcal{S}'$ in $\mathcal{S}$.
\end{thm}
\Proof Let $\mathcal{S}$, $\mathcal{S}'$ have heads $z_0, z_1$ respectively, and suppose that $(P_1,P_2)$ is a switch 
for $\mathcal{S}'$ in $\mathcal{S}$. 
Let $\mathcal{A}$ be the $M$-jetset of $\mathcal{S}$, and let $\mathcal{A}'$ be the $M'$-jetset of $\mathcal{S}'$.
As we saw above,
$$\mathcal{A}'+\{|E(P_1)|,|E(P_1)|+2\}\subseteq \mathcal{A}.$$
Since $\chi(M')>\zeta$ and $\zeta$ is a sidekick for $\sigma$,
it follows that $\mathcal{S}'$ is $\sigma$-complete over $M'$.
Consequently $\mathcal{A}'+\{|E(P_1)|,|E(P_1)+2\}$ is $(\sigma+1)$-complete, and hence so is $\mathcal{A}$, a contradiction.
This proves \ref{noswitch}.~\bbox


\section{The shadow of a wand}

Let $\mathcal{S}=(L_0\l L_k,s)$ be a stable shower.
A {\em wand} $\mathcal{W}$ in $\mathcal{S}$
is a sequence $(W_0\l W_{t})$ with the following properties:
\begin{itemize}
\item $0\le t\le k-2$;
\item $\emptyset\ne W_i\subseteq L_i$ for $0\le i\le t$; and
\item every vertex in $W_i$ is adjacent to every vertex in $W_{i+1}$ for $0\le i\le t-1$.
\end{itemize}
We define
$V(\mathcal{W}) = W_0\cup \cdots\cup W_{t}$.

Let $(W_0\l W_{t})$ be a wand $\mathcal{W}$ in $\mathcal{S}$. If $u\in W_i$ for some $i$, we say that a neighbour $v$ of $u$
is an {\em up-neighbour} of $u$ (relative to $\mathcal{W}$) if
\begin{itemize}
\item $v\notin V(\mathcal{W})$;
\item $v\in L_{i-1}$ (and therefore $i\ge 2$); and
\item every neighbour of $v$ in $V(\mathcal{W})$ belongs to $W_i$ (and therefore $i\ge 3$).
\end{itemize}

For $0\le i\le t-1$, let $T_i$ be the set of all vertices $v\in L_i$
such that $v$ is an up-neighbour of some vertex in $W_{i+1}$.
Let $T = T_0\cup \cdots\cup T_{t-1}$.
For $v\in T$, a {\em post} with {\em top $v$} (in $\mathcal{S}$ for $\mathcal{W}$) is a monotone path between $v$ and $L_k$ such that
no vertex of this path has a parent in $V(\mathcal{W})$ (and consequently no vertex of this path
belongs to $V(\mathcal{W})$).
A post with top $v$ therefore provides an induced path between each neighbour ($u$ say) of $v$ in $V(\mathcal{W})$
and $L_k$, of length two more than a monotone path between $u$ and $L_k$, and both paths can be extended to induced paths between $L_0$
and $L_k$ by adding a path with vertex set within $V(\mathcal{W})$. We shall exploit this later.
For $0\le i\le k$, let $S_i$ be the set of all vertices $v\in L_i$
that belong to a post with top in $T$.
(Thus $S_i\subseteq L_i\setminus V(\mathcal{W})$, and $S_0 = \emptyset$.) If $M$ is a mat for $\mathcal{S}$,
we call $M\cap S_k$ the {\em shadow} (in $\mathcal{S}$, over $M$) of the wand.

Showers in which no wand shadow has large $\chi$ are easier to work with than general showers. 
In this section we prove that their mats have bounded chromatic number.
The proof requires several steps. We begin with:

\begin{thm}\label{nosplit}
Let $\mathcal{S}$ be a stable shower with mat $M$ in a candidate $G$, such that
every wand in $\mathcal{S}$ has shadow over $M$ with chromatic number
at most $\tau$. Let $z\in U(\mathcal{S})$, and let $A,B$ be disjoint sets of children of $z$.
If $\chi(M(A))>  \kappa$ then $\chi(M(B)\setminus M(A))\le (\nu+1)(3\nu^2+1)(\tau+\kappa) +2\kappa$.
\end{thm}
\Proof
Suppose not. Let $\mathcal{S}_1$ be a sublevelling
of $\mathcal{S}$ with head $z$ and base $M(A)$
such that every vertex in its vertex set ($V_1$ say) except $z$ has an $\mathcal{S}$-ancestor in $A$; and 
let $\mathcal{S}_2$ be a sublevelling
of $\mathcal{S}$ with head $z$ and base $M(B)\setminus M(A)$ 
such that every vertex in its vertex set ($V_2$ say) except $z$ has an $\mathcal{S}$-ancestor 
in $B$ and has no $\mathcal{S}$-ancestor in $A$.
Thus $V_1\cap V_2 = \{z\}$, and no vertex in $V_2$ has a parent in $V_1\setminus \{z\}$.
By \ref{splitshower},
since $\chi(M(A))>\kappa$ and $\chi(M(B)\setminus M(A))>2\kappa$,
there is a
monotone path $R$ of $G[V_1]$ between $z$ and $M(A)$ with the following property. Let $X$ 
denote the set of vertices in $V_2\setminus \{z\}$
that have a neighbour in $V(R)\setminus \{z\}$;
then  the set $Y$ of vertices in $M(B)\setminus M(A)$ with an ancestor in $X$
satisfies 
$$(\nu+1)(3\nu^2+1)\chi(Y) \ge \chi(M(B)\setminus M(A)) - 2\kappa>(\nu+1)(3\nu^2+1)(\tau+\kappa),$$
and consequently $\chi(Y)>\tau+\kappa$.

Now no vertex of $R$ different from $z$ belongs to or has a child in $V_2$, and so, since $X\subseteq V_2$, 
every vertex in $X$ has a child in $V(R)$. 
Let $y$ be the vertex of $R$ with height two, and let $R'$ be the subpath of $R$ between $z,y$.
Let $X_1$ be the set of vertices in $X$ with a child in $R'$, and let $X_2$
be the set of vertices in $X$ with a child in $R$ with height at most one. Thus $X=X_1\cup X_2$.
Let $P$ be the union of $R'$ and a monotone path between $L_0$ and $z$. The vertices of $P$ in order form a wand, and
every vertex in $X_1$ is an up-neighbour of a vertex of this wand.
Consequently the set of $\mathcal{S}_2$-descendants in $M$ of $X_1$ is a subset of the shadow in $\mathcal{S}$ over $M$ of this wand, and so has chromatic number
at most $\tau$. But every vertex in $M$ with an ancestor in $X_2$ is at distance at most three from the penultimate 
vertex of $R$, and in particular the set of descendants in $M$ of $X_2$ has chromatic number at most $\kappa$.
Consequently $\chi(M(X))\le \tau+\kappa$, a contradiction since $Y\subseteq M(X)$. This proves \ref{nosplit}.~\bbox

Let $\mathcal{S}$ be a stable shower in a candidate $G$, and let $\xi\ge 0$ be an integer.
A wand $\mathcal{W}=(W_0\l W_t)$ is said to be  {\em $\xi$-diagonal} if
\begin{itemize}
\item every vertex of $U(\mathcal{S})$ with a child in $V(\mathcal{W})$ belongs to $V(\mathcal{W})$; and
\item for $0\le i\le t$, the set of vertices in $M$ that have an ancestor in $W_i$ and no ancestor in $W_{i+1}$
has chromatic number at most $\xi$ (where $W_{t+1}=\emptyset$).
\end{itemize}
Next we need some results about showers that admits $\xi$-diagonal wands, where $\xi$ is bounded. Before we do so, let us 
set up some notation for these things.

If $\mathcal{S}=(L_0\l L_k, s)$ with mat $M$, and $\mathcal{W}$ is a $\xi$-diagonal wand $(W_0\l W_{t})$ in $\mathcal{S}$,
then for every vertex $v$ of $U(\mathcal{S})\cup M$, there is a maximum $i\le t$ such that $W_i$ contains an ancestor of $v$.
We call this number $i$ the {\em reach} of $v$ (with respect to $\mathcal{W}$). 
Let $U= U(\mathcal{S})$, and for $0\le i\le t$ let $M_i$
and $U_i$ be the sets of all vertices with reach $i$ in $M$ and in $U$, respectively.  
It follows that no member of $U_j$ has a child in $U_i$ if $i<j$.
Let $M_i=U_i = W_i = \emptyset$ for $t+1\le i\le k-2$.

\begin{thm}\label{diagonal0}
Let $\mathcal{S}$ be a stable shower with mat $M$ in a candidate $G$, such that
the shadow over $M$ of every wand in $\mathcal{S}$ has chromatic number
at most $\tau$. Let $\mathcal{W}= (W_0\l W_t)$ be a $\xi$-diagonal wand, and
let $P$ be a monotone path between $M$ and $V(\mathcal{W})$,
with no vertex in $V(\mathcal{W})$ except one end. Let $0\le a\le t$, and let
$X= \cup_{0\le i<a} (U_i\cup M_i)$ and $Y= \cup_{a< i\le t} (U_i\cup M_i)$.
Suppose that $V(P)\subseteq Y$.
Let $X(P)$ be the set of vertices in $X\setminus V(\mathcal{W})$ with a neighbour in $V(P)$.
Then the set of vertices in $M\cap X$ with an ancestor in $X(P)$ has chromatic number at most $\tau+\kappa$.
\end{thm}
\Proof  Let $P$ have vertices $p_h\c p_k$ in order, where $p_i\in L_i$ for $h\le i\le k$, and $p_h\in W_h$, and therefore
$a< h$.
Let $Z_1,Z_2$ respectively be the sets of all $v\in X(P)$, such that
\begin{itemize}
\item $v$ has a neighbour in $\{p_h\l p_{k-2}$\};
\item $v$ has a neighbour in $\{p_k,p_{k-1}\}$.
\end{itemize}

If $v\in Z_1$, then $v$ has no parent in $V(P)$ from the definition of ``reach''. Suppose that $v$ has a parent in
$W_0\cup\cdots\cup W_{h-1}$. Then since $v\in N(P)$, it has a neighbour in one of $W_{h-1}, W_{h-2}$. 
But $v$ is not a parent of $p_h$
since $v\notin V(\mathcal{W})$, so $v$ has no neighbour in $W_{h-2}$. Thus $v$ has a parent in $W_{h-1}$, and so
has reach at least $h-1$. But $h-1\ge a$,
contradicting that $v\in Z_1$. This proves that every vertex in $Z_1$ is an up-neighbour of the wand
$$(W_0\l W_{h-1}, \{p_h\},\{p_{h+1}\}\l \{p_{k-2}\}).$$ Moreover, if $R$ is a monotone path between some $v\in Z_1$ and $M\cap X$,
then $V(R)\subseteq X$ from the definition of ``reach'', and so no vertex of $R$ has a parent in this wand.
Consequently every vertex in $M\cap X$ with an ancestor in $Z_1$ belongs to the shadow in $\mathcal{S}$ of this wand over $M$, 
and so the set of such
vertices has chromatic number at most $\tau$.

If $v\in Z_2$ then $v$ has height at most two, and so every descendant of $Z_2$ in $M$ has distance at most three from $p_{k-1}$.
Since $\rho>3$ it follows that the set of such descendants has chromatic number at most $\kappa$.
Summing, this proves \ref{diagonal0}.~\bbox

\bigskip

A monotone path is {\em vertical} if for some $i$, all its vertices belong to $M_i\cup U_i$.
Note that, if $P$ is a monotone path between some vertex in $M_h$ and some vertex in $W_h$, then 
$P$ is vertical. If $X\subseteq U\cup M$, the set of vertices in $M$ joined to a vertex in $X$
by a vertical path is denoted by $X\downarrow M$. The previous result \ref{diagonal0} told us about the chromatic number 
of the descendants
of vertices with neighbours in a monotone path, when we confine ourselves to vertices with smaller reach 
than the vertices of the path (actually,
reach smaller by at least two). The
next result does the same when we confine ourselves to larger reach; except we can only handle descendants reachable by vertical paths,
not general descendants.

\begin{thm}\label{diagonal1}
Let $\zeta$ be a sidekick for $\sigma$.
Let $\mathcal{S}$ be a stable shower with mat $M$ in a candidate $G$, such that
$\mathcal{S}$ is not $(\sigma+1)$-complete over $M$, and the shadow over $M$ of every wand in $\mathcal{S}$ has chromatic number
at most $\tau$. Let $\mathcal{W}= (W_0\l W_t)$ be a $\xi$-diagonal wand, and
let $P$ be a monotone path between $M$ and $V(\mathcal{W})$,
with no vertex in $V(\mathcal{W})$ except one end. Let $0\le a\le t$, and let
$X= \cup_{0\le i<a} (U_i\cup M_i)$ and $Y= \cup_{a< i\le t} (U_i\cup M_i)$.
Suppose that with notation as above, $V(P)\subseteq X$.
Let $Y(P)$ be the set of vertices in $Y\setminus V(\mathcal{W})$ with a neighbour in $V(P)$.
Then $\chi(Y(P)\downarrow M)\le 2\zeta + 2\xi + \kappa.$
\end{thm}
\Proof  Let $P$ have vertices $p_h\c p_k$ in order, where $p_i\in L_i$ for $h\le i\le k$, and $p_h\in W_h$, and therefore 
$h<a$. 

Let $M_0=Y(P)\cap M$. Every vertex in $M_0$ is adjacent to one of $p_{k-1}, p_k$, so $\chi(M_0)\le \kappa$.
We may therefore assume that there are vertices in $Y(P)\cap U$ with reach greater than $a$,
and so there exists $i\in \{h\l k\}$ such that some neighbour
$y$ of $p_i$ belongs to $Y(P)\cap U$ and has reach greater than $a$. Choose $i$ minimum with this property. 
Now $y$ is not a parent of $p_i$ from the definition of ``reach'', and since $y\in U$ it follows that $i<k$ and $p_i\in U$.
Consequently $y$ is a child of $p_i$, and so $i\le k-2$. Let $y\in U_{j_1}$; and we may assume that $y$ is chosen with $j_1$
maximum. The height of $y$ is at most $k-j_1-1$, and so the height of
$p_i$ is at most $k-j_1$, that is, $i\ge j_1$. In particular, since $j_1>a$ and $h<a$, it follows that $i\ge h+2$.
If possible, let $j_2\in \{h+2\l t\}$ be maximum such that $p_{i+1}$ has a child in 
$U_{j_2}\setminus W_{j_2}$, and otherwise $j_2$ is undefined.

Let $Q$ be a vertical path between $y$ and $W_{j_1}$, and let $y'$
be the neighbour of $y$ in $Q$. Then $p_i$ has no neighbour in $V(Q)$ except $y$. Let $\mathcal{S}_1$ be the maximal 
sublevelling of $\mathcal{S}$ with head $p_i$ and with base a subset of $M$, 
such that no child of $y$ or $y'$ belongs to $U(\mathcal{S}_1)$.
Let $M^1$ be its base.
For $0\le j\le j_1$, let $c_j\in W_j$, where $c_h = p_h$ and $c_{j_1}\in V(Q)$.
Consequently $c_0\c c_h\d p_{h+1}\c p_i$
and $c_0\c c_j\d Q\d y\d p_i$ are both induced paths, and so the
pair forms a switch for $\mathcal{S}_1$. From \ref{noswitch}, it follows that $\chi(M^1)\le \zeta$.

If $j_2$ is defined let $y''$ be a child of $p_{i+1}$ in $U_{j_2}$, and let $M^2$ be the set of vertices $v\in M$
such that there is a monotone path of $\mathcal{S}$ between $v,p_{i+1}$ containing no child of $y''$ or of an appropriate parent of $y''$;
then similarly, $\chi(M^2)\le \zeta$. If $j_2$ is undefined let $M^2=\emptyset$.

Let $M^3$ be the set of all $v\in (Y(P)\cap U)\downarrow M$ such that $v\notin M^1\cup M^2$. Let $v\in M^3$ and let $R$
be a vertical path between $v$ and $u\in Y(P)\cap U$ say. Thus $u$ is therefore
a child of $p_{i'}$ for some $i'$ with $i\le i'\le k-2$.
By adding the edge $up_{i'}$ and the path $p_i\c p_{i'}$, we obtain a monotone path between $p_i$ and $v$. Since $v\notin M^1$,
this path contains a child of one of $y,y'$.
Now no vertex of $P$ is a child of $y$ or $y'$ by the definition of ``reach'', and so
this child belongs to $R$; and since $y,y'\in L_i\cup L_{i+1}$, some vertex of $R$ belongs to $L_{i+2}$,
and therefore $i'\le i+1$. On the other hand, $i'\ge i$; so there are two cases, $i'=i$ and $i'=i+1$.
Choose $j$ with $v\in M_{j}$; then $V(R)\subseteq U_{j}\cup M_{j}$, and again, from the definition of ``reach'',
it follows that $j\ge j_1$. Suppose first that $i'=i$; then $j=j_1$ from the choice of $j_1$, and so $v\in M_{j_1}$.
Similarly, if $i'=i+1$, then since $v\notin M^2$, it follows that $v\in M_{j_2}$.
We have shown then that $M^3\subseteq M_{j_1} \cup M_{j_2}$, and so $\chi(M^3)\le 2\xi$.

Now let $v\in Y(P)\downarrow M$, and let $R$ be a vertical path between $v$ and some $u\in Y(P)$.
If $u\in M$ then $u=v$ and $v\in M_0$. If $u\in U$, then $v$ belongs to one of $M^1,M^2,M^3$.
Consequently $\chi(Y(P)\downarrow M)\le \kappa + 2\zeta + 2\xi$. This proves
\ref{diagonal1}.~\bbox

\begin{thm}\label{diagonal2}
Let $\zeta$ be a sidekick for $\sigma$.
Let $\mathcal{S}$ be a stable shower with mat $M$ in a candidate $G$, such that
$\mathcal{S}$ is not $(\sigma+1)$-complete over $M$, and the shadow over $M$ of every wand in $\mathcal{S}$ has chromatic number
at most $\tau$. Let $\mathcal{W}$ be a $\xi$-diagonal wand. In the usual notation, 
let $h<j\le t$, and let $H\subseteq \bigcup_{h< i<j}M_i$ such that $\chi(H)>2\xi+2\kappa+\tau$.
Let $c_h\in W_h$, and $c_j\in W_j$. 
Then there is a set $\mathcal{A}$ of integers, and for each $a\in \mathcal{A}$ there is an induced path $J_a$
of $G$ between $c_h, c_j$, with the following properties:
\begin{itemize}
\item $\mathcal{A}$ has cardinality at most $\nu+1$, and includes a dense set of cardinality $\nu$,
and contains two integers $x,y$ with $y-x\in \{1,3\}$;
\item $|E(J_a)|=a$ for each $a\in \mathcal{A}$;
\item for each $a\in \mathcal{A}$, 
$V(J_a)\subseteq \{c_h,c_j\}\cup U_{h+1}\cup \cdots\cup U_{j-1}\cup H$;
\item for each $a\in \mathcal{A}$, every vertex of $J_a$ belongs either to $V(\mathcal{W})\cup V(H)$
or to a vertical path with one end in $H$; and
\item for each $a\in \mathcal{A}$, there is a set of at most $3\nu^2+2$ $\mathcal{S}$-monotone paths,
each with vertex set a subset of $ \{c_h\}\cup U_{h+1}\cup \cdots\cup U_{j-1}$, 
such that every vertex of $V(J_a)\setminus (H\cup V(\mathcal{W}))$ belongs to one of these paths.
\end{itemize}
\end{thm}
\Proof
We may assume that $H$ is connected, by replacing it by one of its components with maximum chromatic number.
No vertex in $H$ has an ancestor in $W_j$; choose $i<j$ maximum such that $H\cap M_i\ne \emptyset$.
Thus $H\subseteq M_{h+1}\cup\cdots\cup M_i$.
Choose $c_i\in W_i$ with a descendant in $M_i\cap H$, and let 
$Q$ be a vertical path between $c_i$ and $M_i\cap H$. 
Let $N(Q)$ denote the set of vertices in $M\cup (U\setminus V(\mathcal{W}))$ with a neighbour in $V(Q)$.
By \ref{diagonal0}, $N(Q)\downarrow (H\setminus (M_{i-1}\cup M_i))$
has chromatic number at most $\tau+\kappa$, and since $\chi(M_{i-1}\cup M_i)\le 2\xi$, it follows that
there exists $H'\subseteq H\setminus (M_{i-1}\cup M_i)$ such that
$$\chi(H')\ge \chi(H)-(2\xi+\kappa+\tau) > \kappa,$$ 
and no vertical path meets both $H'$ and $N(Q)$.
Thus $H'\subseteq M_{h+1}\cup\cdots\cup M_{i-2}$, and in particular $i\ge h+3$.
Let $\mathcal{S}=(L_0\l L_k,s)$.
Let $X$ be the union of the vertex sets of all vertical paths with one end in $H'$ and the other in $V(\mathcal{W})$,
together with 
$$\{c_h\}\cup W_{h+1}\cup\cdots\cup W_{i-2},$$
and for $h\le j' < k$ let $L_{j'} = L_j\cap X$. Let $L_k'$ be the union of $H$, $V(Q)$, and $W_{i+1}\cup W_{i+2}\cup\cdots\cup W_{j}$.
Then $G[L_k']$ is connected and 
every vertex of $X\setminus L_k$ with a neighbour in $L_k'$ belongs to $L_{k-1}'$. It follows that
$(L_h'\l L_{k-1}', L_k', c_j)$ is a stable shower $\mathcal{S}'$, with mat $H'$; and the result follows from \ref{staircase} and \ref{wastecover}.
(Note: \ref{wastecover}
gives us $3\nu^2+1$ $\mathcal{S}'$-monotone paths containing all the vertices of $J_a$ not in $L_k'$. We can assume 
that none of these paths has
a vertex in $L_k'$, and so they are also $\mathcal{S}$-monotone; but we also need to cover the vertices
of $J_a$ in $L_k'\setminus (H\cup V(\mathcal{W}))$. One more $\mathcal{S}$-monotone path will do this, namely $Q$.) 
This proves \ref{diagonal2}.~\bbox

\begin{thm}\label{diagonal3}
Let $\zeta$ be a sidekick for $\sigma$.
Let $\mathcal{S}$ be a stable shower with mat $M$ in a candidate $G$, such that
$\mathcal{S}$ is not $(\sigma+1)$-complete over $M$, and the shadow over $M$ of every wand in $\mathcal{S}$ has chromatic number
at most $\tau$. Let $\mathcal{W}$ be a $\xi$-diagonal wand. With the usual notation, let $j_0<j_1<j_2\le t$, 
and suppose that $u_1\in M_{j_0}$ and $u_2\in M_{j_2}$
are adjacent. Let $M^1\subseteq \bigcup_{j_0<j<j_1}M_j$ and $M^2\subseteq\bigcup_{j_1< j<j_2}M_j$. If 
$$\chi(M^1)>2\zeta + 5\xi+ 4\kappa+ 2\tau$$
and 
$$\chi(M^2)> ((\nu+1)(3\nu^2+2)+1)(2\zeta+2\xi+\kappa)+ 3\xi+3\kappa+2\tau$$
then there is an edge between $M^1,M^2$.
\end{thm}
\Proof
Let $P_1$ be a vertical path between $u_1,W_{j_0}$, and let $P_2$ be a vertical path between $u_2, W_{j_2}$. 
Let $c_{j_0}$ be the end of $P_1$ in $W_{j_0}$, and let $c_{j_2}$ be the end of $P_2$ in $W_{j_2}$.
Since $u_1,u_2$ are adjacent, there is an induced path $P$ between $c_{j_0}, c_{j_2}$ with $V(P)\subseteq V(P_1\cup P_2)$.
For $i = 1,2$, let $N(P_i)$ be
the set of vertices in $M\cup (U\setminus V(\mathcal{W}))$ with a neighbour in $V(P_i)$.
By \ref{diagonal0},
$\chi(N(P_2)\downarrow M^1)\le \kappa+\tau$. Moreover, by \ref{diagonal1},
$$\chi(N(P_1)\downarrow (M^1\setminus M_{j_1+1}))\le 2\zeta+2\xi+\kappa,$$
and since $\chi(M_{j_1+1})\le \xi$, it follows that
$\chi(N(P_1)\downarrow M^1)\le 2\zeta+3\xi+\kappa$.
Consequently there exists $H_1\subseteq M^1$ with 
$$\chi(H_1)>\chi(M^1)- (2\zeta+3\xi+2\kappa+\tau)\ge 2\xi+2\kappa+\tau,$$ 
such that
no vertex in $H_1$ belongs to a vertical path that intersects
$N(P_1)\cup N(P_2)$. Choose $c_{j_1}\in W_{j_1}$.
By \ref{diagonal2}, 
\\
\\
(1) {\em There is a set $\mathcal{A}$ of integers, and for each $a\in \mathcal{A}$ there is an induced path $J_a$
of $G$ between $c_{j_0}, c_{j_1}$, with the following properties:
\begin{itemize}
\item $\mathcal{A}$ has cardinality at most $\nu+1$, and includes a dense set of cardinality $\nu$,
and contains two integers $x,y$ with $y-x\in \{1,3\}$;
\item $|E(J_a)|=a$ for each $a\in \mathcal{A}$;
\item for each $a\in \mathcal{A}$,
$V(J_a)\subseteq \{c_{j_0},c_{j_1}\}\cup U_{j_0+1}\cup \cdots\cup U_{j_1-1}\cup H_1$;
\item for each $a\in \mathcal{A}$, every vertex of $J_a$ belongs either to $V(\mathcal{W})\cup V(H_1)$
or to a vertical path with one end in $H_1$; and
\item for each $a\in \mathcal{A}$, there is a set of at most $3\nu^2+2$ $\mathcal{S}$-monotone paths,
each with vertex set a subset of $ \{c_{j_0}\}\cup U_{j_0+1}\cup \cdots\cup U_{j_1-1}$,
such that every vertex of $V(J_a)\setminus (H_1\cup V(\mathcal{W}))$ belongs to one of these paths.
\end{itemize}
}
In particular, for each $a\in \mathcal{A}$, $P\cup J_a$ is an induced path between $c_{j_1}$ and $c_{j_2}$, because of the 
fourth bullet above and from the choice of $H_1$.
Now suppose that there are no edges between $M^1,M^2$.
By $(\nu+1)(3\nu^2+2)+1$ applications of \ref{diagonal1}, and one application of \ref{diagonal0},
there exists $H_2\subseteq M_2$ with the following properties:
\begin{itemize}
\item no vertex in $H_2$ belongs to a vertical path that intersects $N(P_1)\cup N(P_2)$;
\item for each $a\in \mathcal{A}$, no vertex in $H_2$ belongs to a vertical path that contains a vertex in $V(J_a)$ or a neighbour
of such a vertex (here we use that there is no edge between $H_1$ and $M^2$); and
\item $\chi(H_2)\ge \chi(M^2)-((\nu+1)(3\nu^2+2)+1)(2\zeta+2\xi+\kappa)- (\kappa+\tau+\xi)> 2\xi+2\kappa+\tau.$
\end{itemize}

We apply \ref{diagonal2} to $H_2$, and thereby obtain a set of paths joining $c_{j_1}$ and $c_{j_2}$.
More precisely:
\\
\\
(2) {\em There is a set $\mathcal{B}$ of integers, and for each $b\in \mathcal{B}$ there is an induced path $K_{b}$
of $G$ between $c_{j_1}, c_{j_2}$, with the following properties:
\begin{itemize}
\item $\mathcal{B}$ has cardinality at most $\nu+1$, and includes a dense set of cardinality $\nu$,
and contains two integers $x,y$ with $y-x\in \{1,3\}$;
\item $|E(K_{b})|=b$ for each $b\in \mathcal{B}$;
\item for each $b\in \mathcal{B}$,
$V(K_{b})\subseteq \{c_{j_1},c_{j_2}\}\cup U_{j_1+1}\cup \cdots\cup U_{j_2-1}\cup H_2$; and
\item for each $b\in \mathcal{B}$, every vertex of $K_b$ belongs either to $V(\mathcal{W})\cup V(H_2)$
or to a vertical path with one end in $H_2$.
\end{itemize}
}
For each $b\in \mathcal{B}$ and each $a\in \mathcal{A}$, it follows from the fourth bullet of (2) 
and the choice of $H_2$ that $J_a\cup K_b\cup P$ is a hole.
It follows as usual that $G$ contains a hole $\nu$-interval,
a contradiction. This proves \ref{diagonal3}.~\bbox

\begin{thm}\label{diagonal4}
Let $\zeta$ be a sidekick for $\sigma$.
Let $\mathcal{S}$ be a stable shower with mat $M$ in a candidate $G$, such that
$\mathcal{S}$ is not $(\sigma+1)$-complete over $M$, and the shadow over $M$ of every wand in $\mathcal{S}$ has chromatic number
at most $\tau$. Let $\mathcal{W}$ be a $\xi$-diagonal wand. In the usual notation, let $j_1<j_2\le t$, and suppose that 
$u_1\in M_{j_1}$ and $u_2\in M_{j_2}$
are adjacent. Let $M^1\subseteq \bigcup_{j_1<j<j_2}M_j$ and $M^2\subseteq\bigcup_{j_2< j\le t}M_j$. If
$$\chi(M^1)> (2\zeta+4\xi+\tau+4\kappa) + (\nu+1)(3\nu^2+1)(\tau+\kappa) $$
and
$$\chi(M^2)>4\zeta+5\xi+3\kappa,$$
then there exist $H_1\subseteq M^1$ and $H_2\subseteq M^2$ such that 
$\chi(H_1) \ge \chi(M_1)-(2\zeta+4\xi+\tau+2\kappa)$ and $\chi(H_2) \ge \chi(M^2) - (4\zeta+5\xi+2\kappa)$
and 
there is no edge between $H_1, H_2$.
\end{thm}
\Proof
For $i = 1,2$, let $P_i$ be a vertical path between $u_i$ and some $c_{j_i}\in W_{j_i}$. Let $P$ be an induced path between
$c_{j_1}, c_{j_2}$ with $V(P)\subseteq V(P_1\cup P_2)$. 
For $i = 1,2$, let $N(P_i)$
be the set of vertices in $M\cup (U\setminus V(\mathcal{W}))$ with a neighbour in $V(P_i)$.

Let $B$ be the set of all vertices that belong to a vertical path $R$ between $M^1\cup M^2$ 
and $V(\mathcal{W})$ such that no
vertex of $R$ belongs to $N(P_1)\cup N(P_2)$. Consequently there are no edges between 
$V(P)$ and $B$. Moreover, there are no edges between the interior of $P$ and $V(\mathcal{W})\setminus (W_{j_1}\cup W_{j_2})$.

By \ref{diagonal1}, 
$\chi(N(P_1)\downarrow (M_1\setminus M_{j_1+1}))\le 2\zeta+2\xi+\kappa,$
and so 
$$\chi(N(P_1)\downarrow M_1)\le 2\zeta+3\xi+\kappa.$$ 
Also,
from \ref{diagonal0},  
$\chi(N(P_2)\downarrow (M_1\setminus M_{j_2-1})\le \tau+\kappa$,
and so 
$$\chi(N(P_2)\downarrow M_1)\le \tau+\xi+\kappa.$$
Consequently
$$\chi(B\cap M_1)>\chi(M_1)-(2\zeta+4\xi+\tau+2\kappa).$$
Choose $H_1\subseteq B\cap M_1$, such that $G[H_1]$ is connected and
$\chi(H_1)=\chi(B\cap M_1)$.
Similarly, we may choose $H_2\subseteq B\cap M_2$ such that
$G[H_2]$ is connected and $\chi(H_2)>\chi(M_2)-(4\zeta+5\xi+2\kappa)$.
For $i = 1,2$, let $B_i$ be the set of vertices in $B$ that belong to a vertical path with one end in $H_i$.

Suppose that there is an edge between $H_1,H_2$, and so $G[H_1\cup H_2]$ is connected.
Let $\mathcal{S}=(L_0\l L_k, s)$. Then $\mathcal{S}' = (L_0\l L_{k-1}, H_1\cup H_2, s')$ 
is also a shower (where $s'\in H_1\cup H_2$ is arbitrary).
We need to define two sublevellings of $\mathcal{S}'$.
\begin{itemize}
\item Let $L_{j_1}^1=\{c_{j_1}\}$, for $j_1< i<j_2$ let $L_i^1=W_i\cup (L_i\cap B_1)$, and for $j_2\le i\le k$ let $L_i^1=L_i\cap B_1$;
then $(L_{j_1}^1\l L_k^1)$ is a sublevelling $\mathcal{S}_1$ of $\mathcal{S}'$ with head $c_{j_1}$ and base $H_1$.
\item 
Let $L_{j_2}^2=\{c_{j_2}\}$, for $j_2< i\le t$ let $L_i^2=W_i\cup (L_i\cap B_2)$, and for $t < i\le k$ let $L_i^2=L_i\cap B_2$;
then $(L_{j_2}^2\l L_k^2)$ is a sublevelling $\mathcal{S}_2$ of $\mathcal{S}'$ with head $c_{j_2}$ and base $H_2$.
\end{itemize}
In particular, there are no edges between the interior of $P$ and $V(\mathcal{S}_i)$ for $i = 1,2$.

Let us apply \ref{splitshower} to the pair $\mathcal{S}_2,\mathcal{S}_1$ of sublevellings of $\mathcal{S}'$ (in this order). 
Since $\chi(H_2)> \kappa$ and $\chi(H_1)>2\kappa$, and the base of $\mathcal{S}'$ is the union of the bases of $\mathcal{S}_1$
and $\mathcal{S}_2$,
we deduce that either
\begin{itemize}
\item there are $\nu$ induced paths
$Q_0\l Q_{\nu-1}$ of $G[V(\mathcal{S}_1)\cup V(\mathcal{S}_2)]$ between $c_{j_1}, c_{j_2}$, 
such that $|E(Q_i)|=|E(Q_0)|+i$ for $0\le i<\nu$; or
\item there is an
$\mathcal{S}_2$-monotone path $R$ between $c_{j_2}$ and $H_2$ such that
$$(\nu+1)(3\nu^2+1)\chi(H_1(X(R))) \ge \chi(H_1) - 2\kappa,$$
where
$X(R)$ denotes the set of vertices in $V(\mathcal{S}_1)$
that have a neighbour in $V(R)$, and $H_1(X(R))$ denotes the set of $\mathcal{S}_0$-descendants in $H_1$ of the members
of $X(R)$.
\end{itemize}

Suppose that $Q_0\l Q_{\nu-1}$ are as in the first statement. Let $0\le j\le \nu-1$; we claim that $P\cup Q_j$ is a hole.
Since $P, Q_j$ are induced paths with the same ends, it is enough to show that no vertex of the interior of $P$ belongs to
or has a neighbour in the interior of $Q_j$. Let $q$ belong to the interior of $Q_j$.
Then $q\in V(\mathcal{S}_i)$ for some $i\in \{1,2\}$, and no vertex of the interior of $P$
belongs to or has a neighbour in $V(\mathcal{S}_i)$, as we saw above. Thus $P\cup Q_j$ is a hole for each $j$, 
and these holes form a hole $\nu$-sequence, which is impossible.

Now suppose that $R$ satisfies the second statement. By \ref{diagonal0}, $\chi(H_1(X(R)))\le \tau+\kappa$, and so
$(\nu+1)(3\nu^2+1)(\tau+\kappa) \ge \chi(H_1) - 2\kappa$, a contradiction. 

It follows that there is no edge between $H_1,H_2$. This proves \ref{diagonal4}.~\bbox

\bigskip

We need the following lemma.

\begin{thm}\label{jumps}
Let $G$ be a graph with chromatic number more than $4N$, and let $M_1\l M_k$ be a partition of $V(G)$
such that $\chi(M_i)\le N$ for $1\le i\le k$. Then there exist $a<b<c<d<e\le k$ such that
there is an edge of $G$ between $M_a$ and $M_c$, and an edge between $M_a$ and $M_e$.
\end{thm}
\Proof Let $J$ be the graph with vertex set $\{1\l k\}$ in which $i,j$ are adjacent if there is an edge of $G$
between $M_i$ and $M_j$. If $J$ is 4-colourable, then $\chi(G)\le 4N$, a contradiction. So $J$ is not 4-colourable,
and consequently there exists $a\in \{1\l k\}$ such that $a$ is adjacent in $J$ to at least four of $a+1\l k$.
Let $b,c,d,e$ be four such neighbours, in order; then the theorem holds. This proves \ref{jumps}.~\bbox

\begin{thm}\label{diagonal}
Let $\zeta$ be a sidekick for $\sigma$.
Let $\mathcal{S}$ be a stable shower with mat $M$ in a candidate $G$, such that
$\mathcal{S}$ is not $(\sigma+1)$-complete over $M$, and the shadow over $M$ of every wand in $\mathcal{S}$ has chromatic number
at most $\tau$. Let
$$\eta = ((\nu+1)(3\nu^2+2)+6)(2\zeta+2\xi+\tau+\kappa)+2\tau.$$ 
Let $\mathcal{W}$ be a $\xi$-diagonal wand.
Then $\chi(M)\le 4(\eta+\xi)+\eta$.
\end{thm}
\Proof
Suppose that $\chi(M)>4(\eta+\xi)+\eta$.
Let $\mathcal{W} = (W_0\l W_t)$.
Let $j_0=-1$, and define $j_1,j_2,\ldots j_r$ and $M^1\l M^{r-1}$ inductively as follows. 
Having defined $j_0\l j_i$ and $M^1\l M^{i-1}$, if $\chi(\cup_{j_i<j\le t}M_j)<\eta$
the sequence terminates; define $r=i$.
Otherwise choose $j_{i+1}\le t$ minimum such that 
$\chi(\cup_{j_i<j\le j_{i+1}}M_j)\ge \eta$. Let $M^{i} = \bigcup_{j_i<j\le j_{i+1}}M_j$.

This completes the inductive definition. We see that the sets $M^1\l M^{r-1}$ are disjoint, and their union has chromatic number
at least $\chi(M)-\eta>4(\eta+\xi)$; and each $M_i$ has chromatic number at least $\eta$, and at most $\eta+\xi$ (from the minimality of
$j_{i+1}$). It follows from \ref{jumps} that
there exist $a<b<c<d<e\le r$ such that
there is an edge of $G$ between $M^a$ and $M^c$, and an edge between $M^a$ and $M^e$. 
Now 
$$\chi(M^b)\ge  \eta > (2\zeta+4\xi+\tau+4\kappa) + (\nu+1)(3\nu^2+1)(\tau+\kappa) $$
and 
$$\chi(M^d)\ge  \eta > 4\zeta+5\xi+3\kappa,$$
so by 
\ref{diagonal4} applied to $M^b, M^d$ and the edge between $M^a,M^c$, there exist $H_1\subseteq M^b$
and $H_2\subseteq M^d$ such that 
$\chi(H_1) \ge \eta-(2\zeta+4\xi+\tau+2\kappa)$ and $\chi(H_2) \ge \eta - (4\zeta+5\xi+2\kappa)$,
and there is no edge between $A_1, A_2$.
But since 
$$\chi(H_1)>2\zeta + 5\xi+ 4\kappa+ 2\tau$$
and
$$\chi(H_2)> ((\nu+1)(3\nu^2+2)+1)(2\zeta+2\xi+\kappa)+ 3\xi+3\kappa+2\tau$$
this contradicts \ref{diagonal3} applied to $H_1,H_2$ and the edge between $M^a, M^e$.
This completes the proof of \ref{diagonal}.~\bbox

\bigskip

Now we can prove the objective of this section, the following.

\begin{thm}\label{noknobbly2}
Let $\zeta$ be a sidekick for $\sigma$. Let $N= (\nu+1)(3\nu^2+2)+9$.
Let $\tau\ge 0$, and let
$\mathcal{S}$ be a stable shower with mat $M$ in a candidate $G$, such that
$\mathcal{S}$ is not $(\sigma+1)$-complete over $M$, and the shadow over $M$ of every wand in $\mathcal{S}$ has chromatic number
at most $\tau$.
Then $\chi(M)\le 40 N^2 \kappa+ 40 N\zeta  + 20N\tau$.
\end{thm}
\Proof
Let 
$$\xi =((\nu+1)(3\nu^2+1)+8)\kappa,$$
and 
$$\eta = ((\nu+1)(3\nu^2+2)+6)(2\zeta+2\xi+\tau+\kappa)+2\tau.$$
Let $\mathcal{S}= (L_0\l L_k,s)$,
and for each $v\in U$, let $M(v)$ denote the set of descendants of $v$ in $M$.
Let $z_0\in L_0$, and recursively, having defined $z_i$, let $z_{i+1}$ be a child of $z_i$ chosen such that $\chi(M(z_{i+1}))>\kappa$
if there is such a child; otherwise the definition terminates, when $i = t$ say. Thus $M=M(z_0)$. Note that since $\chi(M(z_{t}))>\kappa$,
it follows that $z_t$ has height more than $\rho$, and in particular $t\le k-2$, so $(\{z_0\}\l \{z_t\})$ is a wand.
\\
\\
(1) {\em For $0\le i<t$,  $\chi(M(z_{i})\setminus M(z_{i+1})\le \xi$, and $\chi(M(z_t))\le \xi$.}
\\
\\
For $0\le i <t$, since $\chi(M(z_{i+1}))>\kappa$, 
\ref{nosplit} implies that 
$$\chi(M(z_i)\setminus M(z_{i+1}))\le (\nu+1)(3\nu^2+1)(\tau+\kappa) +2\kappa\le \xi.$$ 
We claim that $\chi(M(z_t))\le \xi$; for suppose not. Then by \ref{bigson}, there is a child $z$ of $z_t$
such that $\chi(M(z))\ge \chi(M(z_t))-((\nu + 1)(3\nu ^2 + 1) + 7)\kappa>\kappa$, contrary to the maximality of $t$.
This proves (1).

\bigskip 

For each vertex $v\in M$, choose a monotone path $R_v$ between $v$ and some vertex $x_v$, such that $x_v$ has a neighbour in 
$\{z_0\l z_{t}\}$, with minimum length. Thus no vertex of $R_v$ except $x_v$ has a neighbour in $\{z_0\l z_{t}\}$.
Now $x_v$ might have a parent in $\{z_0\l z_{t}\}$, or a child, or both. Let $X_1$ be the set of vertices in 
$U(\mathcal{S})\setminus V(\mathcal{W})$ with a child and
no parent in $\{z_0\l z_{t}\}$; $X_2$ the set with a parent and no child in $\{z_0\l z_{t}\}$; and $X_3$ the set 
with both a parent and a child in $\{z_0\l z_{t}\}$. For $i = 1,2,3$ let $M^i$ be the set of $u\in M$ such that $x_v\in X_i$.
\\
\\
(2) {\em $\chi(M^1)\le \tau$ and $\chi (M^2)\le 4(\eta+\xi)+\eta$.}
\\
\\
Since $(\{z_0\}\l \{z_t\})$ is a wand, and $M^1$ is a subset of its shadow in $\mathcal{S}$  over $M$, 
it follows that $\chi(M^1)\le \tau$.
Let $V'$ be the union of the vertex sets of the paths $R_v\; (v\in M^2)$, together
with $\{z_0\l z_{t}\}$. Thus no vertex in $V'\setminus \{z_1\l z_t\}$ has a child in $\{z_0\l z_{t}\}$. Now
$$(V'\cap L_0, V'\cap L_1\l V'\cap L_{k-1}, L_k, s)$$ 
is a shower $\mathcal{S}'$ say. Since $\mathcal{S}'$
is included in $\mathcal{S}$, with a one-vertex pipe, it follows that $\mathcal{S}'$ is not 
$(\sigma+1)$-complete over $M$. Moreover, $(\{z_0\}\l \{z_t\})$ is a $\xi$-diagonal wand of $\mathcal{S}'$;
and $M^2$ is a mat for it.
From \ref{diagonal}, it follows that $\chi (M^2)\le 4(\eta+\xi)+\eta$.
This proves (2).

\bigskip

It remains then to bound the chromatic number of $M^3$.
Let $V'$ be the union of the vertex sets of the paths $R_v\; (v\in M^3)$, together
with $\{z_0\l z_{t}\}$; and let $\mathcal{S}'$ be the shower
$$(V'\cap L_0, V'\cap L_1 \l V'\cap L_{k-1}', L_k, s).$$
For $1\le i\le t-1$, let $D_i$ be the set of all vertices of $U(\mathcal{S}')$
(including $z_i$) that are adjacent to both $z_{i+1},z_{i-1}$, and let $D_0 = \{z_0\}$ and $D_t=\{z_t\}$. 
(Note that $z_t$ is the only child of $z_{t-1}$ in $U(\mathcal{S}')$).
For $c = 0,1,2$, let $\mathcal{W}_c$ be the sequence $X_0\l X_t$, where $X_i = D_i$ if $i-c$ is divisible by three, 
and $X_i = \{z_i\}$ otherwise.

Thus each $\mathcal{W}_c$ is a wand, and for each $v\in M^3$, $x_v\in V(\mathcal{W}_c)$ for some $c\in \{0,1,2\}$. 
For $c=0,1,2$, let $H_c$ be the set of $v\in M^3$ such that $x_v\in D_i$ for some $i\in \{0\l t\}$ congruent to $c$ modulo three.
Let $c\in \{0,1,2\}$, and let $v\in H_c$. 
Now no vertex of $R_v\setminus \{x_v\}$ has a parent
in $V(\mathcal{W}_c)$, from the minimality of the length of $R_v$, except for the child of $x_v$ in $R_v$; and the 
latter has no child in $V(\mathcal{W}_c)$ since it has no neighbour in $\{z_0\l z_t\}$. Consequently, if some vertex in
$R_v\setminus \{x_v\}$ has a child in $V(\mathcal{W}_c)$, then $v$ belongs to the shadow in $\mathcal{S}$ of the wand 
$\mathcal{W}_c$ in $\mathcal{S}$; and so the set of all such $v$ has chromatic number at most $\tau$.

Finally, the set of $v\in H_c$ such that no vertex in $R_v\setminus \{x_v\}$ has a child in $V(\mathcal{W}_c)$, has
chromatic number at most 
$4(\eta+\xi)+\eta$, by \ref{diagonal}. 
Thus $\chi(H_c) \le \tau+ 4(\eta+\xi)+\eta$; and so $\chi(M^3)\le 3(\tau+ 4(\eta+\xi)+\eta)$. From (2), it follows
that
$$\chi(M)\le \tau+(4(\eta+\xi)+\eta)+3(\tau+(4(\eta+\xi)+\eta))=20\eta+16\xi+4\tau.$$
Now there is some arithmetic to rewrite this bound in terms of $\tau,\kappa,\nu$, which follows.
Since

$$20\eta = 20(N-3)(2\zeta+2\xi+\tau+\kappa)+40\tau,$$
and $\xi\le N\kappa$,
it follows that
\begin{eqnarray*}
\chi(M)&\le& 20\eta+16\xi+4\tau \\ &=& 20(N-3)(2\zeta+2\xi+\tau+\kappa)+ 16\xi+44\tau\\
&\le&  20N(2\zeta+\kappa)+ (40(N-3)+16)\xi +(20(N-3)+44)\tau \\ &\le& 20N(2\zeta+\kappa) +(40N-104)N\kappa+(20N-16)\tau\\
&\le & 40N^2\kappa + 40 N\zeta  + 20 N\tau.
\end{eqnarray*}


\section{Raising a wand}

Now we turn to general showers, in which a wand shadow may have large chromatic number. We will prove that, if there is such a wand, then we can
use it to construct a new shower, still with large $\chi$, in which no wand shadow has large chromatic number, which we have just 
shown to be impossible.
We begin with:

\begin{thm}\label{knoblemma1}
Let $\mathcal{S}=(L_0\l L_k,s)$ be a stable shower in a candidate $G$, with vertex set $V$, and let $\mathcal{W}=(W_0\l W_{t})$ be a wand
in $\mathcal{S}$.
Let $v$ be a vertex of some post, and let $v\in L_i$ say.
Then there are two induced paths $P_1,P_2$ of $G[V]$ between $v$ and $L_0$, such that $|E(P_2)|=|E(P_1)|+2$, 
and for $j\ge i$ every vertex in
$L_j$ that belongs to either of these paths belongs to $W_{i}\cup W_{i+1}\cup \{v\}$.
\end{thm}
\Proof
Let $P$ be a post containing $v$, with top $t\in T_h$ say; thus $h\le i$. Let $P_0$ be the subpath of $P$ between $v,t$.
Let $u\in W_{h+1}$ be adjacent to $t$.
Let
$P_1$ be the union of $P_0$ and a monotone path between $t$ and $L_0$.
Let $P_2$ be the union of $P_0$, the edge $tu$, and a path between $u$ and $W_0$ with one vertex in each of $W_0\l W_{h+1}$.
This proves \ref{knoblemma1}.~\bbox

\begin{thm}\label{knoblemma2}
Let $\zeta$ be a sidekick for $\sigma$.
Let $\mathcal{S}=(L_0\l L_k,s)$ be a stable shower in a candidate $G$, with mat $M$,
such that $\mathcal{S}$ is not $(\sigma+1)$-complete over $M$. Let
$(W_0\l W_{t})$ be a wand $\mathcal{W}$ in $\mathcal{S}$. Let $0\le i\le t-1$, and let $T_i$ be the set of up-neighbours
of vertices in $W_{i+1}$. Let $M'$ be the set of all $v\in M$ that  belong to a post with top in $T_i$.
Then
$$\chi(M')\le \zeta + 2((\nu+1)(3\nu^2+1)+7)\kappa.$$
\end{thm}
\Proof
For $X\subseteq T_i$, and $j\in \{i\l k\}$, let $L_j(X)$ be the set of all vertices in $L_j$ that belong to a post
with top in $X$.
Then 
$$(W_0, W_1\l W_i, W_{i+1}, X, L_{i+1}(X)\l L_{k-1}(X), L_k, s)$$
is a stable shower $\mathcal{S}(X)$ included in $\mathcal{S}$ (with a one-vertex pipe).
Also $M' = M\cap L_k(T_i)$.
We may assume that $\chi(M')>2((\nu+1)(3\nu^2+1)+7)\kappa$, for otherwise the theorem holds.
By \ref{bigson} applied to $\mathcal{S}(T_i)$ (taking $z_1\in W_i$ and $Y=W_{i+1}$)
there exists $u\in W_{i+1}$ such that 
$$\chi(M\cap L_k(X_0))\ge \chi(M')- ((\nu+1)(3\nu^2+1)+7)\kappa,$$
where
$X_0$ is the set of up-neighbours of $u$. By \ref{bigson} applied to $\mathcal{S}(T(X_0))$ (taking $z_1=u$, and
$Y=X_0$) there exists $x\in X_0$ such that, setting $X=\{x\}$, we have
$$\chi(M\cap L_k(X))\ge \chi(M\cap L_k(X_0))- ((\nu+1)(3\nu^2+1)+7)\kappa;$$
and so
$$\chi(M\cap L_k(X))\ge  \chi(M')- 2((\nu+1)(3\nu^2+1)+7)\kappa.$$
Now
$$(X, L_{i+1}(X)\l L_{k-1}(X), L_k, s)$$
is a shower included in $\mathcal{S}$ (with pipe a monotone path between $x$ and $L_0$), and $M\cap L_k(X)$ is a mat for it.
Since every vertex of $\mathcal{S}(X)$ belongs to a post, it follows that no vertex of $\mathcal{S}(X)$ has a parent in $V(\mathcal{W})$,
and so by \ref{knoblemma1} there is a switch for $\mathcal{S}(X)$ in $\mathcal{S}$.
From \ref{noswitch} it follows
that $\chi(M\cap L_k(X))\le \zeta$. We deduce that
$$ \chi(M\cap L_k(X))\le \zeta + 2((\nu+1)(3\nu^2+1)+7)\kappa.$$
This proves \ref{knoblemma2}.~\bbox

Let $T_0\l T_{t-1}, T$ be as before.
For $0\le i\le k$, let $S_i$ be the set of all vertices $v\in L_i$
that belong to a post with top in $T$.
(Thus $S_i\subseteq L_i\setminus V(\mathcal{W})$, and $S_0 = \emptyset$.) If $M$ is a mat for $\mathcal{S}$,
it follows (since $t\le k-2$) that
$$(W_0, W_1, W_2, W_3\cup S_1, W_4\cup S_2 \l W_t\cup S_{t-2}, S_{t-1}\l S_{k-2}, S_{k-1}, L_k,s)$$
is a stable shower $\mathcal{S}'$ included in $\mathcal{S}$; and we say that $\mathcal{S}'$ is obtained from $\mathcal{S}$
by {\em raising} the wand. Moreover, the shadow $M\cap S_k$ is a mat for $\mathcal{S}'$.

\begin{thm}\label{noparent}
Let $\mathcal{S}=(L_0\l L_k,s)$ be a stable shower in a candidate $G$, and let
$(W_0\l W_{t})$ be a wand in $\mathcal{S}$. Let $\mathcal{S}'$ be obtained from $\mathcal{S}$ by raising the wand.
Then for $0\le i\le t$, if $v\in W_i$ and $v$ is an $\mathcal{S}'$-child of $u$ then $i>0$ and $u\in W_{i-1}$.
\end{thm}
\Proof
In the notation given before, since $v\in W_i$ and $v$ is an $\mathcal{S}'$-child of $u$, it follows that
$i>0$ and $u\in W_{i-1}\cup S_{i-3}$, where $S_{-1}, S_{-2} = \emptyset$.
But $S_{i-3}\subseteq L_{i-3}$ and $v\in W_i\subseteq L_i$, so $u\notin S_{i-3}$, and hence $u\in W_{i-1}$.
This proves \ref{noparent}.~\bbox

\begin{thm}\label{knoblemma3}
Let $\zeta$ be a sidekick for $\sigma$.
Let $\mathcal{S}=(L_0\l L_k,s)$ be a stable shower in a candidate $G$, with mat $M$,
such that $\mathcal{S}$ is not $(\sigma+1)$-complete over $M$. Suppose that
$\mathcal{S}$ is obtained from some stable shower $\mathcal{S}_0$ in $G$ with mat $M_0$ by raising some wand, and
$M$ is the shadow over $M_0$ of this  wand. Let
$\mathcal{W}$ be a wand in $\mathcal{S}$. Then the
shadow $M'$ of $\mathcal{W}$ in $\mathcal{S}$ over $M$ has chromatic number at most
$$3\zeta + 6((\nu+1)(3\nu^2+1)+7)\kappa.$$
\end{thm}
\Proof
Let $\mathcal{W} = (W_0\l W_{t})$,
and for $0\le i\le t-1$, let $T_i$ be the set of up-neighbours
of vertices in $W_{i+1}$ and let $T = T_0\cup\cdots\cup T_{t-1}$.
Thus $M'$ is the set of all $v\in M$ that belong to a post with top in $T$.
Choose $h$ minimum such that $T_h\ne \emptyset$.
Let $M_1,M_2$ be the sets of vertices in $M$ that belong to posts with top in $T_h\cup T_{h+1}$ and with top in
$T\setminus (T_h\cup T_{h+1})$ respectively.
In view of \ref{knoblemma2} it suffices to bound $\chi(M_2)$. For $j=h+2\l k$ let $S_j$ be the set of vertices in $L_j$
that belong to a post with top in $T\setminus (T_h\cup T_{h+1})$. Thus every vertex of every such post belongs to $S_j$
for some $j$. Choose $u\in W_{h+1}$ with a neighbour $v\in T_h$. Consequently
$$(\{u\}, W_{h+2}, W_{h+3}, W_{h+4}\cup S_{h+2}, W_{h+5}\cup S_{h+3}\l W_{t}\cup S_{t-2}, S_{t-1}\l  S_{k-1}, L_k,s)$$
is a shower $\mathcal{S}'$, and $M_2$ is a mat for it. Every vertex of $U(\mathcal{S}')$ belongs to
$L_j$ for some $j\ge h+2$, except $u$. We claim there is a switch for this shower; but in $\mathcal{S}_0$, not in $\mathcal{S}$.

Let $\mathcal{S}_0$ be $(J_0\l J_{k-3}, L_k, s)$. Now $\mathcal{S}$ is obtained from $\mathcal{S}_0$
by raising some wand $\mathcal{D}$ say, where $M$ is the shadow of $\mathcal{D}$ on some mat $M_0$ for
$\mathcal{S}_0$. Let $\mathcal{D}$ be $(D_0\l D_{r})$, and define $D_i=\emptyset$ for $i>r$;
then for $0\le i\le t$, $L_i\subseteq D_i\cup (J_{i-2}\setminus V(\mathcal{D}))$
(where $J_{-1}, J_{-2}=\emptyset$).

Suppose that $u\in V(\mathcal{D})$; then since $u\in L_{h+1}$,  it
follows that $u\in D_{h+1}$. Every vertex of $W_{h-1}$ has distance two from $u$, and so $W_{h-1}\cap J_{h-3} = \emptyset$; so
$W_{h-1}\subseteq D_{h-1}$, since
$$W_{h-1}\subseteq L_{h-1}\subseteq D_{h-1}\cup J_{h-3}.$$
Since $v$ has no neighbour in $W_{h-1}$, and every vertex of $D_h$ is adjacent to every vertex of $D_{h-1}$, it follows that $v\notin D_h$.
But this contradicts \ref{noparent}, since
$v$ is an $\mathcal{S}$-parent of $u$.

This proves that $u\notin V(\mathcal{D})$. Since $u\in W_{h+1}\subseteq L_{h+1}$, it follows that $u\in J_{h-1}$.
By \ref{knoblemma1} applied to $\mathcal{S}_0$, there are two induced paths $P_1,P_2$ of $G$ between $u$ and $L_0$, such that $|E(P_2)|=|E(P_1)|+2$,
and for $j\ge h-1$
every vertex in
$J_j$ that belongs to either of these paths belongs to $D_{h-1}\cup D_{h}\cup \{u\}$.
Suppose that some vertex $x\in V(\mathcal{S}')$ has a neighbour $y$ in one of $P_1,P_2$ where $x,y\ne u$. Let $x\in L_j$;
then $j\ge h+2$. Now $L_j\subseteq D_j\cup (J_{j-2}\setminus V(\mathcal{D}))$. If $x\in D_j$ then $y\in J_i$ for some
$i\ge j-1\ge h+1$,
contradicting that $y\in V(P_1\cup P_2)$. So $x\in J_{j-2}\setminus V(\mathcal{D})$, and so $y\in J_i$ where $i\ge j-3\ge h-1$.
Consequently $y\in D_{h-1}\cup D_{h}\cup \{u\}$, and $y$ is an $\mathcal{S}_0$-parent of $x$. But this is impossible since
$x \in V(\mathcal{S})\setminus V(\mathcal{D})$ and therefore belongs to a post in $\mathcal{S_0}$ for $\mathcal{D}$.

Thus there is no such $x$, and so
$(P_1,P_2)$ is a switch for $\mathcal{S'}$ in
$\mathcal{S}_0$. Hence by \ref{noswitch}, $\chi(M_2)\le \zeta$. Since two applications of \ref{knoblemma2} imply that
$$\chi(M_1)\le 2\zeta + 4((\nu+1)(3\nu^2+1)+7)\kappa,$$
it follows that
$$\chi(M')\le 3\zeta + 4((\nu+1)(3\nu^2+1)+7)\kappa.$$
This proves \ref{knoblemma3}.~\bbox

\begin{thm}\label{knoblemma4}
Let $\zeta$ be a sidekick for $\sigma$.
Let $\mathcal{S}=(L_0\l L_k,s)$ be a stable shower in a candidate $G$, with mat $M$,
such that $\mathcal{S}$ is not $(\sigma+1)$-complete over $M$. 
Let $N= (3\nu^2+2)(\nu+1)+9$.
Then 
$\chi(M)\le 1000N^3\kappa + 1000N^2\zeta$.
\end{thm}
\Proof
Let $\tau= (40 N+176 ) N\kappa+ (40 N+132)\zeta$.
Let $\mathcal{W}$ be a wand in $\mathcal{S}$, let $M'$ be its shadow over $M$, and let $\mathcal{S}'$
be obtained by raising $\mathcal{W}$. Every jet of $\mathcal{S'}$ is a jet of $\mathcal{S}$,
and so $\mathcal{S'}$ is not $(\sigma+1)$-complete. By \ref{knoblemma3}, the shadow over $M'$ of every wand in $\mathcal{S'}$
has chromatic number at most 
$$3\zeta + 4(N-\nu+1)\kappa.$$
By \ref{noknobbly2} applied to $\mathcal{S}'$, it follows that
$\chi(M')\le 40 N^2 \kappa+ 40 N\zeta  + 20N(3\zeta + 4(N-\nu+1)\kappa)\le \tau.$
Thus every wand in $\mathcal{S}$ has shadow over $M$ with chromatic number at most $\tau$; and so 
by another application of \ref{noknobbly2}, 
$\chi(M)\le 40 N^2 \kappa+ 40 N\zeta  + 20N\tau$, and the result follows on substituting for $\tau$.
This proves \ref{knoblemma4}.~\bbox

Let us put these pieces together, to prove \ref{localbound} and hence \ref{bigrad}, in the following strengthened form.

\begin{thm}\label{localbound2}
Let $\nu\ge 2$ and $\kappa\ge 0$ be integers. Let 
$N= (3\nu^2+2)(\nu+1)+9$,
$\zeta_1 = \kappa$, and for $1\le \sigma< \nu$
define
$$\zeta_{\sigma+1}=1000N^2\zeta_{\sigma}+  1000N^3\kappa.$$
Let $G$ be a triangle-free graph such that
$\chi^{\rho}(G)\le \kappa$, where $\rho = 3^{\nu+2}+4$.
If $G$ admits no hole $\nu$-interval then $\chi(G)\le 44 \nu(\kappa+\zeta_{\nu})^{(\nu+1)^2} + 4\kappa$.
\end{thm}
\Proof
By \ref{staircase}, $\zeta_1$ is a sidekick for $1$.
We claim that for $1\le \sigma<\nu$, if $\zeta_{\sigma}$ is a sidekick for $\sigma$ then 
$\zeta_{\sigma+1}$ is a sidekick for $\sigma+1$. For let
$M$ be a mat for a stable shower $\mathcal{S}$ in a candidate $G'$, 
such that $\mathcal{S}$ is not $(\sigma+1)$-complete over $M$. By \ref{knoblemma4},
$\chi(M)\le \zeta_{\sigma+1}$.
This proves the claim that $\zeta_{\sigma+1}$ is a sidekick for $\sigma+1$. Consequently 
$\zeta_{\nu}$ is a sidekick for $\nu$,
and in particular, for every candidate $G$, every $\nu$-incomplete stable shower in $G$ has floor of 
chromatic number at most $\zeta_{\nu}$.
By \ref{getstable}, every candidate has chromatic number at most $44\nu (\kappa+\zeta_{\nu})^{(\nu+1)^2} + 4\kappa$. 
This proves \ref{localbound2}.~\bbox

\section{Acknowledgement}

We would particularly like to thank one of the referees, who gave the paper an immensely thorough checking, resulting
in numerous corrections and improvements.
Thanks also to Maria Chudnovsky, who worked with us on parts of the proof.

\end{document}